\documentclass[11pt,a4paper]{amsart}
\usepackage[utf8]{inputenc}
\usepackage[T1]{fontenc}
\usepackage{amsmath,amssymb,amsthm,mathtools}
\usepackage[margin=2.6cm]{geometry}
\usepackage{enumitem}
\usepackage{booktabs}
\usepackage{hyperref}

\theoremstyle{plain}
\newtheorem{theorem}{Theorem}[section]
\newtheorem{proposition}[theorem]{Proposition}
\newtheorem{lemma}[theorem]{Lemma}
\newtheorem{corollary}[theorem]{Corollary}
\theoremstyle{definition}
\newtheorem{definition}[theorem]{Definition}
\newtheorem{example}[theorem]{Example}
\theoremstyle{remark}
\newtheorem{remark}[theorem]{Remark}

\DeclareMathOperator{\Tr}{Tr}
\DeclareMathOperator{\Div}{div}
\DeclareMathOperator{\Ric}{Ric}
\DeclareMathOperator{\Sch}{Sch}
\DeclareMathOperator{\Rm}{Rm}
\DeclareMathOperator{\Hess}{Hess}
\DeclareMathOperator{\Id}{Id}
\DeclareMathOperator{\Ima}{Im}
\newcommand{\Lr}{\mathring{L}}
\newcommand{\Sr}{\mathring{S}_2}
\newcommand{\TT}{\mathrm{TT}}
\newcommand{\cW}{\mathcal{W}}
\newcommand{\KN}{\mathbin{\bigcirc\mspace{-15mu}\wedge\mspace{3mu}}}
\newcommand{\ip}[2]{\langle #1,#2\rangle}
\newcommand{\Ricr}{\mathring{\Ric}}
\newcommand{\Rev}{\mathcal G}

\title[Natural potentials for TT tensors and linearised initial data]{Natural potentials for TT tensors, divergence-free symmetric tensors and linearised relativistic initial data}
\author{Erwann Delay}

\address{Avignon Universit\'e, Laboratoire de Math\'ematiques d'Avignon (EA 2151), F-84018 Avignon, France}
\address{erwann.delay@univ-avignon.fr}
\address{https://erwanndelay.wordpress.com}
\date{November 6, 2026}

\begin{document}

\begin{abstract}
Let $N$ be a natural tensor which vanishes at a metric $g_0$. By diffeomorphism invariance, $DN(g_0)^*$ takes values in divergence-free symmetric $2$-tensors, and in TT tensors if moreover $DN(g_0)$ kills conformal directions. Such operators are local potentials, which we call machines: they produce smooth compactly supported tensors with the required property and, on closed manifolds, all of them up to a finite-dimensional space when their symbol is injective. On Einstein manifolds of every dimension, the linearised Bach tensor gives the operator $P$ of an earlier paper, which factorises on TT tensors into the massless and partially massless spin-$2$ operators; we also construct machines on Bach-flat, conformally flat and Ricci-parallel metrics. Replacing the diffeomorphisms of $M$ by those of the vacuum development, we obtain potentials, with a Fredholm property, for the linearised Einstein constraint equations on arbitrary spacelike hypersurfaces of Minkowski, de Sitter and anti-de Sitter spacetimes.
\end{abstract}

\maketitle

\tableofcontents

\section{Introduction}\label{sec:intro}

A symmetric $2$-tensor $h$ on a Riemannian manifold $(M,g)$ is called TT, for \emph{transverse traceless}, if $\Div_gh=0$ and $\Tr_gh=0$. These tensors are central in general relativity. In linearised gravity, the TT gauge isolates the physical degrees of freedom of gravitational waves, and in the conformal method of Lichnerowicz and York the second fundamental form of an initial data set is split into a TT part, which belongs to the free data, a part $\Lr W$ and a trace part \cite{Y}; on conformally flat slices, York's conformal covariance of TT tensors allows one to transport flat TT tensors, and Beig and Ó~Murchadha studied the momentum constraint in this setting \cite{B,BOM}. They are just as important in Riemannian geometry: by the Berger--Ebin decomposition \cite{BE}, the infinitesimal Einstein deformations of an Einstein metric, modulo diffeomorphisms and homotheties, are TT tensors \cite[Ch.~12]{Be}, and the stability of Einstein metrics, for the total scalar curvature as well as for the Ricci flow, is read off the Lichnerowicz Laplacian on TT tensors. This is how the instability of the Riemannian Schwarzschild metric appears \cite{GPY}; see also \cite{Ko,Kr}.

The single condition $\Div h=0$ is even more widespread. In elasticity, in the absence of body forces, the Cauchy stress tensor is a symmetric tensor satisfying the equilibrium equation $\Div\sigma=0$. In $\mathbb R^3$, equilibrated stresses are represented by Beltrami stress functions, $\sigma=\operatorname{inc}\Phi$ \cite{Gu}, which expresses the exactness of a piece of the elasticity complex; this complex and its curved version due to Calabi \cite{Ca,Kh} also serve as a model for mixed finite elements \cite{AW,Ea}, and stress potentials on manifolds of constant curvature and pre-complexes for general metrics are studied in \cite{KL1,KL2}. In general relativity, the energy-momentum tensor is divergence free, and in vacuum the momentum constraint reads $\Div(K-(\Tr K)\,g)=0$. Positive divergence-free symmetric tensors also enter compensated compactness estimates in fluid mechanics \cite{Se}.

The constraint equations themselves are a third source of underdetermined linear problems. A vacuum initial data set $(g,K)$ satisfies $\Phi(g,K)=0$, and the solutions of $D\Phi(g_0,K_0)(h,k)=0$ are the tangent vectors to the space of solutions. Compactly supported solutions of these linearised equations are the linear ingredient of gluing constructions \cite{Co,CSc,CD03,DC}; they are usually obtained by solving an underdetermined elliptic problem in weighted spaces, or through solution operators with prescribed support \cite{IMOT,MOT}, and explicit potentials are known in dimension $3$ on space forms \cite{BC20}.

\subsection*{Machines}
What these three situations have in common is that the relevant condition is destroyed by cut-off: $\Div(\chi h)=\chi\Div h-h(\nabla\chi,\cdot)$. We call \emph{TT-tensor machine} (resp.\ \emph{divergence-free tensor machine}) a linear differential operator $\mathcal P$, defined on the sections of a bundle $F$, whose image consists of TT tensors (resp.\ divergence-free symmetric tensors) for $g_0$, and similarly for the linearised constraints. A machine moves the condition to the level of the potential: if $k=\mathcal P\Psi$ and $\chi$ is a cut-off function, then $\mathcal P(\chi\Psi)$ still has the required property, is compactly supported, and equals $k$ where $\chi=1$. In particular, a machine which does not vanish identically on an open set $U$ produces an infinite-dimensional space of smooth tensors with the required property and support in $U$. It also gives smooth compactly supported approximations which preserve the property exactly, for every tensor in its image (Proposition~\ref{prop:approx}), hence, under the conclusion of Theorem~\ref{thm:fredholm}, for every TT tensor up to a finite-dimensional space. Finally, it provides explicit compactly supported test directions, for instance to exhibit unstable directions of an Einstein metric, or for gluing and localised deformation problems in the spirit of \cite{Co,DC}.

All the machines of this paper come from a single principle. Let $N$ be a natural tensor which vanishes at $g_0$. Since $N(\varphi^*g_0)=\varphi^*N(g_0)=0$ for every diffeomorphism $\varphi$, the linearisation $DN(g_0)$ kills the Lie derivatives $L_Xg_0$, so that its adjoint takes values in $\ker\Div$; if $DN(g_0)$ also kills the conformal directions $fg_0$, the adjoint takes values in TT tensors (Proposition~\ref{prop:principe}). This observation is essentially folklore: it is the adjoint form of the linearised Bianchi identity, it is close in spirit to Wald's adjoint method, in which an operator identity is used to produce solutions of perturbation equations from adjoint operators \cite{Wa}, and in dimension $3$ it underlies Beig's construction of TT tensors from the Cotton tensor \cite{B}. Conversely, TT machines are exactly the adjoints of the operators $\mathfrak A$ with $\mathfrak A\circ L=0$ and $\mathfrak A(fg_0)=0$ for every function $f$, that is, of the compatibility operators of $\Lr$, in the sense of Spencer \cite{Sp}, which vanish on pure-trace tensors. The real difficulty is to find, for a given metric, a natural tensor which vanishes there and has good symbolic properties. Infinitesimal naturality turns out to be enough, which allows background structures such as the eigenspaces of a parallel Ricci tensor (Proposition~\ref{prop:relache}). For the linearised constraints, the diffeomorphisms of $M$ are replaced by those of the vacuum spacetime generated by the data, and the Killing operator by the Killing initial data operator. Since only naturality and integration by parts are used, the identities also hold in Lorentzian signature (Remark~\ref{rem:lorentz}).

\subsection*{Main results}
We first show that a machine $DN(g_0)^*$ produces all TT tensors up to a finite-dimensional space on a closed manifold as soon as the symbol of $DN(g_0)$ is injective on TT symbols (Theorem~\ref{thm:fredholm}). Our main example concerns Einstein manifolds of dimension $n\ge3$: there, the operator $P$ of \cite{D} is, up to the factor $-2(n-2)$, the linearised Bach tensor $B=\nabla^kC_{ijk}+A^{kl}W_{ikjl}$, although $B$ is conformally covariant only for $n=4$ (Theorem~\ref{thm:bach}). On TT tensors, $P$ is the product of the massless and the partially massless spin-$2$ operators (Proposition~\ref{prop:PM}), a fact which we compare with Matsumoto's factorisation of the linearised obstruction tensor \cite{Ma} and with critical gravity \cite{DLLPST,LPP}. Beyond Einstein metrics, the Bach tensor still gives machines on Bach-flat metrics (Theorem~\ref{thm:bachplat}), and the Weyl tensor gives a machine of order $2$ on conformally flat metrics (Theorem~\ref{thm:confplat}). Metrics with parallel Ricci tensor require a different idea, based on partial Ricci operators: we obtain a machine of order $4$ under a sign condition on the Ricci eigenvalues (Theorem~\ref{thm:ricpar}) and one of order $5$ with no condition at all (Theorem~\ref{thm:ordre5}). Products of Einstein metrics and cones are treated in Section~\ref{sec:produits}.

Dropping the trace condition makes the conformal requirement disappear, and divergence-free machines are much easier to come by. In constant curvature one recovers the Beltrami stress functions of elasticity, on Einstein metrics the linearised Einstein operator is itself a machine, Ricci-parallel metrics are handled by the partial Ricci operators, and locally symmetric metrics by $\nabla\Rm$. The same ideas apply to the weighted divergence associated with a measure $e^{-f}dV$; on gradient Ricci solitons, we only give a local construction (Theorem~\ref{thm:soliton}).

The last part of the paper deals with the linearised constraint equations. Here a natural tensor of the data which vanishes on all spacelike hypersurfaces of the vacuum development yields a machine (Proposition~\ref{prop:principegen}); the proof relies on the Hamiltonian form of the vacuum evolution equations. Taking the Gauss and Codazzi tensors, we obtain machines with a Fredholm property on arbitrary spacelike hypersurfaces of Minkowski, de Sitter and anti-de Sitter spacetimes (Theorem~\ref{thm:contrgen}). For time-symmetric data the argument becomes elementary and also applies in Riemannian signature; with $\bar\nabla\bar\Rm$, it gives machines at locally symmetric static data such as $S^2\times S^1$, whose static extension is the Nariai solution (Theorem~\ref{thm:contrlocsym}).

\subsection*{Organisation}
The three parts of the paper rest on the same two ingredients: a natural tensor vanishing on the geometry under consideration, whose linearisation is killed by the infinitesimal symmetries (Lie derivatives, conformal directions, or the deformations generated by the vacuum evolution), and the Fredholm Lemma~\ref{lem:fredholmD}, which only requires a symbol computation. After fixing conventions in Section~\ref{sec:conv}, we present the principle, the Fredholm lemma and the compactly supported approximations in Section~\ref{sec:principe}. Sections~\ref{sec:bach} to~\ref{sec:produits} are devoted to TT machines, Sections~\ref{sec:divnulle} and~\ref{sec:divf} to divergence-free and weighted divergence-free machines, and Section~\ref{sec:contraintes} to the linearised constraints. Section~\ref{sec:contraintes} only uses Section~\ref{sec:principe} and the divergence-free machines of Propositions~\ref{prop:divcc} and~\ref{prop:divlocsym}, so that a reader interested in the constraint equations may go directly from Section~\ref{sec:principe} to Section~\ref{sec:contraintes}. The Hamiltonian identity used in Section~\ref{ss:contrgen} is proved in Appendix~\ref{app:hamilton}.

\section{Conventions}\label{sec:conv}

$(M,g)$ is a smooth Riemannian manifold of dimension $n\ge 3$. We denote by $S_2$ (resp.\ $\Sr$) the bundle of symmetric (resp.\ trace-free symmetric) $2$-tensors. We follow the conventions of \cite{D}:
\[
(\Div h)_i=-\nabla^j h_{ij},\qquad (LX)_{ij}=\nabla_iX_j+\nabla_jX_i,\qquad L^*=2\Div,
\]
\[
\Lr=L+\tfrac{2}{n}d^*(\cdot)\,g,\qquad \Lr^*=2\Div+\tfrac2n d\Tr,\qquad \TT=\ker\Div\cap\ker\Tr,
\]
$\nabla^*\nabla=-\nabla^k\nabla_k$, and
$\Delta_Lh_{ij}=\nabla^*\nabla h_{ij}+R_{ik}h^k{}_j+R_{jk}h^k{}_i-2R_{ikjl}h^{kl}$.
The Schouten tensor and its normalised version are
\[
\Sch=\Ric-\frac{R}{2(n-1)}\,g,\qquad A=\frac{1}{n-2}\,\Sch .
\]
The Lichnerowicz Laplacian was introduced in \cite{Li}; the classical formulas used without proof (Ricci and Bianchi identities, decomposition of the curvature tensor, first variation formulas for the curvature, conformal behaviour of $W$ and $\Sch$) can be found for instance in \cite[Ch.~1]{Be} and \cite[Ch.~2]{To}.

\subsection*{Curvature}
We define $R^a{}_{bcd}$ by $[\nabla_c,\nabla_d]X^a=R^a{}_{bcd}X^b$, then $R_{abcd}=g_{am}R^m{}_{bcd}$ and $\Ric_{bd}=g^{ac}R_{abcd}=R^c{}_{bcd}$. For a $1$-form, $[\nabla_c,\nabla_d]\omega_a=-R^m{}_{acd}\,\omega_m$. On a space form of curvature $\kappa$, $R_{abcd}=\kappa(g_{ac}g_{bd}-g_{ad}g_{bc})$.

The Kulkarni--Nomizu product of two symmetric $2$-tensors is
\[
(U\KN V)_{abcd}=U_{ac}V_{bd}+U_{bd}V_{ac}-U_{ad}V_{bc}-U_{bc}V_{ad}.
\]
One has $\Rm=W+A\KN g$. We denote by $\cW$ the bundle of algebraic Weyl tensors (curvature symmetries, algebraic Bianchi identity, all traces zero). For $Z\in\cW$ and $U\in S_2$ we set $(Z\circ U)_{bd}=Z_{abcd}U^{ac}$, a symmetric tensor by the pair symmetry.

For an algebraic curvature tensor $T$, we denote by $\mathrm c(T)_{bd}=g^{ac}T_{abcd}$ its Ricci contraction and by $\Tr\mathrm c(T)$ the trace of the latter. The Weyl projection
\[
\cW(T)=T-\tfrac{1}{n-2}\Big(\mathrm c(T)-\tfrac{\Tr\mathrm c(T)}{2(n-1)}g\Big)\KN g
\]
gives the orthogonal decomposition \cite[Ch.~1]{Be}
\begin{equation}\label{eq:normes}
|T|^2=|\cW(T)|^2+\tfrac{4}{n-2}\Big|\mathrm c(T)-\tfrac{\Tr\mathrm c(T)}{n}g\Big|^2+\tfrac{2}{n(n-1)}\big(\Tr\mathrm c(T)\big)^2 .
\end{equation}

\subsection*{Linearisations, adjoints, symbols}
$D$ denotes the differential with respect to the metric. For $h\in S_2$ we write
\[
\delta\Gamma^m_{jk}=\tfrac12 g^{mp}(\nabla_jh_{pk}+\nabla_kh_{pj}-\nabla_ph_{jk}).
\]
For a covariant tensor $T$ depending on $g$:
\begin{equation}\label{eq:DnablaT}
D(\nabla T)(h)=\nabla\big(DT(h)\big)-\delta\Gamma\cdot T,\qquad (\delta\Gamma\cdot T)_{e\,i_1\dots i_r}=\sum_{s}\delta\Gamma^m_{e i_s}T_{i_1\dots m\dots i_r}.
\end{equation}
Since $\nabla g=0$ for every metric, $\delta\Gamma\cdot g=\nabla h$, and $\delta\Gamma\cdot$ acts as a derivation on tensor products.

Adjoints are formal $L^2(dV_{g_0})$ adjoints, with full contractions. The principal symbol $\sigma_\xi$ is obtained by replacing $\nabla$ with $i\xi$ in the top-order terms; thus $\sigma_\xi(\mathcal D^*)=\sigma_\xi(\mathcal D)^*$, and the symbol of a composition is the composition of the symbols. We denote by
\[
\TT_\xi=\{k\in\Sr:\ k(\xi^\sharp,\cdot)=0\}
\]
the space of TT symbols at $\xi\neq0$.

We shall constantly use the first variation of the curvature tensor, which we recall with our conventions.

\begin{lemma}[{see \cite[1.174]{Be}, \cite[Ch.~2]{To}}]\label{lem:DRm}
For every metric $g$ and every $h\in S_2$:
\[
DR(h)_{abcd}=h_{am}R^m{}_{bcd}+\tfrac12[\nabla_c,\nabla_d]h_{ab}
+\tfrac12\big(\nabla_c\nabla_b h_{ad}+\nabla_d\nabla_a h_{bc}-\nabla_c\nabla_a h_{bd}-\nabla_d\nabla_b h_{ac}\big).
\]
In particular $\sigma_\xi(DR)h=\tfrac12\,\xi\xi\KN h$, independently of the metric.
\end{lemma}

This is the classical formula $\delta R^m{}_{bcd}=\nabla_c\delta\Gamma^m_{db}-\nabla_d\delta\Gamma^m_{cb}$, with the index lowered.

\section{The principle}\label{sec:principe}

Machines will be produced by linearising operators which commute with diffeomorphisms; we first make this notion precise.

\begin{definition}
Let $\mathcal U$ be a diffeomorphism-invariant open set of metrics and $F$ a tensor bundle. An operator $N:\mathcal U\to\Gamma(F)$ is \emph{natural} if $N(\varphi^*g)=\varphi^*N(g)$ for every (local) diffeomorphism $\varphi$. The coefficients of $N$ may involve fixed constants.
\end{definition}

The following elementary observation is the source of all the machines of this paper.

\begin{proposition}[Principle]\label{prop:principe}
Let $N$ be natural and $g_0\in\mathcal U$ be such that
\begin{enumerate}[label=(\roman*)]
\item $N(g_0)=0$;
\item $DN(g_0)(fg_0)=0$ for every function $f$.
\end{enumerate}
Then $\mathcal P_N:=DN(g_0)^*:\Gamma(F)\to\Gamma(S_2)$ satisfies
\[
\Div\mathcal P_N=0,\qquad \Tr\mathcal P_N=0,\qquad DN(g_0)\circ\Lr=0.
\]
\end{proposition}

\begin{proof}
Let $X$ be compactly supported, with flow $\varphi_t$. By naturality, $N(\varphi_t^*g_0)=\varphi_t^*N(g_0)=0$; differentiating at $t=0$ gives $DN(g_0)(LX)=0$. For $\Psi\in\Gamma_c(F)$ we deduce
\[
0=\ip{\Psi}{DN(g_0)LX}_{L^2}=\ip{\mathcal P_N\Psi}{LX}_{L^2}=2\ip{\Div\mathcal P_N\Psi}{X}_{L^2},
\]
hence $\Div\mathcal P_N=0$. Similarly, (ii) gives $\int f\,\Tr(\mathcal P_N\Psi)\,dV=\ip{\Psi}{DN(g_0)(fg_0)}=0$ for every $f$, hence $\Tr\mathcal P_N=0$. Finally $\Lr w=Lw+\frac2nd^*w\,g_0$ is killed by $DN(g_0)$.
\end{proof}

Several identities below are first established on TT tensors and then extended to all trace-free tensors. This relies on a local decomposition, which we prove now.

\begin{lemma}[Local TT decomposition]\label{lem:localTT}
Let $(M,g)$ be a Riemannian manifold, $x\in M$. There is $r_0>0$ such that, on every geodesic ball $U=B(x,r)$ with $r\le r_0$, every smooth $\mathring h\in\Gamma(\Sr)$ can be written $\mathring h=k+\Lr w$ on $U$, with $k$ smooth and TT on $U$ and $w$ smooth on $U$.
\end{lemma}

\begin{proof}
For $w\in C^\infty_c(U)$, using $\int\nabla_iw_j\nabla^jw^i=\int(\nabla^iw_i)^2-\int\Ric(w,w)$ (integration by parts and the Ricci identity) and $|\Lr w|^2=|Lw|^2-\frac4n(\nabla^iw_i)^2$, we get
\[
\int_U|\Lr w|^2dV=2\int_U|\nabla w|^2dV+\Big(2-\frac4n\Big)\int_U(\nabla^iw_i)^2dV-2\int_U\Ric(w,w)\,dV .
\]
By the Poincaré inequality $\int_U|w|^2\le Cr^2\int_U|\nabla w|^2$ on small balls, the quadratic form $w\mapsto\int_U|\Lr w|^2$ is coercive on $H^1_0(U)$ for $r$ small. By the Lax--Milgram theorem there is $w\in H^1_0(U)$ with $\int_U\ip{\Lr w}{\Lr v}=\int_U\ip{\mathring h}{\Lr v}$ for every $v\in H^1_0(U)$, i.e.\ $\Lr^*\Lr w=\Lr^*\mathring h$ weakly. Since $\Lr^*\Lr$ is elliptic with smooth coefficients and $\Lr^*\mathring h$ is smooth, $w$ is smooth in $U$ by interior elliptic regularity, and $k:=\mathring h-\Lr w$ is smooth, trace free, and satisfies $2\Div k=\Lr^*k=0$ on $U$.
\end{proof}

The Fredholm statements of this paper all follow from one lemma, in the spirit of the decomposition theorem of Berger and Ebin \cite{BE}. Its setting is the following: an operator $D$ defines the equation we are interested in (a divergence, or the linearised constraints), and a machine $\mathcal P$ produces solutions of $Du=0$. The lemma says when $\mathcal P$ produces all of them up to a finite-dimensional space.

\begin{lemma}\label{lem:fredholmD}
Let $(M,g_0)$ be closed, $f\in C^\infty(M)$ and $\mu=e^{-f}dV_{g_0}$. Let $E$, $F$, $G$ be tensor bundles over $M$, endowed with the fibre metrics and the Levi-Civita connection $\nabla$ of $g_0$; adjoints are taken in $L^2(\mu)$. Let
\[
\mathcal P:\Gamma(F)\to\Gamma(E),\qquad D:\Gamma(E)\to\Gamma(G)
\]
be differential operators with $D\circ\mathcal P=0$. Assume that $G=G_1\oplus\dots\oplus G_p$ and $F=F_1\oplus\dots\oplus F_q$, that the components $D_\beta:\Gamma(E)\to\Gamma(G_\beta)$ of $D$ have orders $d_\beta\ge1$, and that the restrictions $\mathcal P_\alpha:\Gamma(F_\alpha)\to\Gamma(E)$ of $\mathcal P$ have orders $m_\alpha\ge1$ (in most applications $p=q=1$); the adjoints $\mathcal P_\alpha^*:\Gamma(E)\to\Gamma(F_\alpha)$ are the components of $\mathcal P^*=(\mathcal P_1^*,\dots,\mathcal P_q^*)$. Assume finally that, for every $x\in M$ and $\xi\in T^*_xM\setminus\{0\}$, the only $v\in E_x$ such that $\sigma_\xi(D_\beta)v=0$ for all $\beta$ and $\sigma_\xi(\mathcal P_\alpha^*)v=0$ for all $\alpha$ is $v=0$. Then, denoting by $\ker D$ the space of smooth sections $u$ of $E$ with $Du=0$,
\[
\ker D=\Ima\mathcal P\oplus\big(\ker D\cap\ker\mathcal P^*\big),
\]
the sum being $L^2(\mu)$-orthogonal and the second summand finite dimensional.
\end{lemma}

When $p=q=1$, the assumption says that $\sigma_\xi(\mathcal P^*)$ is injective on $\ker\sigma_\xi(D)$. For instance, for a TT machine $\mathcal P=DN(g_0)^*$, one takes $E=\Sr$ and $D=\Div$, and the assumption is the injectivity of $\sigma_\xi(DN(g_0))$ on $\TT_\xi$.

\begin{proof}
\emph{Step 1: an auxiliary operator.} Let $\nabla^*$ be the $L^2(\mu)$-adjoint of $\nabla$, $\Xi=1+\nabla^*\nabla$, which acts on each tensor bundle, and $c=\max(d_\beta,m_\alpha)$. We set
\[
\mathcal Y=\sum_\alpha\mathcal P_\alpha\Xi^{c-m_\alpha}\mathcal P_\alpha^*+\sum_\beta D_\beta^*\Xi^{c-d_\beta}D_\beta .
\]
Each term is a differential operator of order at most $2c$ on $\Gamma(E)$, so $\mathcal Y$ is of order $2c$, and $\mathcal Y$ is self-adjoint in $L^2(\mu)$, since $\Xi$ is.

\emph{Step 2: ellipticity.} The principal symbol of $\Xi$ is $|\xi|^2$, that of an adjoint is the adjoint of the symbol, and the symbol of a composition is the composition of the symbols. Hence, for $v\in E_x$ and $\xi\in T^*_xM$,
\[
\ip{\sigma_\xi(\mathcal Y)v}{v}=\sum_\alpha|\xi|^{2(c-m_\alpha)}\big|\sigma_\xi(\mathcal P_\alpha^*)v\big|^2+\sum_\beta|\xi|^{2(c-d_\beta)}\big|\sigma_\xi(D_\beta)v\big|^2 .
\]
For $\xi\neq0$, this vanishes only if $v$ lies in all the kernels $\ker\sigma_\xi(\mathcal P_\alpha^*)$ and $\ker\sigma_\xi(D_\beta)$, that is, by assumption, only if $v=0$. Thus $\sigma_\xi(\mathcal Y)$ is positive definite, and $\mathcal Y$ is elliptic. By the standard theory of elliptic operators on closed manifolds \cite{BE}, \cite[Appendix]{Be}, $\ker \mathcal Y$ is a finite-dimensional space of smooth sections, and every smooth section which is $L^2(\mu)$-orthogonal to $\ker \mathcal Y$ is of the form $\mathcal Yu$ with $u$ smooth.

\emph{Step 3: the kernel of $\mathcal Y$.} For a section $v$, $\ip{\Xi v}{v}=\|v\|^2+\|\nabla v\|^2\ge\|v\|^2$, so $\|\Xi v\|\ge\|v\|$ by the Cauchy--Schwarz inequality, and by induction $\ip{\Xi^jv}{v}\ge\|v\|^2$ for every $j\ge0$: for $j=2i$ it equals $\|\Xi^iv\|^2$, and for $j=2i+1$ it equals $\ip{\Xi\Xi^iv}{\Xi^iv}\ge\|\Xi^iv\|^2$. Therefore
\[
\ip{\mathcal Yu}{u}=\sum_\alpha\ip{\Xi^{c-m_\alpha}\mathcal P_\alpha^*u}{\mathcal P_\alpha^*u}+\sum_\beta\ip{\Xi^{c-d_\beta}D_\beta u}{D_\beta u}\ \ge\ \|\mathcal P^*u\|^2+\|Du\|^2 .
\]
If $\mathcal Yu=0$, then $\mathcal P^*u=0$ and $Du=0$; conversely $\mathcal P^*u=Du=0$ obviously implies $\mathcal Yu=0$. Hence $\ker\mathcal Y=\ker\mathcal P^*\cap\ker D$, which is finite dimensional by Step~2.

\emph{Step 4: the image.} Since $D\mathcal P=0$, $\Ima\mathcal P\subset\ker D$; and $\Ima\mathcal P$ is orthogonal to $\ker\mathcal P^*$, since $\ip{\mathcal Pw}{v}=\ip{w}{\mathcal P^*v}$. In particular $\Ima\mathcal P\cap(\ker D\cap\ker\mathcal P^*)=\{0\}$. Let now $k\in\ker D$ be smooth. Subtracting from $k$ its orthogonal projection onto the finite-dimensional space $\ker D\cap\ker\mathcal P^*=\ker \mathcal Y$, which lies in $\ker D$, we may assume that $k$ is orthogonal to $\ker \mathcal Y$. By Step~2, $k=\mathcal Yu$ with $u$ smooth, that is
\[
k=\mathcal Pw+D^*z,\qquad w=\big(\Xi^{c-m_\alpha}\mathcal P_\alpha^*u\big)_\alpha,\quad z=\big(\Xi^{c-d_\beta}D_\beta u\big)_\beta .
\]
Applying $D$ and using $Dk=0$ and $D\mathcal P=0$, we get $DD^*z=0$, hence $\|D^*z\|^2=\ip{DD^*z}{z}=0$. Thus $k=\mathcal Pw\in\Ima\mathcal P$, which proves the decomposition.
\end{proof}

When the symbol of $DN(g_0)$ is injective on TT symbols, the image of the machine is as large as it can be: it contains all TT tensors up to a finite-dimensional space.

\begin{theorem}[Fredholm property]\label{thm:fredholm}
Assume $M$ is compact without boundary and $N$ is as in Proposition~\ref{prop:principe}, with $DN(g_0)$ of order $m\ge1$. If $\sigma_\xi(DN(g_0))$ is injective on $\TT_\xi$ for every $\xi\neq0$, then, in $C^\infty$,
\[
\TT=\Ima\mathcal P_N\ \oplus\ \big(\ker DN(g_0)\cap\TT\big),
\]
the sum being $L^2$-orthogonal and the second summand finite dimensional.
\end{theorem}

\begin{proof}
We apply Lemma~\ref{lem:fredholmD} with $f=0$, $E=\Sr$, $G=T^*M$, $D=\Lr^*=2\Div:\Gamma(\Sr)\to\Gamma(T^*M)$, of order $1$, and $\mathcal P=\mathcal P_N:\Gamma(F)\to\Gamma(\Sr)$, of order $m$; $\mathcal P_N$ takes values in $\Sr$ by Proposition~\ref{prop:principe}. On $\Sr$, $\ker D=\TT$, and $D\circ\mathcal P_N=(DN(g_0)\circ\Lr)^*=0$ by Proposition~\ref{prop:principe}. Since $\mathcal P_N^*$ is the restriction of $DN(g_0)$ to $\Sr$ and $\ker\sigma_\xi(\Lr^*)=\TT_\xi$, the assumption of the lemma is exactly the injectivity of $\sigma_\xi(DN(g_0))$ on $\TT_\xi$.
\end{proof}

\subsection{Compactly supported approximations}

We now make precise how machines produce compactly supported tensors, and in which sense they approximate given ones.

\begin{proposition}\label{prop:approx}
Let $\mathcal P:\Gamma(F)\to\Gamma(S_2)$ be a differential operator of order $m$ with smooth coefficients, with values in TT (resp.\ divergence-free) tensors for $g_0$.
\begin{enumerate}[label=(\alph*)]
\item For every open set $U$, $\mathcal P\big(C^\infty_c(U;F)\big)$ consists of smooth TT (resp.\ divergence-free) tensors supported in $U$.
\item If $k=\mathcal P\Psi$ and $\chi\in C^\infty_c(M)$, then $\mathcal P(\chi\Psi)$ is TT (resp.\ divergence free), supported in $\operatorname{supp}\chi$, and
\[
\mathcal P(\chi\Psi)=\chi\,k+[\mathcal P,\chi]\Psi ,
\]
where the commutator $[\mathcal P,\chi]$ has order $m-1$ and only involves derivatives of $\chi$. In particular $\mathcal P(\chi\Psi)=k$ where $\chi=1$ near the point.
\item Assume $(M,g_0)$ is complete with bounded geometry and the coefficients of $\mathcal P$ are bounded together with their derivatives. If $s\ge0$ is an integer and $k=\mathcal P\Psi$ with $\Psi\in H^{s+m}(F)$, then $k$ is an $H^s$-limit of smooth compactly supported TT (resp.\ divergence-free) tensors.
\end{enumerate}
\end{proposition}

\begin{proof}
(a) and (b) are immediate, $\mathcal P$ being differential with values in TT (resp.\ divergence-free) tensors. (c) On a complete manifold with bounded geometry, $C^\infty_c(F)$ is dense in $H^{s+m}(F)$, see \cite{Au,He} for functions and \cite{Sh} for vector bundles; if $\Psi_j\to\Psi$ in $H^{s+m}$ with $\Psi_j\in C^\infty_c$, then $\mathcal P\Psi_j\to k$ in $H^s$ by continuity of $\mathcal P:H^{s+m}\to H^s$.
\end{proof}

\begin{remark}
Point (c) reduces the approximation of a TT tensor by compactly supported TT tensors to the question whether it lies in the image of the machine with a controlled potential. On a compact manifold, under the assumption of Theorem~\ref{thm:fredholm}, every TT tensor in $H^s$ can be written $\mathcal P\Psi+k_0$ with $k_0$ in a finite-dimensional space of smooth tensors and $\Psi\in H^{s+m}$, by elliptic regularity of the auxiliary operator $\mathcal Y$ of the proof of Lemma~\ref{lem:fredholmD}. On a non-compact manifold, one must also control $\Psi$ at infinity, which requires a weighted analysis adapted to each geometry. An application along these lines is given in a companion paper \cite{D3}, where compactly supported TT unstable directions for the Riemannian Schwarzschild--Tangherlini metrics are obtained with the operator $P$ of \cite{D}.
\end{remark}

\begin{remark}[Weighted spaces on $\mathbb R^n$]\label{rem:poidsRn}
On $\mathbb R^n$, $n\ge3$, consider the weighted Sobolev spaces $W^{k,p}_\delta$ of Bartnik \cite{Bart}, whose elements behave like $o(r^\delta)$. The operator $P$ of \cite{D} equals $(\nabla^*\nabla)^2$ on TT tensors, and the mapping properties of $\Delta^2$ and of the vector Laplacian $\Div\Lr$ in these spaces \cite{Bart} suggest that $P$ maps $W^{k+4,p}_\delta(\Sr)$ onto $\TT\cap W^{k,p}_{\delta-4}$ when $4-n<\delta<3$ and $\delta\notin\mathbb Z$, with an image of finite codimension for every $\delta\notin\mathbb Z$. We do not prove this statement here, and we do not use it. The range $4-n<\delta<3$ would cover the TT tensors $k=O(r^{\delta-4})$ with $-n<\delta-4<-1$, in particular the decay rates $r^{-\tau-1}$, $\tau>\frac{n-2}{2}$, of the second fundamental form for which the ADM mass is defined; since $C^\infty_c$ is dense in $W^{k+4,p}_\delta$, every such TT tensor would then be a limit in $W^{k,p}_{\delta-4}$ of smooth compactly supported TT tensors.
\end{remark}

Machines can be moved within a conformal class, with the weight dictated by York's conformal covariance of TT tensors.

\begin{lemma}[Conformal transport]\label{lem:transport}
If $g_0=e^f\tilde g$ and $\tilde{\mathcal P}$ is a TT machine for $\tilde g$, then $\Psi\mapsto e^{-\frac{n-2}{2}f}\tilde{\mathcal P}\Psi$ is a TT machine for $g_0$.
\end{lemma}

\begin{proof}
By, e.g., \cite[Lemme 7]{D}, $\Div_{g_0}h=e^{-\frac n2f}\Div_{\tilde g}(e^{\frac{n-2}{2}f}h)+\frac12e^{-f}(\Tr_{\tilde g}h)\,df$. For $h=e^{-\frac{n-2}{2}f}k$ with $k$ TT for $\tilde g$ we get $\Tr_{g_0}h=0$ and $\Div_{g_0}h=0$.
\end{proof}

\subsection{Infinitesimal naturality}

Naturality of $N$ enters Proposition~\ref{prop:principe} only through the identity $DN(g_0)(LX)=L_X N(g_0)$. Infinitesimal naturality at $g_0$ is therefore enough, which allows the use of background structures.

\begin{proposition}[Relaxed principle]\label{prop:relache}
Let $\tau$ be a fixed background tensor, $F(g;\tau)$ an operator natural in $(g,\tau)$ jointly (that is, $F(\varphi^*g;\varphi^*\tau)=\varphi^*F(g;\tau)$), and $N(g):=F(g;\tau)$. Then
\[
DN(g_0)(LX)=L_X\big(N(g_0)\big)-D_\tau F(g_0;\tau)\big(L_X\tau\big).
\]
If $N(g_0)=0$, $D_\tau F(g_0;\tau)(L_X\tau)=0$ for every $X$, and $DN(g_0)(fg_0)=0$ for every $f$, then $DN(g_0)^*$ is TT-valued, and the conclusions of Proposition~\ref{prop:principe} and Theorem~\ref{thm:fredholm} remain valid.
\end{proposition}

\begin{proof}
Differentiating $F(\varphi_t^*g_0;\varphi_t^*\tau)=\varphi_t^*F(g_0;\tau)$ at $t=0$ gives
\[
DN(g_0)(LX)+D_\tau F(g_0;\tau)(L_X\tau)=L_XN(g_0).
\]
Under the assumptions, $DN(g_0)\circ L=0$ and $DN(g_0)(fg_0)=0$, and the proof of Proposition~\ref{prop:principe} applies verbatim; the proof of Theorem~\ref{thm:fredholm} only uses these two identities.
\end{proof}

A first application of the relaxed principle is a version of conformal transport in which the conformal factor is an arbitrary fixed function.

\begin{corollary}[Conformal transport without naturality]\label{cor:transportgen}
Let $f$ be an \emph{arbitrary} background function, $G_0:=e^{-f}g_0$ (so that $g_0=e^fG_0$, as in Lemma~\ref{lem:transport}), and $\mathcal N$ a natural operator with $\mathcal N(G_0)=0$ and $D\mathcal N(G_0)(\psi G_0)=0$ for every function $\psi$. Then $N(g):=\mathcal N(e^{-f}g)$ satisfies the assumptions of Proposition~\ref{prop:relache}. In particular $DN(g_0)^*$ is a TT machine for $g_0$, and $DN(g_0)=D\mathcal N(G_0)\circ(e^{-f}\,\cdot)$.
\end{corollary}

\begin{proof}
Here $\tau=f$, and $D_fF(g_0;f)(\delta f)=D\mathcal N(G_0)(-\delta f\,G_0)=0$, in particular for $\delta f=Xf$. Moreover $N(g_0)=\mathcal N(G_0)=0$ and $DN(g_0)(\psi g_0)=D\mathcal N(G_0)(\psi G_0)=0$. The formula for $DN(g_0)$ is immediate.
\end{proof}

Lemma~\ref{lem:transport} is the ``adjoint'' version of this corollary: up to weights $e^{cf}$ depending on the tensor type of the target of $\mathcal N$, the image of $DN(g_0)^*$ is obtained by transporting the image of $D\mathcal N(G_0)^*$ by $k\mapsto e^{-\frac{n-2}{2}f}k$.

\begin{remark}
The important point of the corollary is that $f$ is a fixed function, which does not depend on the metric: the operator $g\mapsto\mathcal N(e^{-f}g)$ is therefore not natural, but its failure of naturality lies in a conformal direction, which is killed by $D\mathcal N(G_0)$. Likewise, the construction of Section~\ref{sec:ricpar} originates in the infinitesimal statement $D\Sch(g_0)(fg_0)=L_{X_f}\Sch_0$ of Lemma~\ref{lem:Xf}, which expresses conformal directions as Lie derivatives of the Schouten tensor.
\end{remark}

Before turning to new constructions, let us list the natural tensors to which the principle applies directly.

\begin{example}\label{ex:liste}
\leavevmode
\begin{enumerate}[label=(\alph*)]
\item \emph{Einstein, $\lambda\neq0$}: $N=\Ricr\circ\Sch$ \cite[\S6]{D}, which gives $P=4\,D\Sch(g)^*D\Ricr(g)^*$, of order $4$.
\item \emph{$n=4$, Bach-flat}: $N=B$. Since $B(e^fg)=e^{-f}B(g)$ \cite[Ch.~4]{Be}, (ii) follows from $B(g_0)=0$; a machine of order $4$. At Einstein metrics, this case extends to every dimension (Section~\ref{sec:bach}).
\item \emph{$n$ even, obstruction-flat}: $N=\mathcal O$, the conformally covariant Fefferman--Graham obstruction tensor \cite{FG,GH}, in the same way; a machine of order $n$. Matsumoto \cite{Ma} showed that, at an obstruction-flat metric, the trace-free part of $D\mathcal O(g_0)$ is a multiple of a conformally covariant GJMS-type operator on trace-free symmetric $2$-tensors, and computed it at Einstein metrics: with $\Ric(g_0)=\lambda g_0$, so that Matsumoto's Schouten tensor is $A(g_0)=\frac{\lambda}{2(n-1)}g_0$, and with $\Delta_L$ normalised as in Section~\ref{sec:conv} (non-negative leading part $\nabla^*\nabla$), it is on TT tensors a multiple of $\prod_{m=0}^{n/2-1}\big(\Delta_L-2\lambda+\frac{2m(n-2m-1)}{n-1}\lambda\big)$.
\item \emph{$n=3$, conformally flat}: $N=$ Cotton tensor; this is Beig's machine of order $3$ \cite{B}.
\end{enumerate}
\end{example}

\begin{remark}[Lorentzian signature]\label{rem:lorentz}
The identities $DN(g_0)\circ\Lr=0$ and Lemma~\ref{lem:transport} only use naturality and integration by parts against compactly supported fields, so they hold verbatim for pseudo-Riemannian metrics; only the symbol and Fredholm statements are lost. For instance, on a vacuum spacetime of dimension $4$ such as Kerr, $DB(g_0)^*\Psi$ is, for $\Psi\in C^\infty_c$, a smooth, compactly supported, trace-free symmetric tensor with vanishing divergence.
\end{remark}

\section{The Bach tensor in dimension $n$ and the operator $P$}\label{sec:bach}

In this section $n\ge3$ and $(M,g_0)$ is Einstein, $\Ric(g_0)=\lambda g_0$. We set
\[
p:=\frac{\lambda}{2(n-1)},\qquad \text{so that}\quad A_0:=A(g_0)=p\,g_0 ,
\]
and write $\mathring W k_{ij}:=W_{ikjl}k^{kl}$, $\mathring R k_{ij}:=R_{ikjl}k^{kl}$.

In arbitrary dimension, we use the following definitions of the Cotton and Bach tensors.

\begin{definition}\label{def:Bach}
The Cotton and Bach tensors of a metric $g$ are
\[
C_{ijk}=\nabla_kA_{ij}-\nabla_jA_{ik},\qquad B_{ij}=\nabla^kC_{ijk}+A^{kl}W_{ikjl}.
\]
In dimension $3$, $W\equiv0$ and $B=\nabla^kC_{ijk}$; in dimension $4$ this is the usual Bach tensor \cite{Ba}, \cite[Ch.~4]{Be}. This definition in arbitrary dimension is the one used, e.g., in \cite{GN,CC}.
\end{definition}

$B$ is a natural tensor of order $4$, symmetric in every dimension $n\ge3$ (the Weyl term is symmetric, and $\nabla^kC_{ijk}-\nabla^kC_{jik}=-\nabla^kC_{kij}$ by the cyclic identity, which vanishes by the Ricci identity, since $\nabla^kA_{ki}=\nabla_i\Tr A$, $\Ric$ commutes with $A$ and $R_{mikj}A^{km}$ is symmetric in $(i,j)$; we shall not need this) and trace free, since $g^{ij}C_{ijk}=0$ by the contracted Bianchi identity and $g^{ij}W_{ikjl}=0$. For $n\neq4$ it is neither conformally covariant, nor divergence free, nor the gradient of a functional. Only the following properties will be used.

\begin{lemma}\label{lem:Bach0}
$C(g_0)=0$, $B(g_0)=0$ and $\nabla^aW(g_0)_{abcd}=0$.
\end{lemma}

\begin{proof}
$A_0=pg_0$ is parallel, so $C(g_0)=0$. Next $A_0^{kl}W_{ikjl}=p\,g^{kl}W_{ikjl}=0$, so $B(g_0)=0$. Finally, the contracted second Bianchi identity gives $\nabla^aR_{abcd}=\nabla_c\Ric_{bd}-\nabla_d\Ric_{bc}=0$, and $\nabla^a(A_0\KN g_0)=0$ since $A_0$ is parallel; hence $\nabla^aW_{abcd}=0$.
\end{proof}

Since the Schouten tensor of $g_0$ is parallel, the linearised Cotton tensor takes a particularly simple form.

\begin{lemma}\label{lem:DnablaA}
For $h\in S_2$, set $\mathcal A(h):=DA(g_0)h-p\,h$. Then
\[
DC(g_0)(h)_{ijk}=\nabla_k\mathcal A_{ij}-\nabla_j\mathcal A_{ik}.
\]
\end{lemma}

\begin{proof}
By \eqref{eq:DnablaT}, $D(\nabla A)(h)=\nabla(DA(h))-\delta\Gamma\cdot A_0=\nabla(DA(h))-p\,\delta\Gamma\cdot g_0=\nabla(DA(h)-p\,h)$.
\end{proof}

The key point is that, at an Einstein metric, the linearised Bach tensor kills conformal directions in every dimension, even though $B$ is not conformally covariant for $n\neq4$.

\begin{proposition}\label{prop:Bconf}
For every function $f$, $DB(g_0)(fg_0)=0$.
\end{proposition}

\begin{proof}
By \cite[Lemme 7]{D}, $D\Sch(g)(fg)=-\frac{n-2}{2}\nabla\nabla f$, so $DA(fg_0)=-\frac12\nabla\nabla f$ and
\[
\mathcal A(fg_0)=-\tfrac12\nabla\nabla f-pf\,g_0 .
\]
Since $\nabla\nabla f$ is symmetric, the Ricci identity for the $1$-form $df$ gives
\[
\nabla_k\nabla_i\nabla_jf-\nabla_j\nabla_i\nabla_kf=[\nabla_k,\nabla_j]\nabla_if=-R^m{}_{ikj}\nabla_mf .
\]
Lemma~\ref{lem:DnablaA} thus yields
\[
DC(fg_0)_{ijk}=\tfrac12R_{mikj}\nabla^mf-p\big(\nabla_kf\,g_{ij}-\nabla_jf\,g_{ik}\big).
\]
Now $R_{mikj}=W_{mikj}+2p(g_{mk}g_{ij}-g_{mj}g_{ik})$, so the terms in $p$ cancel:
\[
DC(fg_0)_{ijk}=\tfrac12W_{mikj}\nabla^mf .
\]
Since $C(g_0)=0$, the variation of $\nabla^kC_{ijk}$ reduces to $\nabla^k DC_{ijk}$ (the terms involving $\delta\Gamma$ and $\delta g^{-1}$ act on $C_0$). With Lemma~\ref{lem:Bach0},
\[
D(\nabla^kC_{ijk})(fg_0)=\tfrac12(\nabla^kW_{mikj})\nabla^mf+\tfrac12W_{mikj}\nabla^k\nabla^mf=\tfrac12\,\mathring W(\nabla\nabla f)_{ij},
\]
since $\nabla^kW_{mikj}=\nabla^kW_{kjmi}=0$ and $W_{mikj}=W_{imjk}$.

For the second term we write $A^{kl}=g^{ka}g^{lb}A_{ab}$. The variation of $A_{ab}$ gives $-\frac12\mathring W(\nabla\nabla f)$. The variations of $g^{-1}$ give $-2pf\,g^{kl}W_{ikjl}=0$. The variation of $W$ gives $p\,g^{kl}DW_{ikjl}(fg_0)$; differentiating $g^{kl}W(g)_{ikjl}\equiv0$ we get $g^{kl}DW_{ikjl}(h)=h^{kl}W_{ikjl}$, which here equals $f\,g^{kl}W_{ikjl}=0$. Altogether $DB(fg_0)=\frac12\mathring W(\nabla\nabla f)-\frac12\mathring W(\nabla\nabla f)=0$.
\end{proof}

\begin{remark}
The cancellation of the terms $\mathring W(\nabla\nabla f)$ occurs in every dimension. What distinguishes dimension $4$ is the \emph{global} conformal covariance of $B$: for $n\neq4$, the transformation law of $B$ under $g\mapsto e^fg$ contains terms proportional to $n-4$, in $C\cdot\nabla f$ (see Proposition~\ref{prop:Bconfgen}) and $W(\nabla f,\nabla f)$, which respectively vanish and are quadratic at an Einstein metric.
\end{remark}

On TT tensors, the linearised Bach tensor factorises into two shifted Lichnerowicz Laplacians.

\begin{proposition}\label{prop:BTT}
For $h$ TT,
\[
DB(g_0)h=-\frac{1}{2(n-2)}\Big(\Delta_L-\frac{n}{n-1}\lambda\Big)(\Delta_L-2\lambda)h .
\]
\end{proposition}

\begin{proof}
\emph{Step 1: the tensor $\mathcal A$.} For $h$ TT, the Bianchi operator of \cite{D} satisfies $2\Div h+d\Tr h=0$, and $DR(g_0)h=d^*\Div h+\nabla^*\nabla\Tr h-\lambda\Tr h=0$, hence, by e.g. \cite{D},
\[
DA(h)=\frac{1}{n-2}D\Sch(h)=\frac{1}{2(n-2)}\Big(\Delta_L-\frac{n}{n-1}\lambda\Big)h .
\]
Since $2(n-2)p+\frac{n}{n-1}\lambda=\frac{(n-2)\lambda+n\lambda}{n-1}=2\lambda$,
\[
\mathcal A:=\mathcal A(h)=DA(h)-ph=\frac{1}{2(n-2)}(\Delta_L-2\lambda)h .
\]
$\mathcal A$ is TT, since $\Delta_L$ preserves TT tensors on an Einstein manifold, see e.g. \cite[Proposition 4]{D}.

\emph{Step 2: two identities on $\Sr$.} For $k$ TT,
\[
[\nabla_k,\nabla_j]k_i{}^k=\lambda k_{ij}-\mathring Rk_{ij},
\]
by the Ricci identity: the term $R^k{}_{mkj}k_i{}^m$ gives $\lambda k_{ij}$, and
\[
-R^m{}_{ikj}k_m{}^k=-R_{mikj}k^{mk}=-R_{imjk}k^{mk}=-\mathring Rk_{ij}.
\] Moreover, the definition of $\Delta_L$ gives $\nabla^*\nabla=\Delta_L-2\lambda+2\mathring R$ on $S_2$, and, for $k$ trace free, $\mathring Rk=\mathring Wk+(A_0\KN g_0)_{ikjl}k^{kl}=\mathring Wk-2pk$.

\emph{Step 3: the Cotton term.} Since $C(g_0)=0$, the variation of $\nabla^kC_{ijk}$ is $\nabla^kDC_{ijk}$, and by Lemma~\ref{lem:DnablaA}
\[
\nabla^kDC_{ijk}=\nabla^k\nabla_k\mathcal A_{ij}-\nabla^k\nabla_j\mathcal A_{ik}
=-\nabla^*\nabla \mathcal A_{ij}-\nabla_j\nabla^k\mathcal A_{ik}-[\nabla_k,\nabla_j]\mathcal A_i{}^k .
\]
The middle term vanishes since $\mathcal A$ is TT. With Step~2,
\[
\nabla^kDC=-(\Delta_L-2\lambda+2\mathring R)\mathcal A-\lambda \mathcal A+\mathring R\mathcal A=-\Delta_L\mathcal A+\lambda \mathcal A-\mathring W\mathcal A+2p\mathcal A .
\]

\emph{Step 4: the Weyl term.} We differentiate $A^{kl}W_{ikjl}=g^{ka}g^{lb}A_{ab}W_{ikjl}$. The variation of $A_{ab}$ gives $\mathring W(DA(h))$, the variations of $g^{-1}$ give $-2p\,\mathring Wh$, and the variation of $W$ gives $p\,g^{kl}DW_{ikjl}(h)=p\,h^{kl}W_{ikjl}=p\,\mathring Wh$, as in the proof of Proposition~\ref{prop:Bconf}. Altogether,
\[
D(A^{kl}W_{ikjl})(h)=\mathring W\big(DA(h)-ph\big)=\mathring W\mathcal A .
\]

\emph{Step 5: sum.} The terms in $\mathring W$ cancel exactly:
\[
DB(h)=-\Delta_L\mathcal A+(\lambda+2p)\mathcal A=-\Big(\Delta_L-\frac{n}{n-1}\lambda\Big)\mathcal A,
\]
since $\lambda+2p=\lambda+\frac{\lambda}{n-1}=\frac{n}{n-1}\lambda$. We conclude with Step~1.
\end{proof}

Comparing with the formula of \cite{D}, we obtain the identification announced in the introduction.

\begin{theorem}\label{thm:bach}
On every Einstein manifold of dimension $n\ge3$, as differential operators,
\[
P=-2(n-2)\,DB(g_0)\quad\text{on }\Sr,\qquad DB(g_0)(ug_0)=0,
\]
where $P$ is the operator of \cite[Théorème 1]{D}. In particular $DB(g_0)$ is self-adjoint, TT-valued, and
\[
DB(g_0)h=-\frac{1}{2(n-2)}P(\mathring h)\qquad\text{for all }h\in S_2 .
\]
\end{theorem}

\begin{proof}
The identity $DB(g_0)(ug_0)=0$ is Proposition~\ref{prop:Bconf}. By naturality and Lemma~\ref{lem:Bach0}, $DB(g_0)(LX)=L_XB(g_0)=0$; with Proposition~\ref{prop:Bconf}, $DB(g_0)\Lr=0$. Also $P\Lr=0$ \cite[Théorème 1]{D}.

On TT tensors, $P=(\Delta_L-2\lambda)(\Delta_L-\frac{n}{n-1}\lambda)$ \cite{D}, and the two factors commute; Proposition~\ref{prop:BTT} thus gives $P=-2(n-2)DB(g_0)$ on TT tensors.

Finally let $\mathring h\in\Sr$ and $U$ a small ball. By Lemma~\ref{lem:localTT}, $\mathring h=k+\Lr w$ on $U$ with $k$ TT. Both operators are local, agree on TT tensors and vanish on $\Lr w$, hence agree on $\mathring h$ in $U$. TT-valuedness follows from that of $P$. Finally, $P$ is formally self-adjoint on $\Sr$ \cite[Théorème 1]{D}: in the formula
\[
P=(\Delta_L-2\lambda)\Big(\Delta_L-\tfrac{n}{n-1}\lambda\Big)-\Lr\Big(d^*d+\tfrac{n}{2(n-1)}dd^*-\tfrac{n}{n-1}\lambda\Big)\Div
\]
of \cite{D}, the first term is self-adjoint since $\Delta_L$ is and preserves $\Sr$, and so is the second one, since $\Div^*=\frac12\Lr$ and $\Lr^*=2\Div$ on $\Sr$. As $DB(g_0)h=-\frac{1}{2(n-2)}P(\mathring h)$ and $h\mapsto\mathring h$ is self-adjoint, $DB(g_0)$ is self-adjoint on $S_2$.
\end{proof}

\begin{remark}
\leavevmode
\begin{enumerate}[label=(\alph*)]
\item For $n=4$ we recover $P=-4\,DB(g_0)$: the machine of \cite{D} is the linearised Bach tensor, as suggested in the introduction of \cite{D}. The theorem shows that this identification persists in every dimension, provided one uses the above definition of $B$, which remains of order $4$.
\item The theorem gives another proof that $P$ maps $\Sr$ to TT tensors, via Proposition~\ref{prop:principe} applied to $N=B$: $DB(g_0)^*=DB(g_0)$ is TT-valued.
\item Neither the gradient of $\int|W|^2$ \cite{GV}, which contains terms in $W\ast W$ proportional to $n-4$ and does not vanish at a generic Einstein metric, nor the obstruction tensor, of order $n$, give this order-$4$ identification for $n\neq4$. In Matsumoto's factorisation of $D\mathcal O(g_0)$ on TT tensors (Example~\ref{ex:liste}(c)), the factor $m=0$ is $\Delta_L-2\lambda$ and the last factor $m=\frac n2-1$ is $\Delta_L-\frac{n}{n-1}\lambda$: these are exactly the two factors of $P$, and the $\frac n2-2$ intermediate factors are absent from $P$. Theorem~\ref{thm:bach} shows that the order-$4$ tensor $B$ isolates these two factors in every dimension, even or odd.
\end{enumerate}
\end{remark}

\subsection{The factorisation of $P$ seen through the Bach tensor}

To understand this identification factor by factor, it is convenient to separate, in the Bach tensor, the dependence on the Schouten tensor from the dependence on the metric.

\begin{definition}
For a metric $g$, we define the linear second-order operator with values in $T^0_2$
\[
(\mathfrak B_gk)_{ij}:=\nabla^k\big(\nabla_kk_{ij}-\nabla_jk_{ik}\big)+W_{ikjl}k^{kl}=\nabla^k(d^\nabla k)_{ijk}+\mathring Wk_{ij},\qquad k\in S_2 ,
\]
where $(d^\nabla k)_{ijk}=\nabla_kk_{ij}-\nabla_jk_{ik}$, so that $C=d^\nabla A$. Its coefficients depend on $g$, but it is linear in $k$.
\end{definition}

The Bach tensor is obtained by applying this operator to $A$, and the metric itself is invisible to it.

\begin{lemma}\label{lem:BBfrak}
For every metric $g$ and every constant $q$, $B(g)=\mathfrak B_g\big(A(g)-q\,g\big)$.
\end{lemma}

\begin{proof}
By definition $B(g)=\mathfrak B_g(A(g))$. Since $\nabla g=0$ and $g^{kl}W_{ikjl}=0$, $\mathfrak B_g(g)=0$, and we conclude by linearity.
\end{proof}

Linearising this identity at an Einstein metric yields a factorisation of $DB(g_0)$.

\begin{proposition}\label{prop:factBach}
On an Einstein manifold, with $p=\frac{\lambda}{2(n-1)}$, for every $h\in S_2$,
\[
DB(g_0)h=\mathfrak B_{g_0}\big(DA(g_0)h-p\,h\big)=\mathfrak B_{g_0}\big(\mathcal A(h)\big).
\]
Moreover, for $k$ TT, $\mathfrak B_{g_0}k=-\big(\Delta_L-\frac{n}{n-1}\lambda\big)k$.
\end{proposition}

\begin{proof}
Differentiating the identity of Lemma~\ref{lem:BBfrak} with $q=p$ at $g_0$:
\[
DB(g_0)h=\mathfrak B_{g_0}\big(DA(g_0)h-ph\big)+\big(D_g\mathfrak B\big)(h)\big[A(g_0)-pg_0\big],
\]
and the second term vanishes since $A(g_0)=pg_0$. For $k$ TT, the computation of Step~3 in the proof of Proposition~\ref{prop:BTT}, which only uses that $\mathcal A$ is TT, gives $\nabla^k(d^\nabla k)_{ijk}=\big(-\Delta_Lk+\lambda k-\mathring Wk+2pk\big)_{ij}$; adding $\mathring Wk$ we get $\mathfrak B_{g_0}k=-\Delta_Lk+(\lambda+2p)k$, and $\lambda+2p=\frac{n}{n-1}\lambda$.
\end{proof}

\begin{remark}
Proposition~\ref{prop:factBach} gives a shorter proof of Steps~3 and~4 of Proposition~\ref{prop:BTT}: they reduce to the computation of $\mathfrak B_{g_0}$ on TT tensors.
\end{remark}

In terms of the Schouten tensor shifted by its value at $g_0$, we recover the factorisation of $P$ of \cite{D}, factor by factor.

\begin{corollary}\label{cor:factP}
Let $\varpi:=\frac{(n-2)\lambda}{2(n-1)}$, so that $\Sch(g_0)=\varpi g_0$, and
\[
\mathcal S:=D\big(\Sch-\varpi\,g\big)(g_0)=(n-2)\mathcal A .
\]
Then, on $\Sr$,
\[
P=-2\,\mathfrak B_{g_0}\circ\mathcal S,\qquad
\mathcal S=D\Ricr(g_0)+\frac{n-2}{2n(n-1)}\,DR(g_0)(\cdot)\,g_0 .
\]
On TT tensors the factors agree one by one with those of \cite[Proposition 16]{D}:
\[
\mathcal S=D\Ricr(g_0)^*=\tfrac12(\Delta_L-2\lambda),\qquad
-2\,\mathfrak B_{g_0}=4\,D\Sch(g_0)^*=2\Big(\Delta_L-\tfrac{n}{n-1}\lambda\Big).
\]
\end{corollary}

\begin{proof}
The first formula follows from Theorem~\ref{thm:bach} and Proposition~\ref{prop:factBach}, since $\mathcal S=(n-2)\mathcal A$. For the second, write
\[
\Sch-\varpi g=\Ricr+\Big(\frac{n-2}{2n(n-1)}R-\varpi\Big)g
\]
and differentiate at $g_0$; since $\frac{n-2}{2n(n-1)}R(g_0)=\varpi$, the term in $h$ disappears. The values on TT tensors follow from Step~1 of Proposition~\ref{prop:BTT}, Proposition~\ref{prop:factBach} and the formulas of \cite{D} for $D\Ricr^*$ and $D\Sch^*$.
\end{proof}

\begin{remark}[Interpretation]
\leavevmode
\begin{enumerate}[label=(\alph*)]
\item In $P=4D\Sch^*D\Ricr^*$, the inner factor $D\Ricr^*$ plays the role of a \emph{shifted} linearised Schouten tensor $\mathcal S$, and the outer factor $D\Sch^*$ that of the linear Bach operator $\mathfrak B_{g_0}$ (divergence of Cotton plus Weyl). Off TT tensors, the factors differ by gauge and trace terms: on $\Sr$, $\mathcal S-D\Ricr(g_0)=\frac{n-2}{2n(n-1)}(d^*\Div\,\cdot)\,g_0$ by Corollary~\ref{cor:factP}, and $\mathfrak B_{g_0}$ contains the non-symmetrised term $\nabla_j(\Div k)_i$, so it is not symmetric-valued on arbitrary $k$. These differences cancel in both compositions, which equal $P$ on $\Sr$.
\item Both constructions have the same structure: a curvature operator applied to the Schouten tensor. In \cite[\S6]{D} this is the non-linear composition $\Ricr(\Sch(g))$, which makes sense because $\Sch(g_0)$ is Einstein; here it is $\mathfrak B_g(\Sch(g)-\varpi g)$, linear in the Schouten tensor, which makes sense because $\mathfrak B_g$ kills $g$. In both cases the proportionality $\Sch(g_0)\propto g_0$ makes the term where the outer operator is differentiated disappear.
\item For a Bach-flat non-Einstein metric we only have
\[
DB(g_0)h=\mathfrak B_{g_0}(DA(g_0)h)+(D_g\mathfrak B)(h)[A(g_0)],
\]
and the second term no longer vanishes: there is no factorisation of this type.
\end{enumerate}
\end{remark}

\subsection{The partially massless spin-$2$ field in every dimension}

On an Einstein manifold (or spacetime) with $\Ric=\lambda g$, a spin-$2$ field of mass $m$ satisfies, on TT tensors and with the sign conventions of Section~\ref{sec:conv}, $(\Delta_L-2\lambda+m^2)h=0$; $m=0$ is the massless graviton. The value
\[
m^2_{\mathrm{PM}}=\frac{n-2}{n-1}\,\lambda
\]
is the \emph{partially massless} point, where the theory acquires the scalar gauge symmetry $\delta h=\nabla\nabla\alpha+\frac{\lambda}{n-1}\alpha\,g$ \cite{DWa,DJW}.

The factorisation of $P$ has a direct interpretation in these terms.

\begin{proposition}\label{prop:PM}
Let $g_0$ be Einstein, $\Ric(g_0)=\lambda g_0$, $n\ge3$, and $\mathcal S$, $\mathfrak B_{g_0}$ as in Corollary~\ref{cor:factP}.
\begin{enumerate}[label=(\roman*)]
\item On TT tensors, $P=(\Delta_L-2\lambda)\big(\Delta_L-2\lambda+m^2_{\mathrm{PM}}\big)$: $P$ is the product of the massless and the partially massless spin-$2$ operators.
\item For every function $f$,
\[
\mathcal S(fg_0)=-\frac{n-2}{2}\Big(\nabla\nabla f+\frac{\lambda}{n-1}\,f\,g_0\Big),
\]
that is, $\mathcal S$ maps conformal directions onto partially massless pure-gauge modes.
\item The linear Bach operator annihilates partially massless pure-gauge modes: for every function $\alpha$,
\[
\mathfrak B_{g_0}\Big(\nabla\nabla\alpha+\frac{\lambda}{n-1}\,\alpha\,g_0\Big)=0 .
\]
\end{enumerate}
These statements are algebraic and hold in any signature.
\end{proposition}

\begin{proof}
(i) $\frac{n}{n-1}\lambda=2\lambda-\frac{n-2}{n-1}\lambda$, so the second factor in Theorem~\ref{thm:bach} is $\Delta_L-2\lambda+m^2_{\mathrm{PM}}$. (ii) By \cite[Lemme 7]{D}, $D\Sch(fg_0)=-\frac{n-2}{2}\nabla\nabla f$, and $D(\varpi g)(fg_0)=\varpi fg_0$ with $\varpi=\frac{(n-2)\lambda}{2(n-1)}$. (iii) By Proposition~\ref{prop:factBach} and Corollary~\ref{cor:factP}, $\mathfrak B_{g_0}\circ\mathcal S=(n-2)\,DB(g_0)$, and $DB(g_0)(fg_0)=0$ by Proposition~\ref{prop:Bconf}; apply (ii). The proofs only use algebraic identities and naturality, which do not depend on the signature.
\end{proof}

\begin{remark}
In dimension $4$, the relation between the Bach tensor and the partially massless graviton is known from conformal (Weyl) gravity, whose linearisation on an Einstein background contains a massless and a partially massless mode \cite{DJW}. Proposition~\ref{prop:PM} extends the structural statement to every dimension, with the order-$4$ tensor $B$ which is not conformally covariant for $n\neq4$: the conformal invariance needed is only the infinitesimal one at Einstein metrics (Proposition~\ref{prop:Bconf}), and it becomes the partially massless gauge invariance of $\mathfrak B_{g_0}$. In $D$-dimensional quadratic (critical) gravity, the linearised field equations around an Einstein background are, on TT tensors, of the form $(\Delta_L-2\lambda)(\Delta_L-2\lambda+M^2)h=0$, where the mass $M$ depends on the couplings \cite{DLLPST,LPP}; the partially massless value $M^2=m^2_{\mathrm{PM}}$ is reached for particular values of the couplings, and for $D=4$ it is the case of conformal gravity \cite{LPP,DJW}. Gover, Waldron and coauthors relate masses, and in particular partial masslessness, to Weyl invariance through tractor calculus \cite{GSW}. Proposition~\ref{prop:PM} is of a different nature: the PM factor is produced, in every dimension, by a single natural tensor, which for $n\neq4$ is not conformally covariant; for even $n$, Matsumoto's factorisation of the linearised obstruction tensor \cite{Ma}, recalled in Example~\ref{ex:liste}(c), contains both the massless and the PM factors (for $m=0$ and $m=\frac n2-1$), together with $\frac n2-2$ other factors. For $n\neq4$, $B$ is not the Euler--Lagrange tensor of $\int|W|^2$, so the statement is structural rather than about a variational theory; $P$ is nevertheless formally self-adjoint on $\Sr$ (Theorem~\ref{thm:bach}).
\end{remark}

\section{Bach-flat metrics}\label{sec:bachplat}

We no longer assume that $g_0$ is Einstein.

We start with the linearised Cotton tensor in conformal directions, which can be computed for an arbitrary metric.

\begin{lemma}\label{lem:DCconf}
For every metric $g_0$ and every function $f$,
\[
DC(g_0)(fg_0)_{ijk}=\tfrac12\,W_{mikj}\nabla^mf .
\]
\end{lemma}

\begin{proof}
By \cite[Lemme 7]{D}, $DA(fg_0)=-\frac12\nabla\nabla f$. By \eqref{eq:DnablaT}, and since $\delta\Gamma$ is symmetric in its lower indices,
\[
DC_{ijk}=\nabla_k\delta A_{ij}-\nabla_j\delta A_{ik}-\delta\Gamma^m_{ki}A_{mj}+\delta\Gamma^m_{ji}A_{mk},
\]
with $\delta A=-\frac12\nabla\nabla f$ and $\delta\Gamma^m_{ki}=\frac12(f_k\delta^m_i+f_i\delta^m_k-f^mg_{ki})$, where $f_i=\nabla_if$. On the one hand, the Ricci identity for $df$ gives
\[
\nabla_k\delta A_{ij}-\nabla_j\delta A_{ik}=-\tfrac12[\nabla_k,\nabla_j]\nabla_if=\tfrac12R_{mikj}f^m .
\]
On the other hand, the terms in $f_i$ cancel and
\[
-\delta\Gamma^m_{ki}A_{mj}+\delta\Gamma^m_{ji}A_{mk}=\tfrac12\big(f_jA_{ik}-f_kA_{ij}-A(\nabla f)_kg_{ij}+A(\nabla f)_jg_{ik}\big).
\]
With $R=W+A\KN g$, $\frac12(A\KN g)_{mikj}f^m=\frac12\big(A(\nabla f)_kg_{ij}+A_{ij}f_k-A(\nabla f)_jg_{ik}-A_{ik}f_j\big)$, which exactly cancels the last term.
\end{proof}

This leads to an exact formula for the linearised Bach tensor in conformal directions.

\begin{proposition}\label{prop:Bconfgen}
For every metric $g_0$, every function $f$ and every $n\ge3$,
\[
DB(g_0)(fg_0)_{ij}=-f\,B(g_0)_{ij}+\frac{n-4}{2}\big(C_{ijm}+C_{jim}\big)\nabla^mf,\qquad C=C(g_0).
\]
\end{proposition}

\begin{proof}
We write $f_i=\nabla_if$, $h=fg_0$, so that $\delta g^{kl}=-fg^{kl}$ and $\delta\Gamma^m_{li}=\frac12(f_l\delta^m_i+f_i\delta^m_l-f^mg_{li})$. We use three facts: the cyclic identity $C_{ijk}+C_{jki}+C_{kij}=0$; the vanishing of the traces $g^{ij}C_{ijk}=g^{ik}C_{ijk}=0$, which follows from $\nabla^iA_{ij}=\nabla_j\Tr A$; and
\[
\nabla^aW_{abcd}=(n-3)\,C_{bdc},
\]
which follows from $\nabla^aR_{abcd}=\nabla_c\Ric_{bd}-\nabla_d\Ric_{bc}$, $R=W+A\KN g$ and $\Ric=(n-2)A+(\Tr A)g$.

\emph{Cotton term.} By \eqref{eq:DnablaT},
\[
D\big(\nabla^kC_{ijk}\big)(h)=-f\,\nabla^kC_{ijk}+\nabla^kDC(h)_{ijk}-g^{kl}\big(\delta\Gamma^m_{li}C_{mjk}+\delta\Gamma^m_{lj}C_{imk}+\delta\Gamma^m_{lk}C_{ijm}\big).
\]
Using the vanishing of the traces of $C$ and its antisymmetry in the last two indices, the three $\delta\Gamma$ terms are respectively $\frac12(f^mC_{ijm}-f^mC_{mji})$, $f^mC_{ijm}$ and $\frac{2-n}{2}f^mC_{ijm}$. By Lemma~\ref{lem:DCconf} and the identity above,
\[
\nabla^kDC(h)_{ijk}=\tfrac12\nabla^k\big(W_{mikj}f^m\big)=\tfrac{n-3}{2}C_{jim}f^m+\tfrac12\mathring W(\nabla\nabla f)_{ij},
\]
since $\nabla^kW_{mikj}=\nabla^kW_{kjmi}=(n-3)C_{jim}$ and $W_{mikj}=W_{imjk}$. Hence
\[
D\big(\nabla^kC_{ijk}\big)(h)=-f\,\nabla^kC_{ijk}+\tfrac12\mathring W(\nabla\nabla f)_{ij}+\tfrac12f^m\big((n-3)C_{jim}+C_{mji}+(n-5)C_{ijm}\big).
\]
By the cyclic identity, $C_{mji}=-C_{jim}-C_{imj}=-C_{jim}+C_{ijm}$, so the last term equals $\frac{n-4}{2}(C_{ijm}+C_{jim})f^m$.

\emph{Weyl term.} In $A^{kl}W_{ikjl}=g^{ka}g^{lb}A_{ab}W_{ikjl}$, the variation of $A$ gives $-\frac12\mathring W(\nabla\nabla f)$ since $DA(g_0)(fg_0)=-\frac12\nabla\nabla f$ \cite[Lemme 7]{D}, those of $g^{-1}$ give $-2f\,\mathring WA_0$, and, since $W(e^fg)=e^fW(g)$, the variation of $W$ gives $f\,\mathring WA_0$.

Summing, the terms $\mathring W(\nabla\nabla f)$ cancel, and $-f\nabla^kC_{ijk}-f\,\mathring WA_0=-fB(g_0)$.
\end{proof}

Away from dimension $4$, the vanishing of the Cotton tensor is therefore exactly what is needed.

\begin{corollary}\label{cor:Bconfgen}
Let $n\neq4$ and $B(g_0)=0$. Then $DB(g_0)(fg_0)=0$ for every function $f$ if and only if $C(g_0)=0$.
\end{corollary}

\begin{proof}
By Proposition~\ref{prop:Bconfgen}, the condition is equivalent to $C_{ijm}+C_{jim}=0$. A $3$-tensor which is antisymmetric in its first two and in its last two indices is totally antisymmetric; the cyclic identity then gives $3C_{ijk}=0$.
\end{proof}

Since $B$ is natural of order $4$, its symbol does not depend on the background metric and can be read off from the flat case.

\begin{lemma}\label{lem:symbB}
For every metric $g_0$, the principal symbol of $DB(g_0)$ satisfies
\[
\sigma_\xi\big(DB(g_0)\big)k=-\frac{|\xi|^4}{2(n-2)}\,k\qquad\text{for all }k\in\TT_\xi .
\]
\end{lemma}

\begin{proof}
$B$ is natural of order $4$: in normal coordinates, the fourth-order part of $DB(g_0)h$ involves only $g_0$ at the point, hence agrees with the flat case. In particular, the term $A^{kl}W_{ikjl}$, quadratic in the curvature, contributes to $DB(g_0)h$ with at most two derivatives of $h$, and does not enter the symbol of order $4$. For the flat metric (Einstein with $\lambda=0$), Theorem~\ref{thm:bach} gives $DB=-\frac{1}{2(n-2)}P$ on $\Sr$, and the symbol of $P$ on $\TT_\xi$ is $|\xi|^4$, since the term $\Lr(\cdots)\Div$ of $P$ has vanishing symbol on $\TT_\xi$.
\end{proof}

Putting these facts together with the principle, we obtain machines on Bach-flat metrics.

\begin{theorem}\label{thm:bachplat}
Let $(M,g_0)$ be such that
\begin{enumerate}[label=(\alph*)]
\item either $n=4$ and $B(g_0)=0$;
\item or $n\ge3$ and $B(g_0)=0$, $C(g_0)=0$.
\end{enumerate}
Then $DB(g_0)^*:\Gamma(S_2)\to\TT$ is a TT machine of order $4$. If $M$ is compact without boundary,
\[
\TT=\Ima DB(g_0)^*\oplus\big(\ker DB(g_0)\cap\TT\big),
\]
the second summand, consisting of infinitesimal Bach-flat TT deformations, being finite dimensional. In dimension $4$, $DB(g_0)$ is moreover self-adjoint, so that $DB(g_0)$ itself is TT-valued.
\end{theorem}

\begin{proof}
We apply Proposition~\ref{prop:principe} with $N=B$. Condition (i) is $B(g_0)=0$. Condition (ii) follows from Proposition~\ref{prop:Bconfgen}, since $B(g_0)=0$ and either $n=4$ or $C(g_0)=0$; for $n=4$ it is also a consequence of the conformal covariance $B(e^{f}g)=e^{-f}B(g)$ \cite{Be}. By Corollary~\ref{cor:Bconfgen}, for $n\neq4$ the assumption $C(g_0)=0$ in case (b) cannot be removed from condition (ii). Theorem~\ref{thm:fredholm} applies by Lemma~\ref{lem:symbB}. In dimension $4$, $B$ is, up to a constant, the gradient of $\int|W|^2dV$, and $DB(g_0)$ is the Hessian of this functional at a critical point, hence self-adjoint.
\end{proof}

\begin{remark}
At an Einstein metric, $DB(g_0)$ is self-adjoint in every dimension, since it is a multiple of $P$ (Theorem~\ref{thm:bach}). In dimension $4$ it is self-adjoint at every Bach-flat metric, $B$ being the gradient of $\int|W|^2$. For $n\neq4$ and a non-Einstein metric in case (b), $B$ is not a gradient, and we only use the fact that $DB(g_0)^*$ is a machine.
\end{remark}

\begin{example}
\leavevmode
\begin{enumerate}[label=(\alph*)]
\item Case (a): besides conformally Einstein metrics, self-dual metrics (for instance LeBrun's metrics on $\#k\,\mathbb{CP}^2$) and Bach-flat Kähler metrics. These are in general not conformally Einstein, and the machine $DB(g_0)$ is then not obtained by conformal transport of $P$.
\item Case (b): Einstein metrics (Section~\ref{sec:bach}) and conformally flat metrics in every dimension $n\ge3$ (for $n\ge4$, $W=0$ implies $C=0$ and then $B=0$). On a space form, $P_W=-2(n-3)DB(g_0)$ on $\Sr$ by Theorems~\ref{thm:formes} and~\ref{thm:bach}.
\end{enumerate}
\end{example}

\section{Conformally flat metrics, $n\ge4$}\label{sec:confplat}

Throughout this section $W(g_0)=0$, so that $\Rm_0=A\KN g_0$.

On a conformally flat manifold, the Weyl tensor is a natural tensor vanishing at $g_0$, and its linearisation kills conformal directions.

\begin{lemma}\label{lem:Wconf}
The operator $N=W$ satisfies the assumptions of Proposition~\ref{prop:principe}. Moreover $DW(g_0)h\in\cW$, and for $Z\in\cW$,
\[
\ip{DW(g_0)h}{Z}=\ip{DR(h)}{Z}-\ip{A\KN h}{Z}.
\]
\end{lemma}

\begin{proof}
$W$ is natural and $W(e^fg)=e^fW(g)$ (as a $(0,4)$-tensor, \cite[Ch.~1]{Be}), so $DW(g_0)(fg_0)=fW(g_0)=0$. Differentiating the algebraic identities satisfied by $W(g)$, in particular $\Tr_gW(g)=0$, gives $DW(h)\in\cW$ since $W(g_0)=0$. Finally $DW(h)=DR(h)-DA(h)\KN g_0-A\KN h$, and $\ip{U\KN g_0}{Z}=4U_{ac}g_{bd}Z^{abcd}=0$ for $Z\in\cW$.
\end{proof}

The resulting machine can be written explicitly.

\begin{proposition}\label{prop:DWstar}
For $Z\in\Gamma(\cW)$,
\[
\big(DW(g_0)^*Z\big)_{ad}=2\big(\nabla^b\nabla^c Z_{abcd}+A^{bc}Z_{abcd}\big).
\]
This expression is symmetric in $(a,d)$, trace free and divergence free; more precisely, the commutator $[\nabla^b,\nabla^c]Z_{abcd}$ vanishes.
\end{proposition}

\begin{proof}
We pair the formula of Lemma~\ref{lem:DRm} with a compactly supported $Z\in\Gamma_c(\cW)$.

\emph{Second-order terms.} By the antisymmetries of $Z$ in $(ab)$ and in $(cd)$, the four terms give the same contribution, namely $2\int\nabla_c\nabla_bh_{ad}Z^{abcd}=2\int h_{ad}\nabla_b\nabla_cZ^{abcd}$.

\emph{Commutator.} $[\nabla_c,\nabla_d]h_{ab}$ is symmetric in $(a,b)$, whereas $Z$ is antisymmetric there; the contribution vanishes.

\emph{The term $h_{am}R^m{}_{bcd}$.} With $R_{mbcd}=(A\KN g)_{mbcd}$ and $Z$ totally trace free, only two terms survive, and they are equal:
\[
R_{mbcd}Z^{abcd}=A_{bd}Z^{ab}{}_m{}^d-A_{bc}Z^{abc}{}_m=2A_{bd}Z^{ab}{}_m{}^d .
\]
Since $Z_{abmd}A^{bd}=(Z\circ A)_{am}$, the contribution is $2\ip{h}{Z\circ A}$.

\emph{The term $-\ip{A\KN h}{Z}$.} The four terms of the Kulkarni--Nomizu product give the same contribution, whence $-4A_{ac}h_{bd}Z^{abcd}=-4\ip{h}{Z\circ A}$.

Altogether $\ip{DWh}{Z}=\int h_{ad}\big(2\nabla_b\nabla_cZ^{abcd}-2(Z\circ A)^{ad}\big)$, and since $(Z\circ A)_{ad}=Z_{bacd}A^{bc}=-Z_{abcd}A^{bc}$, we obtain the formula.

\emph{The commutator $[\nabla^b,\nabla^c]Z_{abcd}$.} We use the Ricci identity on each of the four indices, with $R^m{}_{epq}=A^m{}_pg_{eq}+A_{eq}\delta^m_p-A^m{}_qg_{ep}-A_{ep}\delta^m_q$. Since $Z$ is trace free and antisymmetric, most terms vanish, and the four indices contribute respectively
\[
(Z\circ A)_{ad},\quad -(n-2)(Z\circ A)_{ad},\quad (n-2)(Z\circ A)_{ad},\quad -(Z\circ A)_{ad},
\]
whose sum is zero.

\emph{Symmetry.} By the pair symmetry, $Z_{dbca}=Z_{acbd}$, so $\nabla^b\nabla^cZ_{dbca}=\nabla^c\nabla^bZ_{abcd}$, which equals $\nabla^b\nabla^cZ_{abcd}$ by the previous point. The term $A^{bc}Z_{abcd}$ is symmetric for the same reason.

\emph{Trace and divergence.} The trace vanishes since $Z$ is trace free; the divergence vanishes by Proposition~\ref{prop:principe}.
\end{proof}

Its symbol is injective on TT symbols as soon as $n\ge4$.

\begin{proposition}[Symbol]\label{prop:symbW}
For $\xi\neq0$ and $k\in\TT_\xi$,
\[
|\sigma_\xi(DW(g_0))k|^2=\tfrac{n-3}{n-2}\,|\xi|^4|k|^2 .
\]
In particular this symbol is injective on $\TT_\xi$ for $n\ge4$.
\end{proposition}

\begin{proof}
By Lemma~\ref{lem:DRm}, $\sigma_\xi(DW)k=\frac12\cW(T)$ with $T=\xi\xi\KN k$. We compute
\[
\mathrm c(T)_{bd}=|\xi|^2k_{bd}+\xi_b\xi_d\Tr k-\xi_d k(\xi^\sharp)_b-\xi_bk(\xi^\sharp)_d=|\xi|^2k_{bd},
\]
which is trace free, so that $\Tr\mathrm c(T)=0$. Take $|\xi|=1$ and $\xi=e_1$, so that $k$ lives on $e_1^\perp$. The only non-zero components of $T$ are
\[
T_{1a1b}=T_{a1b1}=k_{ab},\qquad T_{1ab1}=T_{a11b}=-k_{ab},
\]
so $|T|^2=4|k|^2$. By \eqref{eq:normes}, $|\cW(T)|^2=4|k|^2-\frac{4}{n-2}|k|^2$, which gives the claim.

Directly: if $\cW(T)=0$, then $T=\frac{1}{n-2}k\KN g$, and the $(1,a,1,b)$ component gives $k_{ab}=\frac{1}{n-2}k_{ab}$, hence $k=0$ for $n\neq3$.
\end{proof}

We can now state the result for conformally flat manifolds.

\begin{theorem}\label{thm:confplat}
Let $(M,g_0)$ be conformally flat, $n\ge4$. Then $DW(g_0)^*:\Gamma(\cW)\to\TT$ is a TT machine of order $2$. If moreover $M$ is compact without boundary,
\[
\TT=\Ima DW(g_0)^*\oplus\big(\ker DW(g_0)\cap\TT\big),
\]
the second summand being finite dimensional. The operator $P_W:=DW(g_0)^*DW(g_0)$ is a TT machine of order $4$ on $S_2$; it is self-adjoint and satisfies $P_W\Lr=0$ and $P_W(ug_0)=0$.
\end{theorem}

\begin{proof}
Combine Lemma~\ref{lem:Wconf}, Proposition~\ref{prop:principe}, Theorem~\ref{thm:fredholm} and Proposition~\ref{prop:symbW}.
\end{proof}

\subsection{Space forms}

Assume now $R_{abcd}=\kappa(g_{ac}g_{bd}-g_{ad}g_{bc})$ and set $\lambda=(n-1)\kappa$. Then $A=\frac\kappa2g$, $\Rm_0=\frac\kappa2g\KN g$, and $\Delta_L=\nabla^*\nabla+2n\kappa$ on $\Sr$. We denote by $\Theta(h)$ the second-order part of $DR(h)$ in Lemma~\ref{lem:DRm}, that is,
\[
\Theta(h)_{abcd}=\tfrac12\big(\nabla_c\nabla_bh_{ad}+\nabla_d\nabla_ah_{bc}-\nabla_c\nabla_ah_{bd}-\nabla_d\nabla_bh_{ac}\big).
\]

To compare $P_W$ with $P$, we need two lemmas on the linearised curvature at TT tensors.

\begin{lemma}\label{lem:FE1}
\leavevmode
\begin{enumerate}[label=(\alph*)]
\item For every $h\in S_2$, $DR(h)=\Theta(h)+\frac\kappa2h\KN g$, and $\Theta(h)$ is an algebraic curvature tensor.
\item For $h$ TT, $c(\Theta(h))=\frac12(\Delta_L-n\kappa)h=:\rho$, and
\[
DW(h)=\Theta(h)-\tfrac{1}{n-2}\rho\KN g .
\]
\item For $k$ TT, $[\nabla_b,\nabla_a]k^b{}_d=n\kappa\,k_{ad}$.
\end{enumerate}
\end{lemma}

\begin{proof}
(a) We have $h_{am}R^m{}_{bcd}=\kappa(h_{ac}g_{bd}-h_{ad}g_{bc})$ and
\[
\tfrac12[\nabla_c,\nabla_d]h_{ab}=-\tfrac\kappa2\big(h_{bc}g_{ad}-h_{bd}g_{ac}+h_{ac}g_{bd}-h_{ad}g_{bc}\big),
\]
whose sum is $\frac\kappa2(h\KN g)$. Since $DR(h)$ and $h\KN g$ are algebraic curvature tensors, so is $\Theta(h)$.

(c) By the Ricci identity, $[\nabla_b,\nabla_a]k^b{}_d=\Ric_{ma}k^m{}_d-R^m{}_{dba}k^b{}_m=(n-1)\kappa k_{ad}+\kappa k_{ad}$, using $\Tr k=0$.

(b) For $h$ TT, $g^{ac}\Theta_{abcd}=\frac12\big(\nabla_c\nabla_bh^c{}_d+\nabla^*\nabla h_{bd}\big)$, since the terms involving $\Div h$ or $\Tr h$ vanish. The first term equals $\nabla_b\nabla_ch^c{}_d+[\nabla_c,\nabla_b]h^c{}_d=n\kappa\,h_{bd}$ by (c). Hence $c(\Theta)=\frac12(\nabla^*\nabla+n\kappa)h=\frac12(\Delta_L-n\kappa)h$. Moreover $DW=DR-DA\KN g-A\KN h=\Theta-DA\KN g$, since $A\KN h=\frac\kappa2 g\KN h$. As $DW$ is trace free, $c(\Theta)=c(DA\KN g)=(n-2)DA+(\Tr DA)\,g$, and since $\rho$ is trace free, $DA=\rho/(n-2)$.
\end{proof}

\begin{proof}
The second Bianchi identity $\nabla_eR_{abcd}+\nabla_cR_{abde}+\nabla_dR_{abec}=0$ \cite[Ch.~1]{Be}, contracted in $(e,a)$, gives
\[
\nabla^aR_{abcd}=\nabla_c\Ric_{bd}-\nabla_d\Ric_{bc}.
\]
We linearise it at $g_0$ using \eqref{eq:DnablaT}, $\nabla\Rm_0=0$, $\delta\Gamma\cdot(g\KN g)=2\nabla(h\KN g)$ and $\delta\Gamma\cdot\Ric_0=(n-1)\kappa\nabla h$:
\[
\nabla^a\big(DR(h)-\kappa h\KN g\big)_{abcd}=\nabla_c\big(D\Ric(h)-(n-1)\kappa h\big)_{bd}-(c\leftrightarrow d).
\]
For $h$ TT, $D\Ric(h)=\frac12\Delta_Lh$ \cite{D}, and $\nabla^a(h\KN g)_{abcd}=\nabla_ch_{bd}-\nabla_dh_{bc}$. With Lemma~\ref{lem:FE1}(a) we obtain the first identity.

Since $\Theta$ is an algebraic curvature tensor, the pair symmetry gives $\nabla^c\Theta_{abcd}=\nabla^c\Theta_{cdab}=\nabla_a\phi_{db}-\nabla_b\phi_{da}$. Hence
\[
\nabla^b\nabla^c\Theta_{abcd}=\nabla^b\nabla_a\phi_{bd}+\nabla^*\nabla \phi_{ad}=(\nabla^*\nabla+n\kappa)\phi_{ad},
\]
by Lemma~\ref{lem:FE1}(c), using that $\phi$ is TT since $\Delta_L$ preserves TT tensors \cite[Proposition 4]{D}.

For $k\KN g$, only the terms $k_{bd}g_{ac}$ and $-k_{ad}g_{bc}$ survive; they give $\nabla^b\nabla_ak_{bd}+\nabla^*\nabla k_{ad}=(\nabla^*\nabla+n\kappa)k_{ad}$.
\end{proof}

On space forms, the two machines turn out to be proportional.

\begin{theorem}\label{thm:formes}
On a space form of dimension $n\ge4$,
\[
P_W\big|_{\TT}=\tfrac{n-3}{n-2}\,(\Delta_L-2\lambda)\Big(\Delta_L-\tfrac{n}{n-1}\lambda\Big),
\]
and, as differential operators on $\Sr$, $P_W=\frac{n-3}{n-2}\,P$, where $P$ is the operator of \cite[Théorème 1]{D}.
\end{theorem}

\begin{proof}
Let $h$ be TT. By Proposition~\ref{prop:DWstar} and Lemma~\ref{lem:FE1}(b), and since $A^{bc}Z_{abcd}=\frac\kappa2g^{bc}Z_{abcd}=0$,
\[
P_Wh=2\nabla^b\nabla^c\Big(\Theta-\tfrac{1}{n-2}\rho\KN g\Big)_{abcd}=2(\Delta_L-n\kappa)\Big(\phi-\tfrac{\rho}{n-2}\Big),
\]
by Lemma~\ref{lem:FE2}. Since $(2n-3)(n-2)-n=2(n-1)(n-3)$,
\[
\phi-\frac{\rho}{n-2}=\frac12\cdot\frac{n-3}{n-2}\big(\Delta_L-2(n-1)\kappa\big)h ,
\]
and we conclude with $2(n-1)\kappa=2\lambda$ and $n\kappa=\frac{n}{n-1}\lambda$.

Both $P_W$ and $P$ vanish on $\Ima\Lr$ (Theorem~\ref{thm:confplat} and \cite[Théorème 1]{D}). By Lemma~\ref{lem:localTT}, on small balls every $\mathring h\in\Sr$ can be written $h_{\TT}+\Lr w$ with $h_{\TT}$ TT. Both operators being local, the equality holds on $\Sr$.
\end{proof}

\begin{remark}
For $\kappa=0$ we recover $P_W|_{\TT}=\frac{n-3}{n-2}(\nabla^*\nabla)^2$, consistently with \cite[Remarque 15]{D}. The factor $n-3$ reflects the vanishing of $W$ in dimension $3$.
\end{remark}

\section{Metrics with parallel Ricci tensor}\label{sec:ricpar}

Throughout this section, $\Ric(g_0)$ is parallel, with distinct eigenvalues $\lambda_1,\dots,\lambda_{\ell}$ ($\ell\ge1$). The eigenspaces $E_a$ are parallel, mutually orthogonal distributions; we denote by $\Pi_a$ the corresponding orthogonal projections. They are parallel, $\nabla\Pi_a=0$, since they are polynomials in the parallel endomorphism $\Ric_0^\sharp$ (by Lagrange interpolation), and we set $g_a:=g_0(\Pi_a\cdot,\Pi_a\cdot)$, $g_a^{-1}:=\Pi_a\,g_0^{-1}$, that is $g_a^{kp}=(\Pi_a)^k{}_mg_0^{mp}$. For $h\in S_2$, the \emph{partial trace}
\[
\Tr_{g_a}h:=g_a^{kp}h_{kp}=\sum_ih(e_i,e_i),\qquad (e_i)\ \text{an orthonormal basis of }E_a,
\]
is the trace of the restriction of $h$ to $E_a$; thus $\sum_a\Tr_{g_a}h=\Tr_{g_0}h$, and $\Tr_{g_a}$ commutes with $\nabla$ since $g_a^{-1}$ is parallel. By de Rham's theorem \cite{dR}, \cite[Ch.~IV]{KN}, $g_0$ is locally a Riemannian product $g_0=g_1\oplus\dots\oplus g_{\ell}$ of Einstein manifolds $(M_a,g_a)$, $\Ric(g_a)=\lambda_ag_a$; the factors may themselves be reducible. We identify $g_a$ with the (degenerate) $2$-tensor $g_0(\Pi_a\cdot,\Pi_a\cdot)$ on $M$, so that $g_1\oplus\dots\oplus g_{\ell}=\sum_ag_a$, and similarly write $\Sch_0=\mu_1g_1\oplus\dots\oplus\mu_{\ell}g_{\ell}=\sum_a\mu_ag_a$. In particular, for $u\in E_c$, $R(u,v)w$ depends only on the $E_c$-components of $v,w$ and lies in $E_c$. We set $n_a=\dim E_a$ and
\[
\mu_a=\lambda_a-\frac{R}{2(n-1)},\qquad \Lambda:=\prod_c\lambda_c,\qquad \Lambda_a:=\prod_{c\neq a}\lambda_c ,
\]
with the convention $\Lambda_1=1$ if $\ell=1$. No assumption is made on the $\mu_a$: the Schouten tensor may be degenerate, or even zero.

\subsection{Linearised partial Ricci tensors}

For $h\in S_2$, let $\delta\Gamma(h)^m_{jk}=\frac12g_0^{mp}(\nabla_jh_{pk}+\nabla_kh_{pj}-\nabla_ph_{jk})$, so that, with the conventions of Section~\ref{sec:conv},
\[
D\Ric(g_0)h_{jl}=\nabla_k\delta\Gamma(h)^k_{lj}-\nabla_l\delta\Gamma(h)^k_{kj}.
\]

We split the linearised Ricci tensor according to the eigenspaces of $\Ric_0$.

\begin{definition}
The \emph{linearised partial Ricci tensor} of index $a$ is the second-order operator
\[
\mathcal R_a(h)_{jl}:=(\Pi_a)^k{}_m\Big(\nabla_k\delta\Gamma(h)^m_{lj}-\nabla_l\delta\Gamma(h)^m_{kj}\Big),
\]
obtained by restricting to $E_a$ the contraction defining $D\Ric$. One has $\sum_a\mathcal R_a=D\Ric(g_0)$.
\end{definition}

These operators have the following elementary properties.

\begin{lemma}\label{lem:Ea}
\leavevmode
\begin{enumerate}[label=(\alph*)]
\item $\mathcal R_a(h)$ is symmetric, and
\[
\mathcal R_a(h)_{jl}=\tfrac12g_a^{kp}\nabla_k\big(\nabla_lh_{pj}+\nabla_jh_{pl}-\nabla_ph_{lj}\big)-\tfrac12\nabla_l\nabla_j\big(\Tr_{g_a}h\big).
\]
\item For $\xi\in T^*M$ and $h\in S_2$, with $s_a=|\xi_a|^2$, $\tau_a=\Tr_{g_a}h$, $\eta_a=h(\xi_a^\sharp,\cdot)$, where $\xi_a=\xi\circ \Pi_a$,
\[
\sigma_\xi(\mathcal R_a)h=\tfrac12\big(s_ah+\tau_a\,\xi\otimes\xi-\xi\odot\eta_a\big),\qquad \xi\odot\eta=\xi\otimes\eta+\eta\otimes\xi .
\]
\end{enumerate}
\end{lemma}

\begin{proof}
(a) Since $(\Pi_a)^k{}_mg_0^{mp}=g_a^{kp}$, the first term gives the first part of the formula. In the second, $g_a^{kp}(\nabla_kh_{pj}-\nabla_ph_{kj})=0$ by symmetry of $g_a^{kp}$, and what remains is $\frac12\nabla_l\nabla_j\Tr_{g_a}h$, since $g_a^{-1}$ is parallel. Both terms are symmetric in $(j,l)$.
(b) Replace $\nabla_p\nabla_q$ with $-\xi_p\xi_q$ in (a).
\end{proof}

The construction rests on the following identity, which says that $\mathcal R_a$ only sees Lie derivatives along the factor $a$.

\begin{lemma}[Key identity]\label{lem:cle}
For every vector field $Y$ and all $a,b$,
\[
\mathcal R_a\big(L_Yg_b\big)=\delta_{ab}\,\lambda_a\,L_Yg_a .
\]
\end{lemma}

\begin{proof}
For $c=(c_1,\dots,c_{\ell})\in(0,\infty)^\ell$, the metric $G_c:=c_1g_1\oplus\dots\oplus c_{\ell}g_{\ell}$ is Riemannian and $\nabla$-parallel; its Levi-Civita connection is therefore $\nabla$, and its tensor $R(G_c)^m{}_{jkl}$ agrees with that of $g_0$. With the projection $\Pi_a$ as background structure, we consider
\[
C_a(G;\Pi_a)_{jl}:=(\Pi_a)^k{}_m\,R(G)^m{}_{jkl},
\]
which is natural in $(G,\Pi_a)$ jointly. At $G=G_c$, $C_a(G_c;\Pi_a)$ is the Ricci tensor of the factor $a$, that is $\lambda_ag_a$, independently of $c$.

\emph{Linearisation.} By Lemma~\ref{lem:DRm} in the form $\delta R^m{}_{jkl}=\nabla_k\delta\Gamma^m_{lj}-\nabla_l\delta\Gamma^m_{kj}$, with $\delta\Gamma$ computed for $G_c$, i.e.\ with $G_c^{-1}$ instead of $g_0^{-1}$, and since $\Pi_aG_c^{-1}=c_a^{-1}g_a^{-1}$, we get $D_GC_a(G_c;\Pi_a)=c_a^{-1}\mathcal R_a$.

\emph{Failure of naturality.} Since $\nabla\Pi_a=0$, $L_Y\Pi_a=[\Pi_a,\nabla Y]$, where $(\nabla Y)(u)=\nabla_uY$. This endomorphism maps $E_a$ into $E_a^\perp$ and $E_a^\perp$ into $E_a$, hence each $E_c$ into $E_c^\perp$. Now
\[
D_{\Pi}C_a(G_c;\Pi_a)(L_Y\Pi_a)_{jl}=(L_Y\Pi_a)^k{}_mR^m{}_{jkl}=\Tr\big(u\mapsto(L_Y\Pi_a)R(u,e_l)e_j\big),
\]
and for $u\in E_c$, $R(u,e_l)e_j\in E_c$, so $(L_Y\Pi_a)R(u,e_l)e_j\in E_c^\perp$: the trace, computed in an adapted basis, vanishes.

\emph{Conclusion.} Differentiating $C_a(\varphi_t^*G_c;\varphi_t^*\Pi_a)=\varphi_t^*C_a(G_c;\Pi_a)$ gives
\[
c_a^{-1}\sum_bc_b\,\mathcal R_a(L_Yg_b)=\lambda_a\,L_Yg_a\qquad\text{for all }c\in(0,\infty)^\ell .
\]
After multiplication by $c_a$ both sides are linear in $c$ on an open set, whence the identity coefficient by coefficient.
\end{proof}

\subsection{The machine}

Weighting the partial Ricci operators by products of Ricci eigenvalues, we obtain a machine of order $4$.

\begin{theorem}\label{thm:ricpar}
Let
\[
Q:=\Lambda\,\Id-\sum_a\Lambda_a\,\mathcal R_a,\qquad \mathfrak A:=Q\circ D\Sch(g_0):\Gamma(S_2)\to\Gamma(S_2).
\]
\begin{enumerate}[label=(\roman*)]
\item For every metric $g_0$ with parallel Ricci tensor, with no assumption on the Schouten tensor or on the $\lambda_a$, $\mathfrak A\circ L=0$ and $\mathfrak A(fg_0)=0$ for every function $f$. Thus $\mathfrak A^*:\Gamma(S_2)\to\TT$ is a TT machine of order $4$.
\item If $\sum_a\Lambda_a|\xi_a|^2\neq0$ for every $\xi\neq0$, that is, if $\ell=1$, or if $\ell\ge2$ and all $\lambda_a$ are non-zero and of the same sign, then $\sigma_\xi(\mathfrak A)$ is injective on $\TT_\xi$ for every $\xi\neq0$. On a compact $M$ without boundary,
\[
\TT=\Ima\mathfrak A^*\oplus\big(\ker\mathfrak A\cap\TT\big),
\]
the second summand being finite dimensional.
\end{enumerate}
\end{theorem}

\begin{proof}
(i) By naturality of $\Sch$, $D\Sch(g_0)(L_Xg_0)=L_X\Sch_0=\sum_b\mu_bL_Xg_b$, the $\mu_b$ being constant. Lemma~\ref{lem:cle} gives
\[
\mathfrak A(L_Xg_0)=\Lambda\sum_b\mu_bL_Xg_b-\sum_a\Lambda_a\mu_a\lambda_aL_Xg_a=0,
\]
since $\Lambda_a\lambda_a=\Lambda$. Next, $D\Sch(g_0)(fg_0)=-\frac{n-2}{2}\nabla\nabla f$ \cite[Lemme 7]{D}, and $\nabla\nabla f=\frac12L_{\nabla f}g_0=\frac12\sum_bL_{\nabla f}g_b$, so similarly
\[
Q(\nabla\nabla f)=\tfrac12\Big(\Lambda\sum_bL_{\nabla f}g_b-\sum_a\Lambda_a\lambda_aL_{\nabla f}g_a\Big)=0 .
\]
Thus $\mathfrak A\circ L=0$ and $\mathfrak A(fg_0)=0$, and the proof of Proposition~\ref{prop:principe} shows that $\mathfrak A^*$ is TT-valued.

(ii) Let $\xi\neq0$ and $k\in\TT_\xi$. By Lemma~\ref{lem:Ea}(b) and $\sum_a\mathcal R_a=D\Ric$, $\sigma_\xi(D\Ric)k=\frac12|\xi|^2k$ and $\sigma_\xi(DR)k=0$, so $\sigma_\xi(D\Sch(g_0))k=\frac12|\xi|^2k$. The term $\Lambda\,\Id$ has order $0$, and
\[
\sigma_\xi(\mathfrak A)k=-\tfrac14|\xi|^2\Big(c\,k+\tau\,\xi\otimes\xi-\xi\odot v\Big),\qquad c=\sum_a\Lambda_as_a,\quad \tau=\sum_a\Lambda_a\tau_a,\quad v=\sum_a\Lambda_a\eta_a .
\]
Assume $c\neq0$ and $\sigma_\xi(\mathfrak A)k=0$. Then $k=c^{-1}(\xi\odot v-\tau\,\xi\otimes\xi)$, and $k(\xi^\sharp)=0$ reads $|\xi|^2v+(\ip v\xi-\tau|\xi|^2)\xi=0$. Hence $v=\beta\xi$ with $2\beta=\tau$, and $k=c^{-1}(2\beta-\tau)\,\xi\otimes\xi=0$. If all $\Lambda_a$ are non-zero and of the same sign, $c\neq0$ for every $\xi\neq0$; for $\ell\ge2$ this is equivalent to the $\lambda_a$ being non-zero and of the same sign, since $\Lambda_a=\Lambda/\lambda_a$ and at most one $\lambda_a$ vanishes. Theorem~\ref{thm:fredholm} applies, its proof only using these two identities.
\end{proof}

\begin{remark}\label{rem:ricpar}
\leavevmode
\begin{enumerate}[label=(\alph*)]
\item \emph{Converse for the symbol.} If $c(\xi)=0$ for some $\xi\neq0$, every $k\in\TT_\xi$ with vanishing diagonal blocks and $\eta_a(k)=0$ for all $a$ lies in the kernel of the symbol, for instance $k=\zeta\odot\tilde\zeta$ with $\zeta\in E_a^*$, $\tilde\zeta\in E_b^*$, $a\neq b$, $\zeta\perp\xi_a$, $\tilde\zeta\perp\xi_b$; such a non-zero $k$ exists as soon as $n_a,n_b\ge2$. The sign condition is thus necessary for the injectivity of $\sigma_\xi(\mathfrak A)$ on $\TT_\xi$, that is, for the method of proof of (ii). We do not claim that it is necessary for the conclusion of (ii), namely that $\Ima\mathfrak A^*$ has finite codimension in $\TT$; in any case, the order-$5$ machine of Section~\ref{sec:ordre5} has this property without any sign condition.
\item \emph{The Einstein case ($\ell=1$).} Then $Q=\lambda-D\Ric(g_0)=-D\Ricr(g_0)-\frac1n\,DR(g_0)(\cdot)\,g_0$, so for trace-free $\Psi$, $\mathfrak A^*\Psi=-D\Sch(g_0)^*D\Ricr(g_0)^*\Psi=-\frac14P\Psi$. The machine of the theorem is thus the direct generalisation of that of \cite{D}, valid for every $\lambda$, including $\lambda=0$.
\item \emph{Degenerate Schouten tensor.} The theorem involves neither $\Sch_0^{-1}$ nor spectral projections of a perturbed metric. For instance, $S^2(r_1)\times S^2(r_2)$ is covered for all radii, including the ratios $r_1/r_2=\sqrt2^{\pm1}$ where the Schouten tensor is degenerate, as are $S^{p_1}(r_1)\times\dots\times S^{p_{\ell}}(r_{\ell})$ ($p_a\ge2$), $H^{p_1}\times\dots\times H^{p_{\ell}}$ and their compact quotients, for all radii.
\item When some $\lambda_a$ vanishes, or the signs are mixed, $\mathfrak A^*$ is still a TT machine, but control of the image by Theorem~\ref{thm:fredholm} is lost. The order-$5$ machine of Section~\ref{sec:ordre5} covers all these cases.
\end{enumerate}
\end{remark}

\subsection{Origin of the construction}

When $\Sch_0$ is invertible, conformal directions become Lie derivatives of $\Sch_0$.

\begin{lemma}\label{lem:Xf}
If $\Sch_0$ is invertible and $X_f:=-\frac{n-2}{4}\Sch_0^{-1}(df)$, then $D\Sch(g_0)(fg_0)=L_{X_f}\Sch_0$.
\end{lemma}

\begin{proof}
By \cite[Lemme 7]{D}, $D\Sch(g_0)(fg_0)=-\frac{n-2}{2}\nabla\nabla f$, and, $\Sch_0$ being parallel, $(L_X\Sch_0)_{ij}=(\Sch_0)_{ik}\nabla_jX^k+(\Sch_0)_{jk}\nabla_iX^k$, which equals $-\frac{n-2}{2}\nabla_i\nabla_jf$ for $X=X_f$.
\end{proof}

The machine of Theorem~\ref{thm:ricpar} was found by regarding the Schouten tensor as a metric. If $\Sch_0$ is invertible and all $\lambda_a$ are non-zero, the operator $\mathcal N(G):=G-\sum_a\frac{\mu_a}{\lambda_a}\,\mathrm{Sym}\,C_a(G;\Pi_a)$ vanishes and is infinitesimally natural at $G_0=\Sch_0$ (proof of Lemma~\ref{lem:cle}), with $D\mathcal N(G_0)=\Lambda^{-1}Q$; by Lemma~\ref{lem:Xf} and Proposition~\ref{prop:relache}, $N:=\mathcal N\circ\Sch$ gives $DN(g_0)=\Lambda^{-1}\mathfrak A$. The $\mu_a$ disappear from the final formula, which explains why Theorem~\ref{thm:ricpar} needs no assumption on the Schouten tensor.

\section{An order-5 machine for metrics with parallel Ricci tensor}\label{sec:ordre5}

We keep the setting and notation of Section~\ref{sec:ricpar}. The machine of this section has the advantage of being written only in terms of standard operators ($D\Ric$, $D\Sch$, $\nabla$), without the projections $\Pi_a$, and of producing all TT tensors up to a finite-dimensional space with no condition on the $\lambda_a$ or on the Schouten tensor.

\subsection{Reduction to a product}\label{ss:reduction}

The statements of this section are of two kinds. The identities $\mathfrak A\circ L=0$ and $\mathfrak A(fg_0)=0$ for a differential operator $\mathfrak A$ are \emph{local}: it suffices to check them on $Lw$ and $fg_0$ for $w$, $f$ supported in an arbitrarily small open set. The symbol computation is \emph{pointwise}. Only Theorem~\ref{thm:fredholm} is global, and it uses only these two facts. By de Rham's theorem \cite{dR}, \cite[Ch.~IV]{KN}, every point of $M$ has a neighbourhood $U$ isometric to a Riemannian product
\[
(U,g_0)\simeq(U_1\times\dots\times U_{\ell},\ g_1\oplus\dots\oplus g_{\ell}),\qquad\Ric(g_a)=\lambda_ag_a ,
\]
where the $E_a$ are tangent to the factors and the $\Pi_a$ are the corresponding projections. \emph{In what follows we may therefore assume that $g_0$ is such a product}; the objects $\Pi_a$, $g_a$, $\mathcal R_a$, $F_a$ below are then those of the product, and they patch together globally since they only depend on the decomposition of $TM$ into eigenspaces of $\Ric_0$.

\subsection{The key identity}

For $h\in S_2$ let $\delta\Gamma_a(h)^m_{jk}:=(\Pi_a)^m{}_p\,\delta\Gamma(h)^p_{jk}$, so that $\sum_a\delta\Gamma_a=\delta\Gamma$, and, with the notation of \eqref{eq:DnablaT},
\[
F_a(h):=\nabla\big(\mathcal R_a(h)\big)-\lambda_a\,\delta\Gamma_a(h)\cdot g_a .
\]
Since $\Ric_0=\sum_a\lambda_ag_a$ and $\delta\Gamma_b(h)\cdot g_a=0$ for $b\neq a$, we have $\delta\Gamma(h)\cdot\Ric_0=\sum_a\lambda_a\,\delta\Gamma_a(h)\cdot g_a$, hence, by \eqref{eq:DnablaT},
\begin{equation}\label{eq:sommeFa}
\sum_aF_a=D(\nabla\Ric)(g_0),\qquad D(\nabla\Ric)(g_0)h=\nabla\big(D\Ric(g_0)h\big)-\delta\Gamma(h)\cdot\Ric_0 .
\end{equation}

The operators $F_a$ arise as components of the linearisation of $\nabla\Ric$ at the metrics obtained by rescaling the factors.

\begin{lemma}\label{lem:lineaireinv}
For $c=(c_1,\dots,c_{\ell})$ with all $c_a\neq0$, the tensor $G_c:=c_1g_1\oplus\dots\oplus c_{\ell}g_{\ell}$ is non-degenerate, $\nabla$-parallel, with Levi-Civita connection $\nabla$, and
\[
D(\nabla\Ric)(G_c)=\sum_ac_a^{-1}F_a .
\]
\end{lemma}

\begin{proof}
Since $\nabla G_c=0$, $\nabla$ is torsion free and $G_c$-metric, hence it is the Levi-Civita connection of $G_c$; in particular $\Ric(G_c)=\Ric_0$. By \eqref{eq:DnablaT}, $D(\nabla\Ric)(G_c)h=\nabla\big(D\Ric(G_c)h\big)-\delta\Gamma_{G_c}(h)\cdot\Ric_0$, and $\delta\Gamma_{G_c}(h)^m_{jk}=\frac12(G_c^{-1})^{mp}(\nabla_jh_{pk}+\nabla_kh_{pj}-\nabla_ph_{jk})$ depends on $G_c$ only linearly through $G_c^{-1}=\sum_ac_a^{-1}g_a^{-1}$, so $\delta\Gamma_{G_c}=\sum_ac_a^{-1}\delta\Gamma_a$. Likewise $D\Ric(G_c)=\sum_ac_a^{-1}\mathcal R_a$, since $(\Pi_a)^k{}_mG_c^{mp}=c_a^{-1}g_a^{kp}$. Finally $\delta\Gamma_a(h)\cdot\Ric_0=\lambda_a\,\delta\Gamma_a(h)\cdot g_a$.
\end{proof}

When the Schouten tensor is invertible, the heuristic of Section~\ref{sec:ricpar} gives a first machine of order $5$ directly.

\begin{proposition}[Invertible Schouten tensor]\label{prop:schinv}
Assume all $\mu_a$ are non-zero, and let $G_0:=\Sch_0=\mu_1g_1\oplus\dots\oplus\mu_{\ell}g_{\ell}$. Then
\[
\mathfrak A_{\Sch}:=D(\nabla\Ric)(G_0)\circ D\Sch(g_0)=\Big(\sum_a\mu_a^{-1}F_a\Big)\circ D\Sch(g_0)
\]
satisfies $\mathfrak A_{\Sch}\circ L=0$ and $\mathfrak A_{\Sch}(fg_0)=0$. Its symbol on $\TT_\xi$ is $\frac i4|\xi|^2\,\xi\otimes\big(c\,k+\tau\,\xi\otimes\xi-\xi\odot v\big)$ with $c=\sum_a\mu_a^{-1}|\xi_a|^2$, which is injective if $\Sch_0$ is definite.
\end{proposition}

\begin{proof}
The (possibly indefinite) metric $G_0$ is parallel, so, by Lemma~\ref{lem:lineaireinv}, $\nabla\Ric(G_0)=\nabla\Ric_0=0$. By naturality of $G\mapsto\nabla^G\Ric(G)$, $D(\nabla\Ric)(G_0)(L_YG_0)=L_Y\big(\nabla\Ric(G_0)\big)=0$ for every $Y$. Now $D\Sch(g_0)(L_Xg_0)=L_X\Sch_0=L_XG_0$ by naturality of $\Sch$, and $D\Sch(g_0)(fg_0)=L_{X_f}G_0$ by Lemma~\ref{lem:Xf}. Hence $\mathfrak A_{\Sch}\circ L=0$ and $\mathfrak A_{\Sch}(fg_0)=0$. The formula follows from Lemma~\ref{lem:lineaireinv} with $c_a=\mu_a$, and the symbol is computed as in the proof of Theorem~\ref{thm:ricpar}.
\end{proof}

\begin{remark}\label{rem:continuite}
The weights $\mu_a^{-1}$ of Proposition~\ref{prop:schinv} blow up when some $\mu_a$ vanishes. A continuity argument in the scales of the factors only yields the regularised machine with weights $\prod_{d\neq a}\mu_d$, which vanishes as soon as two $\mu_a$ vanish. The weight-$1$ machine of Theorem~\ref{thm:ordre5} requires instead the identities $F_b(L_Xg_b)=0$ factor by factor; this is the content of Lemma~\ref{lem:cle5}, where the projection $\Pi_a$ is used as background structure.
\end{remark}

The weights $\mu_a^{-1}$ can be removed thanks to the following identity, the analogue of Lemma~\ref{lem:cle} one derivative higher.

\begin{lemma}\label{lem:cle5}
For every vector field $Y$ and all $a,b$, $F_a(L_Yg_b)=0$. In particular, for every $b$,
\[
D(\nabla\Ric)(g_0)\big(L_Yg_b\big)=0 .
\]
\end{lemma}

\begin{proof}
By \S\ref{ss:reduction}, we work on a product. We take $G_c=c_1g_1\oplus\dots\oplus c_{\ell}g_{\ell}$, $c_a>0$, and $C_a(G;\Pi_a)$ as in the proof of Lemma~\ref{lem:cle}, and consider the tensor $\nabla^GC_a(G;\Pi_a)$, natural in $(G,\Pi_a)$. It vanishes at $G=G_c$, since $C_a(G_c;\Pi_a)=\lambda_ag_a$ is parallel.

\emph{Linearisation.} As in Lemma~\ref{lem:lineaireinv}, $\delta\Gamma_{G_c}(h)=\sum_bc_b^{-1}\delta\Gamma_b(h)$. By \eqref{eq:DnablaT} and the proof of Lemma~\ref{lem:cle},
\[
D_G\big(\nabla^GC_a\big)(G_c)h=\nabla\big(c_a^{-1}\mathcal R_a(h)\big)-\sum_bc_b^{-1}\delta\Gamma_b(h)\cdot(\lambda_ag_a)=c_a^{-1}F_a(h),
\]
since $\delta\Gamma_b(h)\cdot g_a=0$ for $b\neq a$.

\emph{Failure of naturality.} The derivative with respect to $\Pi_a$ in the direction $L_Y\Pi_a$ is
\[
\nabla^{G_c}\big(D_{\Pi}C_a(G_c;\Pi_a)(L_Y\Pi_a)\big).
\]
Now $D_{\Pi}C_a(G_c;\Pi_a)(L_Y\Pi_a)$ is the tensor field $(L_Y\Pi_a)^k{}_mR^m{}_{jkl}$, which vanishes identically by the proof of Lemma~\ref{lem:cle}; its covariant derivative therefore vanishes.

\emph{Conclusion.} Differentiating at $t=0$ the identity
\[
\nabla^{\varphi_t^*G_c}C_a(\varphi_t^*G_c;\varphi_t^*\Pi_a)=\varphi_t^*\big(\nabla^{G_c}C_a(G_c;\Pi_a)\big)=0
\]
gives $c_a^{-1}\sum_bc_b\,F_a(L_Yg_b)=0$ for all $c\in(0,\infty)^\ell$, hence $F_a(L_Yg_b)=0$ for all $a,b$. The second statement follows from \eqref{eq:sommeFa}.
\end{proof}

\subsection{The machine}

Composing $D(\nabla\Ric)(g_0)$ with $D\Sch(g_0)$ then gives a machine without any condition.

\begin{theorem}\label{thm:ordre5}
Let $g_0$ be a metric with parallel Ricci tensor, and
\begin{gather*}
\mathfrak A_5:=D(\nabla\Ric)(g_0)\circ D\Sch(g_0):\Gamma(S_2)\to\Gamma(T^0_3),\\
\mathfrak A_5h=\nabla\big(D\Ric(g_0)q\big)-\delta\Gamma(q)\cdot\Ric_0,\qquad q=D\Sch(g_0)h .
\end{gather*}
\begin{enumerate}[label=(\roman*)]
\item $\mathfrak A_5\circ L=0$ and $\mathfrak A_5(fg_0)=0$ for every function $f$, with no assumption on the Schouten tensor or on the $\lambda_a$. Thus $\mathfrak A_5^*=D\Sch(g_0)^*\circ D(\nabla\Ric)(g_0)^*:\Gamma(T^0_3)\to\TT$ is a TT machine of order $5$.
\item For $\xi\neq0$ and $k\in\TT_\xi$, $\sigma_\xi(\mathfrak A_5)k=\frac i4|\xi|^4\,\xi\otimes k$, which is injective on $\TT_\xi$. Consequently, if $M$ is compact without boundary,
\[
\TT=\Ima\mathfrak A_5^*\oplus\big(\ker\mathfrak A_5\cap\TT\big),
\]
the second summand being finite dimensional, for \emph{every} metric with parallel Ricci tensor.
\end{enumerate}
\end{theorem}

\begin{proof}
(i) By naturality, $D\Sch(g_0)(L_Xg_0)=L_X\Sch_0=\sum_b\mu_bL_Xg_b$, and $D\Sch(g_0)(fg_0)=-\frac{n-2}{2}\nabla\nabla f=-\frac{n-2}{4}\sum_bL_{\nabla f}g_b$. Lemma~\ref{lem:cle5} gives $\mathfrak A_5(L_Xg_0)=0$ and $\mathfrak A_5(fg_0)=0$.

(ii) The term $\delta\Gamma(q)\cdot\Ric_0$ has order $1$ in $q$, so $\sigma_\xi(D(\nabla\Ric)(g_0))q=i\,\xi\otimes\sigma_\xi(D\Ric(g_0))q$. For $k\in\TT_\xi$, $\sigma_\xi(D\Sch(g_0))k=\frac12|\xi|^2k$ (proof of Theorem~\ref{thm:ricpar}) and $\sigma_\xi(D\Ric(g_0))k=\frac12|\xi|^2k$, whence the formula. We conclude with Theorem~\ref{thm:fredholm}, whose proof only uses these two identities.
\end{proof}

\begin{remark}
\leavevmode
\begin{enumerate}[label=(\alph*)]
\item \emph{The Einstein case.} If $\ell=1$, $\delta\Gamma(q)\cdot\Ric_0=\lambda\,\delta\Gamma(q)\cdot g_0=\lambda\nabla q$, so $\mathfrak A_5=\nabla\circ\big(D\Ric(g_0)-\lambda\big)\circ D\Sch(g_0)=-\nabla\circ\mathfrak A$, with $\mathfrak A$ from Theorem~\ref{thm:ricpar}. Thus $\mathfrak A_5^*\Psi=-\mathfrak A^*\big(\mathrm{Sym}\,\nabla^*\Psi\big)$, which equals $\frac14P\big(\mathrm{Sym}\,\nabla^*\Psi\big)$ when the latter is trace free (Remark~\ref{rem:ricpar}(b)). For $\ell\ge2$, $\mathfrak A_5$ is not of the form $\nabla\circ\mathfrak A$.
\item \emph{Comparison with order $4$.} The machine $\mathfrak A$ of Theorem~\ref{thm:ricpar} has lower order, but its ellipticity requires the $\lambda_a$ to be non-zero and of the same sign. The machine $\mathfrak A_5$ needs no condition, at the price of one more derivative and of the source bundle $T^0_3$.
\item \emph{Variants.} By Lemma~\ref{lem:cle5}, $\big(\sum_aw_aF_a\big)\circ D\Sch(g_0)$ is a TT machine for all constants $w_a$. Its symbol on $\TT_\xi$ is $\frac i4|\xi|^2\,\xi\otimes\big(c\,k+\tau\,\xi\otimes\xi-\xi\odot v\big)$ with $c=\sum_aw_a|\xi_a|^2$, as in the proof of Theorem~\ref{thm:ricpar}: it is injective as soon as the $w_a$ are all non-zero and of the same sign. The choice $w_a=1$ is the only one not involving the projections.
\end{enumerate}
\end{remark}

\section{Bach-flat products of Einstein metrics and cones}\label{sec:produits}

\subsection{The Bach tensor of a product of Einstein metrics}

Let $g_0=g_1\oplus g_2$ on $M_1\times M_2$, with $\dim M_a=n_a$, $n=n_1+n_2\ge3$ and $\Ric(g_a)=\lambda_ag_a$. Let $\mu_a=\lambda_a-\frac{R}{2(n-1)}$, with $R=n_1\lambda_1+n_2\lambda_2$, and $\alpha_a=\frac{\mu_a}{n-2}$, so that $A_0=\alpha_1g_1+\alpha_2g_2$.

The Bach tensor of such a product can be computed explicitly.

\begin{proposition}\label{prop:Bproduit}
$C(g_0)=0$ and
\[
B(g_0)=(\alpha_1^2-\alpha_2^2)\,(n_2\,g_1-n_1\,g_2).
\]
In particular $B(g_0)=0$ if and only if $\lambda_1=\lambda_2$ or
\begin{equation}\label{eq:equilibre}
(n_2-1)\lambda_1+(n_1-1)\lambda_2=0 .
\end{equation}
\end{proposition}

\begin{proof}
$\Ric(g_0)$ is parallel, hence so is $A_0$, and $C(g_0)=0$. Consequently $B(g_0)=\mathring WA_0=\alpha_1T+\alpha_2\tilde T$, with $T_{ij}=g_1^{kl}W_{ikjl}$ and $\tilde T_{ij}=g_2^{kl}W_{ikjl}$. Since $g^{kl}W_{ikjl}=0$, $\tilde T=-T$, so $B(g_0)=(\alpha_1-\alpha_2)T$.

We compute $T$ on the $M_2$ block, i.e.\ for $i,j$ tangent to $M_2$ and $k,l$ tangent to $M_1$. The product being Riemannian, $R_{ikjl}=0$. In $(A_0\KN g_0)_{ikjl}=A_{ij}g_{kl}+A_{kl}g_{ij}-A_{il}g_{kj}-A_{kj}g_{il}$ the last two terms vanish, being mixed. Contracting with $g_1^{kl}$ gives $n_1\alpha_2g_2+n_1\alpha_1g_2$, so $T|_{M_2}=-n_1(\alpha_1+\alpha_2)g_2$.

The mixed blocks of $T$ vanish, since all mixed components involved vanish. On the $M_1$ block, $g_1^{kl}R_{ikjl}=\lambda_1(g_1)_{ij}$ and $g_1^{kl}(A_0\KN g_0)_{ikjl}=2(n_1-1)\alpha_1(g_1)_{ij}$, so $T|_{M_1}=c_1g_1$ with $c_1=\lambda_1-2(n_1-1)\alpha_1$. Since $\Tr_{g_0}T=g_1^{kl}g^{ij}W_{ikjl}=0$, $n_1c_1=n_1n_2(\alpha_1+\alpha_2)$, i.e.\ $c_1=n_2(\alpha_1+\alpha_2)$. This gives the formula.

Finally $\alpha_1+\alpha_2=0\iff\lambda_1+\lambda_2=\frac{n_1\lambda_1+n_2\lambda_2}{n-1}\iff(n_2-1)\lambda_1+(n_1-1)\lambda_2=0$, and $\alpha_1=\alpha_2\iff\lambda_1=\lambda_2$.
\end{proof}

For $n=4$ and $n_1=n_2=2$, this recovers the known fact that a product of two surfaces of constant curvatures $\kappa_1,\kappa_2$ is Bach-flat if and only if $\kappa_1=\pm\kappa_2$; for $n\ge5$, with the normalisation of Definition~\ref{def:Bach}, we have not found Proposition~\ref{prop:Bproduit} in the literature.

\begin{remark}
Comparing the two expressions of $c_1$ gives $\lambda_1=2(n_1-1)\alpha_1+n_2(\alpha_1+\alpha_2)$, which can also be checked directly from the definitions. For $S^2\times S^2$ one finds $B\propto(\lambda_1^2-\lambda_2^2)(g_1-g_2)$, which vanishes only for equal radii.
\end{remark}

\begin{remark}[Products of $\ell$ factors]\label{rem:kfacteurs}
For $g_0=g_1\oplus\dots\oplus g_{\ell}$, with the same notation, the same computation gives $T_a:=g_a^{kl}W_{ikjl}$ with $T_a|_{E_b}=-n_a(\alpha_a+\alpha_b)g_b$ for $b\neq a$, vanishing mixed blocks, and $T_a|_{E_a}=c_{aa}g_a$, where the trace condition gives $c_{aa}=\sum_{b\neq a}n_b(\alpha_a+\alpha_b)$. Hence
\[
B(g_0)=\sum_a\alpha_aT_a=\sum_b\big(n\,\alpha_b^2-\textstyle\sum_an_a\alpha_a^2\big)\,g_b .
\]
Thus $B(g_0)=0$ if and only if all $|\mu_a|$ are equal. The $\lambda_a$ then take at most two values $\lambda_\pm$; grouping the factors with the same constant into two Einstein factors of dimensions $m_\pm$, the condition becomes $(m_--1)\lambda_++(m_+-1)\lambda_-=0$. The Bach tensor therefore gives nothing new beyond two Ricci eigenvalues; an example is $\mathbb R\times N_1\times\dots\times N_j$ with Einstein $N_i$ of the same constant.
\end{remark}

For products satisfying \eqref{eq:equilibre}, the Bach machine is therefore available.

\begin{theorem}\label{thm:produits}
If $g_0=g_1\oplus g_2$ is a product of Einstein metrics satisfying \eqref{eq:equilibre}, then $DB(g_0)^*$ is a TT machine of order $4$ for $g_0$. If $M$ is compact without boundary,
\[
\TT=\Ima DB(g_0)^*\oplus\big(\ker DB(g_0)\cap\TT\big),
\]
the second summand being finite dimensional.
\end{theorem}

\begin{proof}
By Proposition~\ref{prop:Bproduit}, $B(g_0)=0$ and $C(g_0)=0$; apply Theorem~\ref{thm:bachplat}(b).
\end{proof}

\begin{example}
\leavevmode
\begin{enumerate}[label=(\alph*)]
\item \emph{$\mathbb R\times N$ for every Einstein $N$}, since then $n_1=1$ and $\lambda_1=0$; likewise $S^1\times N$.
\item \emph{Factors with Einstein constants of opposite signs.} If $\lambda_1>0>\lambda_2$ and $n_1,n_2\ge2$, a unique ratio of scales between the factors achieves \eqref{eq:equilibre}, namely $\frac{\lambda_1}{n_1-1}=-\frac{\lambda_2}{n_2-1}$; for instance $\mathbb{CP}^2\times\Sigma$, with $\mathbb{CP}^2$ normalised by $\lambda_1=6$ and $\Sigma$ a hyperbolic surface of curvature $-2$. When both factors are space forms of opposite curvatures this is $S^p\times H^q$ with equal radii, which is conformally flat; otherwise these products are not conformally flat.
\end{enumerate}
\end{example}

\begin{remark}[Summary for products of two Einstein factors]
A product $g_1\oplus g_2$ of Einstein factors has parallel Ricci tensor, so it is always covered by Sections~\ref{sec:ricpar} and~\ref{sec:ordre5}. The order-$5$ machine of Theorem~\ref{thm:ordre5} exists and controls the image for all $\lambda_a$; the other machines only serve to lower the order.
\begin{center}
\begin{tabular}{lll}
\toprule
Einstein constants & Order $4$, image controlled & Order $5$\\
\midrule
$\lambda_1=\lambda_2$ & $P$ of \cite{D} & Theorem~\ref{thm:ordre5}\\
$\lambda_1\lambda_2>0$ & Theorem~\ref{thm:ricpar} & Theorem~\ref{thm:ordre5}\\
$(n_2-1)\lambda_1+(n_1-1)\lambda_2=0$ & $DB(g_0)^*$ (Theorem~\ref{thm:produits}) & Theorem~\ref{thm:ordre5}\\
other cases & no (TT only) & Theorem~\ref{thm:ordre5}\\
\bottomrule
\end{tabular}
\end{center}
\end{remark}

\subsection{Cones over an Einstein section}

Let $(N,\hat g)$ have dimension $n-1\ge2$ with $\Ric(\hat g)=(n-2)\mu\,\hat g$, and consider the cone $g_c=dr^2+r^2\hat g$ on $(0,\infty)\times N$.

A cone is conformal to a product with a line, which reduces it to the previous case.

\begin{theorem}\label{thm:cone}
With $t=\ln r$, $g_c=r^2g_{\mathrm{p}}$, where $g_{\mathrm p}=dt^2\oplus\hat g$ is a Bach- and Cotton-flat product $\mathbb R\times N$. For every $\mu\in\mathbb R$,
\[
\mathcal P_c\,\Psi:=r^{-(n-2)}\,DB(g_{\mathrm p})^*\Psi
\]
is a TT machine of order $4$ on the cone, defined globally on $(0,\infty)\times N$.
\end{theorem}

\begin{proof}
The decomposition $g_c=r^2(dt^2\oplus\hat g)$ is immediate. Proposition~\ref{prop:Bproduit} with $n_1=1$ and $\lambda_1=0$ gives $B(g_{\mathrm p})=0$ and $C(g_{\mathrm p})=0$. Theorem~\ref{thm:bachplat}(b) gives a TT machine on $g_{\mathrm p}$, and we conclude with Lemma~\ref{lem:transport} with $e^f=r^2$.
\end{proof}

When $\mu\ge0$, one may instead transport the machine $P$ of an Einstein metric.

\begin{proposition}\label{prop:conepos}
The cases $\mu\ge0$ also admit a realisation by transport of $P$.
\begin{enumerate}[label=(\alph*)]
\item \emph{$\mu>0$.} With $\rho=r^{\sqrt\mu}$, $g_c=\mu^{-1}r^{2-2\sqrt\mu}\,g_E$, where $g_E=d\rho^2+\rho^2(\mu\hat g)$ is the Ricci-flat cone over $(N,\mu\hat g)$. Hence
\[
\Psi\longmapsto\mu^{\frac{n-2}{2}}\,r^{(n-2)(\sqrt\mu-1)}\,P_{g_E}\Psi,\qquad P_{g_E}|_{\TT}=\Delta_L^2 ,
\]
is a global TT machine on the cone.
\item \emph{$\mu=0$.} $g_{\mathrm p}=dt^2\oplus\hat g$ is Ricci flat, and $\Psi\mapsto r^{-(n-2)}P_{g_{\mathrm p}}\Psi$ works.
\end{enumerate}
\end{proposition}

\begin{proof}
(a) $d\rho^2=\mu r^{2\sqrt\mu-2}dr^2$ and $\mu\rho^2\hat g=\mu r^{2\sqrt\mu-2}\,r^2\hat g$, so $g_E=\mu r^{2\sqrt\mu-2}g_c$. Since the Ricci tensor is scale invariant, $\Ric(\mu\hat g)=\Ric(\hat g)=(n-2)\mu\hat g=(n-2)(\mu\hat g)$: the section $(N,\mu\hat g)$, of dimension $n-1$, is Einstein with constant $n-2$, which is exactly the condition for its cone to be Ricci flat \cite[Ch.~9]{Be}. Hence $g_E$ is Ricci flat, and $r\mapsto\rho$ is a diffeomorphism of $(0,\infty)$. Lemma~\ref{lem:transport} with $e^f=\mu^{-1}r^{2-2\sqrt\mu}$ gives the weight $e^{-\frac{n-2}{2}f}=\mu^{\frac{n-2}{2}}r^{(n-2)(\sqrt\mu-1)}$. The formula $P|_{\TT}=\Delta_L^2$ is that of \cite{D} for $\lambda=0$.

(b) For $\mu=0$ the product $\mathbb R\times N$ is Ricci flat; apply Lemma~\ref{lem:transport} with $e^f=r^2$.
\end{proof}

\begin{remark}
\leavevmode
\begin{enumerate}[label=(\alph*)]
\item \emph{$\mu<0$.} A warped product $ds^2+w(s)^2\hat g$ is Einstein (see \cite[Ch.~9]{Be}, \cite[Ch.~7]{ON}) with constant $(n-1)\kappa$ if and only if, writing $\dot w=dw/ds$, $\ddot w=-\kappa w$ and $\dot w^2+\kappa w^2=\mu$. For $\mu<0$ this forces $\kappa<0$ and $w=\sqrt{|\mu|/|\kappa|}\cosh(\sqrt{|\kappa|}\,s)$ up to translation. With $dt=ds/w$, the variable $t$ then ranges over a bounded interval only. Thus, with a conformal factor depending only on $r$, the cone is conformally Einstein only on slabs $\{r_1<r<r_2\}$; a global machine is provided by Theorem~\ref{thm:cone}, or by conformal transport, with weight $r^{-(n-2)}$, of the order-$5$ machine of Theorem~\ref{thm:ordre5} for the product $\mathbb R\times N$, which has parallel Ricci tensor.
\item For $\mu>0$, the machines of Theorem~\ref{thm:cone} and Proposition~\ref{prop:conepos} are both valid. In dimension $4$, writing $g_E=e^{f}g_{\mathrm p}$ with $e^{f}=\mu r^{2\sqrt\mu}$, the conformal covariance of $B$ gives $DB(g_E)^*\Psi=e^{-f}DB(g_{\mathrm p})^*(e^{-f}\Psi)$, and, with $P_{g_E}=-4DB(g_E)$ (Theorem~\ref{thm:bach}), the machine of Proposition~\ref{prop:conepos} is $\Psi\mapsto-4\,\mathcal P_c(e^{-f}\Psi)$.
\item The same arguments apply to every warped product $I\times_wN$ with $N$ Einstein, which is conformal to a piece of $\mathbb R\times N$ via $dt=ds/w$.
\end{enumerate}
\end{remark}

\section{Machines for divergence-free symmetric tensors}\label{sec:divnulle}

If one only asks for $\Div h=0$, the conformal condition (ii) of Proposition~\ref{prop:principe} disappears. We denote by
\[
\Sigma_\xi:=\{k\in S_2:\ k(\xi^\sharp,\cdot)=0\}\supset\TT_\xi
\]
the space of divergence-free symbols at $\xi\neq0$.

\subsection{Principle}

For divergence-free tensors, diffeomorphism invariance alone is needed.

\begin{proposition}\label{prop:principediv}
Let $N$ be natural, or infinitesimally natural in the sense of Proposition~\ref{prop:relache}, with $N(g_0)=0$. Then $DN(g_0)\circ L=0$, and $DN(g_0)^*$ takes values in divergence-free symmetric tensors.
\end{proposition}

\begin{proof}
This is the first part of the proof of Proposition~\ref{prop:principe}, where only diffeomorphism invariance is used.
\end{proof}

The Fredholm statement keeps the same form, with $\Sigma_\xi$ in place of $\TT_\xi$.

\begin{theorem}\label{thm:fredholmdiv}
Assume $M$ is compact without boundary, $\mathfrak A:\Gamma(S_2)\to\Gamma(F)$ is a differential operator of order $m\ge1$ with $\mathfrak A\circ L=0$, and $\sigma_\xi(\mathfrak A)$ is injective on $\Sigma_\xi$ for every $\xi\neq0$. Then, in $C^\infty$,
\[
\ker\Div=\Ima\mathfrak A^*\oplus\big(\ker\mathfrak A\cap\ker\Div\big),
\]
the sum being $L^2$-orthogonal and the second summand finite dimensional.
\end{theorem}

\begin{proof}
Apply Lemma~\ref{lem:fredholmD} with $f=0$, $E=S_2$, $D=L^*=2\Div:\Gamma(S_2)\to\Gamma(T^*M)$ and $\mathcal P=\mathfrak A^*:\Gamma(F)\to\Gamma(S_2)$: then $D\mathcal P=(\mathfrak A\circ L)^*=0$, $\mathcal P^*=\mathfrak A$, and $\ker\sigma_\xi(L^*)=\Sigma_\xi$.
\end{proof}

\begin{remark}
\leavevmode
\begin{enumerate}[label=(\alph*)]
\item Every TT machine is a divergence-free machine; the converse is false.
\item Lemma~\ref{lem:transport} has no analogue: $\Div_{e^f\tilde g}$ involves $\Tr_{\tilde g}h\,df$, so the divergence-free condition is preserved by conformal changes only for trace-free tensors. Only homotheties are allowed.
\end{enumerate}
\end{remark}

\subsection{Examples}

We begin with constant curvature, where the machine is related to elasticity.

\begin{proposition}[Constant curvature, elasticity]\label{prop:divcc}
Let $g_0$ have constant sectional curvature $\kappa$, and $N(g):=\Rm(g)-\frac\kappa2\,g\KN g$, with values in the bundle $\mathcal C$ of algebraic curvature tensors. Then $DN(g_0)^*:\Gamma(\mathcal C)\to\ker\Div$ is a machine of order $2$, and Theorem~\ref{thm:fredholmdiv} applies.
\begin{enumerate}[label=(\alph*)]
\item If $g_0$ is flat, in Cartesian coordinates, $\big(DN(g_0)^*Z\big)_{ad}=2\,\partial^b\partial^cZ_{abcd}$.
\item In $\mathbb R^3$, every $Z\in\mathcal C$ can be written $Z_{abcd}=\varepsilon_{abi}\varepsilon_{cdj}\Phi_{ij}$ with $\Phi$ symmetric, and one obtains, up to a constant, the Beltrami representation $\sigma=\operatorname{inc}\Phi$ \cite{Gu}, where $(\operatorname{inc}\Phi)_{ad}=\varepsilon_{abi}\varepsilon_{dcj}\partial_b\partial_c\Phi_{ij}$.
\end{enumerate}
\end{proposition}

\begin{proof}
$N$ is natural (the constant $\kappa$ being fixed) and vanishes at $g_0$; Proposition~\ref{prop:principediv} applies. For the symbol, Lemma~\ref{lem:DRm} gives $\sigma_\xi(DN(g_0))k=\frac12\,\xi\xi\KN k$. With $|\xi|=1$, $\xi=e_1$ and $k\in\Sigma_\xi$, hence supported on $e_1^\perp$, the $(1,a,1,b)$ component of $\xi\xi\KN k$ is $k_{ab}$; the symbol is therefore injective on $\Sigma_\xi$. (a) In the flat case only the second-order terms of Lemma~\ref{lem:DRm} remain, and the computation of the second-order terms in the proof of Proposition~\ref{prop:DWstar} gives the formula; one checks directly that $\partial^a\partial^b\partial^cZ_{abcd}=0$ by antisymmetry in $(a,b)$. (b) In dimension $3$, $\mathcal C$ has rank $6$ and $\Phi\mapsto\varepsilon\varepsilon\Phi$ is an isomorphism onto it; substitute in (a).
\end{proof}

On Einstein metrics, the linearised Einstein operator is itself a divergence-free machine.

\begin{proposition}[Einstein metrics]\label{prop:divein}
Let $g_0$ be Einstein, $\Ric(g_0)=\lambda g_0$, $n\ge3$, and $\operatorname{Ein}(g):=\Ric(g)-\frac{R(g)}2g+\Lambda g$ with $\Lambda=\frac{(n-2)\lambda}{2}$, so that $\operatorname{Ein}(g_0)=0$. Then the linearised Einstein operator $D\operatorname{Ein}(g_0)$ (see \cite{FM}) is self-adjoint, of order $2$, with values in $\ker\Div$, and, if $M$ is compact,
\[
\ker\Div=\Ima D\operatorname{Ein}(g_0)\oplus\big(\ker D\operatorname{Ein}(g_0)\cap\ker\Div\big).
\]
\end{proposition}

\begin{proof}
This is essentially \cite[Proposition 8]{D}; we recall the short argument. $\operatorname{Ein}$ is natural and, up to sign, the $L^2$-gradient of $\int(R-2\Lambda)\,dV$ \cite[Ch.~4]{Be}; at its critical point $g_0$, $D\operatorname{Ein}(g_0)$ is therefore the Hessian up to sign, hence self-adjoint. Proposition~\ref{prop:principediv} gives $\Div D\operatorname{Ein}(g_0)=0$. For the symbol, for $k\in\Sigma_\xi$ and $\tau=\Tr k$, $\sigma_\xi(D\Ric)k=\frac12(|\xi|^2k+\tau\,\xi\otimes\xi)$ and $\sigma_\xi(DR)k=|\xi|^2\tau$, whence
\[
\sigma_\xi(D\operatorname{Ein})k=\tfrac12\big(|\xi|^2k+\tau\,\xi\otimes\xi-|\xi|^2\tau\,g_0\big).
\]
If it vanishes, taking the trace gives $(n-2)\tau=0$, so $\tau=0$ and $k=0$. We conclude with Theorem~\ref{thm:fredholmdiv}.
\end{proof}

The machines of Sections~\ref{sec:ricpar} and~\ref{sec:ordre5} have divergence-free counterparts of lower order.

\begin{proposition}[Parallel Ricci tensor]\label{prop:divricpar}
In the setting of Section~\ref{sec:ricpar}:
\begin{enumerate}[label=(\alph*)]
\item $Q^*$, with $Q=\Lambda\,\Id-\sum_a\Lambda_a\mathcal R_a$, is a divergence-free machine of order $2$; Theorem~\ref{thm:fredholmdiv} applies if $\ell=1$, or if the $\lambda_a$ are non-zero and of the same sign.
\item $D(\nabla\Ric)(g_0)^*:\Gamma(T^0_3)\to\ker\Div$ is a machine of order $3$, and Theorem~\ref{thm:fredholmdiv} applies with no condition.
\end{enumerate}
\end{proposition}

\begin{proof}
(a) By Lemma~\ref{lem:cle}, $Q(L_Xg_0)=\Lambda\sum_bL_Xg_b-\sum_a\Lambda_a\lambda_aL_Xg_a=0$. The symbol computation in the proof of Theorem~\ref{thm:ricpar}(ii) only uses $k(\xi^\sharp)=0$, hence holds on $\Sigma_\xi$. (b) By Lemma~\ref{lem:cle5}, $D(\nabla\Ric)(g_0)(L_Xg_0)=\sum_bD(\nabla\Ric)(g_0)(L_Xg_b)=0$. The symbol is $k\mapsto i\,\xi\otimes\sigma_\xi(D\Ric)k=\frac i2\,\xi\otimes(|\xi|^2k+\tau\,\xi\otimes\xi)$ on $\Sigma_\xi$; if it vanishes, $k=-\tau|\xi|^{-2}\xi\otimes\xi$, and $k(\xi^\sharp)=0$ forces $\tau=0$, hence $k=0$.
\end{proof}

On locally symmetric metrics, the covariant derivative of the curvature tensor gives a machine with no condition.

\begin{proposition}[Locally symmetric metrics]\label{prop:divlocsym}
If $\nabla\Rm(g_0)=0$, then $D(\nabla\Rm)(g_0)^*:\Gamma(T^0_5)\to\ker\Div$ is a machine of order $3$, and Theorem~\ref{thm:fredholmdiv} applies with no condition.
\end{proposition}

\begin{proof}
$N=\nabla\Rm$ is natural and vanishes at $g_0$. The symbol is $k\mapsto i\,\xi\otimes\frac12\,\xi\xi\KN k$, injective on $\Sigma_\xi$ by the proof of Proposition~\ref{prop:divcc}.
\end{proof}

\begin{remark}[Stress functions]\label{rem:stress}
In dimension $3$, an Einstein metric has constant curvature, where the Calabi complex already gives complete stress-function representations \cite{Ca,Kh,KL1}. The new three-dimensional cases covered by Propositions~\ref{prop:divricpar} and~\ref{prop:divlocsym} are the metrics locally isometric to $S^2\times\mathbb R$ or $H^2\times\mathbb R$: on their compact quotients, every divergence-free symmetric tensor is, up to a finite-dimensional space, of the form $D(\nabla\Rm)(g_0)^*\Psi$. Manifolds with boundary, which are needed in elasticity, are not treated here; see \cite{KL1,KL2}.
\end{remark}

\begin{remark}
For a generic metric, no explicit natural tensor vanishing at $g_0$ is known. An operator $\mathfrak A$ with $\mathfrak A\circ L=0$ exists (the compatibility operator of the Killing operator; in constant curvature it is the first operator of the Calabi complex \cite{Ca,Kh}, and the locally symmetric case is studied by Gasqui--Goldschmidt \cite{GG}), but it is in general neither natural nor explicit. A notable exception is the Kerr spacetime, for which a complete compatibility complex of the Killing operator has been constructed explicitly by Aksteiner, Andersson, B\"ackdahl, Khavkine and Whiting \cite{AABKW}: its first operator $\mathfrak A$ satisfies $\mathfrak A\circ L=0$, so that $\mathfrak A^*$ is a divergence-free machine on Kerr in the Lorentzian sense of Remark~\ref{rem:lorentz}. As for TT tensors, the approach is then non-local: $h\mapsto h-L(\Div L)^{-1}\Div h$ on a compact manifold without Killing fields.
\end{remark}

\section{Machines for tensors with vanishing weighted divergence; Ricci solitons}\label{sec:divf}

Let $f$ be a fixed smooth function and $d\mu_f=e^{-f}dV_{g_0}$. The \emph{weighted divergence} is
\[
(\Div_fh)_i:=-e^{f}\nabla^j\big(e^{-f}h_{ij}\big)=-\nabla^jh_{ij}+h_{ij}\nabla^jf ,
\qquad\text{that is}\qquad \Div_fh=\Div h+h(\nabla f,\cdot).
\]
Up to a factor, it is the adjoint of $L$ in $L^2(d\mu_f)$: $\int\ip{h}{LX}\,d\mu_f=2\int\ip{\Div_fh}{X}\,d\mu_f$. Tensors with vanishing weighted divergence appear in the study of gradient Ricci solitons
\[
\Ric+\Hess f=\lambda g .
\]
For a compact shrinking soliton, normalised as in \cite{CHI,CZ} by $\Ric+\Hess f=\frac{1}{2\tau}g$ (the instability criterion below is invariant under scaling), the second variation of Perelman's entropy $\nu$ \cite{Pe} is read, modulo diffeomorphisms, on $\ker\Div_f$: on this subspace, the stability operator is $\frac12\Delta_f+\Rm(\cdot)$, up to a scalar rank-one term \cite{CHI,CZ}. Here $\Delta_f=-\nabla^*\nabla-\nabla_{\nabla f}$ and $\Rm(h)_{bd}=R_{abcd}h^{ac}$. Instability of some Kähler--Ricci solitons is proved by Hall and Murphy \cite{HM}.

Only non-trivial potentials matter here. Indeed, if $f$ is constant, $\Div_f=\Div$ and the soliton is Einstein; Proposition~\ref{prop:divein} then already provides a divergence-free machine of order $2$ which produces all divergence-free tensors up to a finite-dimensional space, for \emph{every} Einstein metric, whereas the machine constructed below has order $3$ (and, on space forms, a non-injective symbol, see Remark~\ref{rem:symbsol}). The results of this section are therefore only of interest when $\nabla f\not\equiv0$. In that case none of the previous machines applies in general: on a compact gradient shrinking soliton with $\nabla f\not\equiv0$, the scalar curvature is not constant (otherwise the trace $-\nabla^*\nabla f=n\lambda-R$ of the soliton equation would be a constant with zero mean, hence $0$, and $f$ would be constant), so that $\Ric$ is not parallel.

\subsection{Principle}

Multiplication by $e^f$ intertwines the two divergences.

\begin{lemma}\label{lem:poids}
For every symmetric $2$-tensor $k$, $\Div_f(e^fk)=e^f\Div k$ and $\Tr(e^fk)=e^f\Tr k$. In particular, if $\mathcal P$ is a divergence-free (resp. TT) machine for $g_0$, then $e^f\mathcal P$ takes values in $\ker\Div_f$ (resp. $\ker\Div_f\cap\ker\Tr$).
\end{lemma}

\begin{proof}
$-e^f\nabla^j(e^{-f}e^fk_{ij})=-e^f\nabla^jk_{ij}$.
\end{proof}

The lemma reduces the construction of machines with vanishing weighted divergence to that of divergence-free machines for $g_0$. But the machines of the previous sections assume $g_0$ Einstein, Ricci-parallel, and so on, which excludes non-trivial solitons. We therefore need an operator $\mathfrak A$ \emph{adapted to the soliton}, with $\mathfrak A\circ L=0$. We set
\[
\mathfrak A^{*_f}\psi:=e^f\mathfrak A^*(e^{-f}\psi),
\]
which is the adjoint of $\mathfrak A$ in $L^2(d\mu_f)$.

These weighted adjoints are machines for the weighted divergence, with the same Fredholm property as before.

\begin{proposition}\label{prop:fredholmf}
Let $\mathfrak A:\Gamma(S_2)\to\Gamma(F)$ be a differential operator of order $m\ge1$ with $\mathfrak A\circ L=0$. Then $\mathfrak A^{*_f}$ takes values in $\ker\Div_f$. If moreover $M$ is closed and $\sigma_\xi(\mathfrak A)$ is injective on $\Sigma_\xi$ for every $\xi\neq0$, then, in $C^\infty$,
\[
\ker\Div_f=\Ima\mathfrak A^{*_f}\oplus\big(\ker\mathfrak A\cap\ker\Div_f\big),
\]
the sum being $L^2(d\mu_f)$-orthogonal and the second summand finite-dimensional.
\end{proposition}

\begin{proof}
For compactly supported $X$, $\int\ip{\mathfrak A^{*_f}\psi}{LX}\,d\mu_f=\int\ip{\psi}{\mathfrak ALX}\,d\mu_f=0$, so $\Div_f\mathfrak A^{*_f}\psi=0$. For the decomposition, we apply Lemma~\ref{lem:fredholmD} with the function $f$ of this section, so that $\mu=d\mu_f$, with $E=S_2$, $D=L^{*_f}=2\Div_f$ and $\mathcal P=\mathfrak A^{*_f}$, whose $L^2(d\mu_f)$-adjoint is $\mathfrak A$; $D\mathcal P=(\mathfrak A\circ L)^{*_f}=0$, and the principal symbol of $\Div_f$ is that of $\Div$, whose kernel is $\Sigma_\xi$.
\end{proof}

\subsection{The compatibility operator of Hessians}

For a metric $g$ with invertible Ricci tensor, let $\Ric^\sharp$ be the endomorphism of $T^*M$ given by $(\Ric^\sharp\omega)_i=\Ric_i{}^m\omega_m$. For $S\in\Gamma(S_2)$ we set
\[
\Theta_g(S):=-(\Ric_g^\sharp)^{-1}\big(\Div S+d\Tr S\big),\qquad
\mathcal L_g(S)_{ijk}:=\nabla_kS_{ij}-\nabla_jS_{ik}+R^m{}_{ikj}\,\Theta_g(S)_m .
\]
The tensor $\mathcal L_g(S)$ is skew-symmetric in $(j,k)$; we denote by $F$ the corresponding bundle.

As its name indicates, $\mathcal L_g$ annihilates Hessians.

\begin{lemma}\label{lem:Lhess}
\leavevmode
\begin{enumerate}[label=(\roman*)]
\item $\mathcal L_g$ is natural in $g$, on the open set of metrics with invertible Ricci tensor.
\item For every function $u$, $\mathcal L_g(\Hess_gu)=0$. Moreover $\mathcal L_g(g)=0$.
\end{enumerate}
\end{lemma}

\begin{proof}
(i) is clear. (ii) Let $H=\Hess u$ and $\omega=du$. Since $\nabla_kH_{ij}=\nabla_k\nabla_j\omega_i$, the Ricci identity gives
\[
\nabla_kH_{ij}-\nabla_jH_{ik}=[\nabla_k,\nabla_j]\omega_i=-R^m{}_{ikj}\,\omega_m .
\]
On the other hand, $\nabla^j\nabla_i\omega_j=\nabla_i\nabla^j\omega_j+\Ric_i{}^m\omega_m$. Hence $\Div H+d\Tr H=-\Ric^\sharp\omega$, so $\Theta_g(H)=\omega$, and the two terms of $\mathcal L_g(H)$ cancel. For $S=g$, we have $\nabla g=0$, $\Div g=0$ and $d\Tr g=0$.
\end{proof}

\subsection{The machine}

The Ricci tensor of a gradient soliton differs from a multiple of the metric by a Hessian, so that $\mathcal L_g(\Ric_g)$ vanishes on solitons; its linearisation is the machine we are looking for.

\begin{theorem}\label{thm:soliton}
Let $g_0$ be a metric with invertible Ricci tensor such that $\mathcal L_{g_0}(\Ric_{g_0})=0$. By Lemma~\ref{lem:Lhess}(ii), this holds for every gradient Ricci soliton $\Ric_0+\Hess f=\lambda g_0$; the case of interest is $\nabla f\not\equiv0$. Set $N(g):=\mathcal L_g(\Ric_g)$ and $\mathfrak A:=DN(g_0)$. Then:
\begin{enumerate}[label=(\alph*)]
\item $\mathfrak A\circ L=0$. Hence $\mathfrak A^*$ is a divergence-free machine of order $3$, and $\mathfrak A^{*_f}$ is a machine with vanishing weighted divergence.
\item For a gradient soliton,
\[
\mathfrak A=\mathcal L_{g_0}\circ\mathcal T,\qquad \mathcal T(k)_{ij}:=D\Ric(g_0)(k)_{ij}-\delta\Gamma(k)^m_{ij}\,\nabla_mf-\lambda k_{ij},
\]
where $\mathcal T$ is the linearisation in $g$ of $\Ric_g+\Hess_gf-\lambda g$, with $f$ fixed.
\item $g_0\in\ker\mathfrak A$. For a gradient soliton, $\Ric_0\in\ker\mathfrak A\cap\ker\Div_f$, so that $\Ric_0$ is $L^2(d\mu_f)$-orthogonal to $\Ima\mathfrak A^{*_f}$.
\end{enumerate}
\end{theorem}

\begin{proof}
(a) $N$ is natural and vanishes at $g_0$; apply Proposition~\ref{prop:principediv}, then the first part of Proposition~\ref{prop:fredholmf}.

(b) We have $\mathfrak A(k)=\mathcal L_{g_0}(D\Ric(g_0)k)+(D\mathcal L(g_0)k)(\Ric_0)$. Differentiating the identities $\mathcal L_g(g)=0$ and $\mathcal L_g(\Hess_gf)=0$, which hold for every $g$, gives
\[
(D\mathcal L(g_0)k)(g_0)=-\mathcal L_{g_0}(k),\qquad
(D\mathcal L(g_0)k)(\Hess f)=-\mathcal L_{g_0}\big(D(\Hess f)(k)\big)=\mathcal L_{g_0}\big(\delta\Gamma(k)^m_{ij}\nabla_mf\big).
\]
Conclude with $\Ric_0=\lambda g_0-\Hess f$.

(c) By scale invariance of $\Ric$, $D\Ric(g_0)(g_0)=0$, and $\delta\Gamma(g_0)=0$. Hence $\mathcal T(g_0)=-\lambda g_0$ and $\mathfrak A(g_0)=-\lambda\mathcal L_{g_0}(g_0)=0$. For a gradient soliton, $\Ric_0=\lambda g_0-\frac12L(\nabla f)$, so $\mathfrak A(\Ric_0)=0$ by (a). Finally $\Div_f\Ric_0=0$, since $\Div\Ric=-\frac12dR$ and $dR=2\Ric(\nabla f,\cdot)$ for a gradient soliton \cite{CHI}.
\end{proof}

\begin{remark}\label{rem:symbsol}
The assumption that $\Ric_0$ is invertible excludes the cylinders $S^k\times\mathbb R^{n-k}$ and the Gaussian soliton. It is local: on an open set $U$ where $\Ric_0$ is invertible, $\mathfrak A$ is defined and Proposition~\ref{prop:approx}(a)--(b) applies to $\mathfrak A^{*_f}$. By Proposition~\ref{prop:fredholmf}, a Fredholm property would follow from the injectivity of $\sigma_\xi(\mathfrak A)$ on $\Sigma_\xi$; this is a pointwise algebraic condition on $(\Rm_0,\Ric_0)$, which fails on space forms, and which we do not discuss further.
\end{remark}

\begin{remark}[Origin]
If $\Ric+\Hess f-\frac1mdf\otimes df=\lambda g$, a warped product $g\oplus e^{-2f/m}g_F$ with a suitable Einstein fibre $F^m$ is Einstein \cite{CSW}, and for $H=h\oplus w\,e^{-2f/m}g_F$ the TT condition reads $w=-\frac1m\Tr h$, $\Div_fh+\frac1m\Tr h\,df=0$. As $m\to\infty$ it becomes $\Div_fh=0$, and conformal directions become translations of $f$; condition (ii) of Proposition~\ref{prop:principe} becomes $\partial_fN=0$ for an operator $N(g,f)$ natural in $(g,f)$. The operator $\mathcal L_g(\Ric_g+\Hess_gf-\lambda g)=\mathcal L_g(\Ric_g)$ satisfies it; this is how Theorem~\ref{thm:soliton} was found.
\end{remark}

\subsection{Application to the stability of shrinking solitons}

Let $(M,g_0,f)$ be a compact shrinking soliton with $\nabla f\not\equiv0$. By \cite{CHI,CZ}, $g_0$ is $\nu$-linearly unstable as soon as there is $h\in\ker\Div_f$, $L^2(d\mu_f)$-orthogonal to $\Ric_0$, with
\[
\int\Big(-\tfrac12|\nabla h|^2+\ip{\Rm(h)}{h}\Big)\,d\mu_f>0 .
\]
For $h\in\ker\Div_f$, $\int\ip{h}{\Hess f}\,d\mu_f=\frac12\int\ip{h}{L\nabla f}\,d\mu_f=\int\ip{\Div_fh}{\nabla f}\,d\mu_f=0$, so that, since $\Ric_0=\lambda g_0-\Hess f$, the orthogonality to $\Ric_0$ is equivalent to $\int\Tr h\,d\mu_f=0$; this is why the rank-one scalar term of \cite{CZ} disappears on this subspace. Theorem~\ref{thm:soliton} provides candidates $h=\mathfrak A^{*_f}\psi$, with $\psi$ arbitrary, on the open set where $\Ric_0$ is invertible. They automatically have vanishing weighted divergence, compact support if $\psi$ has, and they are orthogonal to $\Ric_0$ (Theorem~\ref{thm:soliton}(c)), hence satisfy $\int\Tr h\,d\mu_f=0$.

\section{Machines for the linearised constraint equations}\label{sec:contraintes}

In this section a \emph{machine for the linearised constraints} at an initial data set $(g_0,K_0)$ is a linear differential operator whose image consists of solutions $(h,k)$ of the linearised constraint equations at $(g_0,K_0)$. As for TT tensors, it produces smooth compactly supported solutions from compactly supported potentials, which is the linear ingredient of gluing constructions \cite{Co,DC}; non-local solution operators with prescribed support are constructed in \cite{IMOT,MOT}. A vacuum initial data set generates by evolution a vacuum spacetime, its maximal globally hyperbolic development \cite{CB,CBG}, in which it is a spacelike hypersurface. The guiding principle replaces diffeomorphisms of $M$ by diffeomorphisms of this spacetime, that is, the Killing operator by the Killing initial data operator $D\Phi^*$: a natural tensor of the data which vanishes on all hypersurfaces of the development yields a machine (Proposition~\ref{prop:principegen}). We first treat time-symmetric data, $K_0=0$, for which the argument is elementary and also applies in Riemannian signature, as soon as the data are a totally geodesic hypersurface of an Einstein manifold, for instance of the static extension $g_0\oplus\varepsilon V^2d\theta^2$ when there is a static potential $V$ (Proposition~\ref{prop:principecontr}, Theorems~\ref{thm:contrformes} and~\ref{thm:contrlocsym}); we then treat general vacuum data (Section~\ref{ss:contrgen}).

\subsection{Setting}

Let $\varepsilon\in\{1,-1\}$ and $\bar\lambda\in\mathbb R$. For a metric $g$ and a symmetric $2$-tensor $K$ on $M^n$ we set
\[
\Phi(g,K)=\Big(R_g+\varepsilon\big(|K|_g^2-(\Tr_gK)^2\big)-(n-1)\bar\lambda,\ \Div_g\big(K-(\Tr_gK)\,g\big)\Big).
\]
If $\Sigma$ is a hypersurface of an Einstein manifold $(\bar M^{n+1},\bar g)$, $\Ric(\bar g)=\bar\lambda\bar g$, with unit normal $\nu$, $\bar g(\nu,\nu)=\varepsilon$, induced metric $g$ and second fundamental form $K(X,Y)=\bar g(\bar\nabla_X\nu,Y)$, the Gauss and Codazzi equations give $\Phi(g,K)=0$: for $\varepsilon=-1$ these are the vacuum constraint equations with cosmological constant, for $\varepsilon=1$ their Riemannian analogue. At time-symmetric data, the terms quadratic in $K$ do not contribute and
\begin{equation}\label{eq:DPhi}
D\Phi(g_0,0)(h,k)=\big(DR(g_0)h,\ \Div_{g_0}(k-(\Tr_{g_0}k)\,g_0)\big),
\end{equation}
independently of $\varepsilon$. We use
\[
DR(g)^*\alpha=\Hess \alpha+(\nabla^*\nabla \alpha)\,g-\alpha\Ric,\qquad \Rev(q):=q-\tfrac1{n-1}(\Tr q)\,g,
\]
where $\Rev$ is self-adjoint and invertible on $S_2$, with $\Rev(q)-\Tr\Rev(q)\,g=q$. With $\mathcal E:=\Ric-\frac{R}{n-1}g$, $R=R(g)$, a direct computation gives
\begin{equation}\label{eq:KDR}
\Rev\big(DR(g)^*\alpha\big)=\Hess \alpha-\alpha \mathcal E .
\end{equation}

A function $V\not\equiv0$ with $DR(g_0)^*V=0$ is a \emph{static potential}. Then $\Hess V=V\mathcal E$ by \eqref{eq:KDR}, and, taking the trace of $DR(g_0)^*V=0$, $\Delta V=-\frac{R}{n-1}V$ with $\Delta=-\nabla^*\nabla$; the scalar curvature is constant \cite{FM}. Hence, where $V\neq0$, the metric $g_0\oplus\varepsilon V^2d\theta^2$ has Ricci tensor $\Ric_0-V^{-1}\Hess V=\frac{R}{n-1}g_0$ in the directions of $M$ and $-V^{-1}\Delta V=\frac{R}{n-1}$ in the direction of $\partial_\theta$: it is Einstein with $\bar\lambda=\frac{R}{n-1}$, for both signatures, and $M=\{\theta=0\}$ is totally geodesic. We call $g_0\oplus\varepsilon V^2d\theta^2$ the \emph{static extension} of the data.

\subsection{The principle}

In the time-symmetric case, the diffeomorphisms of the ambient manifold act through normal deformations of $\Sigma$, whose effect on the data we compute first.

\begin{lemma}\label{lem:variationK}
Let $\Sigma$ be a totally geodesic hypersurface of an Einstein manifold $(\bar M,\bar g)$, $\Ric(\bar g)=\bar\lambda\bar g$, $\bar g(\nu,\nu)=\varepsilon$, with induced metric $g_0$. Then $R(g_0)=(n-1)\bar\lambda$. For $\alpha\in C^\infty_c(\Sigma)$ and small $s$, let $g_s,K_s$ be the metric and second fundamental form of $\psi_s(\Sigma)$, $\psi_s(x)=\exp_x(s\alpha(x)\nu(x))$, pulled back to $\Sigma$ by $\psi_s$. Then
\[
\frac{d}{ds}\Big|_{s=0}(g_s,K_s)=\Big(0,\ -\varepsilon\big(\Hess \alpha-\alpha(\Ric_0-\bar\lambda g_0)\big)\Big).
\]
\end{lemma}

\begin{proof}
Since $\Sigma$ is totally geodesic, the Gauss equation gives $\bar R_{kilj}=R_{kilj}$ for tangent indices, so $\bar\lambda g_{ij}=\Ric(\bar g)_{ij}=\varepsilon\bar R_{\nu i\nu j}+\Ric(g_0)_{ij}$; taking the trace and using $g^{ij}\bar R_{\nu i\nu j}=\Ric(\bar g)(\nu,\nu)=\varepsilon\bar\lambda$ gives $R(g_0)=(n-1)\bar\lambda$, and
\begin{equation}\label{eq:Rnunu}
\bar R_{\nu i\nu j}=\varepsilon\big(\bar\lambda g_0-\Ric_0\big)_{ij}.
\end{equation}
Let $W=\alpha\nu$ be the variation field and $e_i=\partial_i$. Then $\dot g_{ij}=\bar g(\bar\nabla_{e_i}W,e_j)+\bar g(e_i,\bar\nabla_{e_j}W)=2\alpha K_0(e_i,e_j)=0$. Differentiating $\bar g(\nu_s,\psi_{s*}e_i)=0$ and $\bar g(\nu_s,\nu_s)=\varepsilon$ gives $\bar\nabla_s\nu=-\varepsilon\nabla \alpha$ at $s=0$, since $\bar\nabla_s\psi_{s*}e_i=\bar\nabla_{e_i}W=(\partial_i\alpha)\nu+\alpha K_0(e_i)$. Finally $K_s(e_i,e_j)=\bar g(\bar\nabla_{e_i}\nu_s,\psi_{s*}e_j)$, and, since $\bar\nabla_{e_i}\nu=0$ at $s=0$,
\[
\dot K_{ij}=\bar g\big(\bar\nabla_{e_i}\bar\nabla_s\nu+\bar R(W,e_i)\nu,e_j\big)=-\varepsilon\nabla_i\nabla_j\alpha-\alpha\bar R_{\nu i\nu j},
\]
where we used that $\Sigma$ is totally geodesic, so that the tangential part of $\bar\nabla_{e_i}\nabla \alpha$ is $\nabla_{e_i}\nabla \alpha$, and $\bar g(\bar R(\nu,e_i)\nu,e_j)=\bar R_{j\nu\nu i}=-\bar R_{\nu i\nu j}$ with the conventions of Section~\ref{sec:conv}. We conclude with \eqref{eq:Rnunu}.
\end{proof}

A natural tensor of the data which is insensitive to these deformations yields machines for both parts of the linearised constraints.

\begin{proposition}[Principle]\label{prop:principecontr}
Let $\Sigma$, $\bar g$, $g_0$ be as in Lemma~\ref{lem:variationK}, and let $T(g,K)$ be a differential operator, with values in a tensor bundle $F$, natural in $(g,K)$ (that is, $T(\varphi^*g,\varphi^*K)=\varphi^*T(g,K)$ for local diffeomorphisms $\varphi$ of $\Sigma$), such that $T$ vanishes on the data induced by $\bar g$ on all the hypersurfaces $\psi_s(\Sigma)$ of Lemma~\ref{lem:variationK}, including $T(g_0,0)=0$. Then, with $\mathcal E_0=\Ric_0-\bar\lambda g_0$,
\[
D_gT(g_0,0)\circ L=0,\qquad D_KT(g_0,0)\big(\Hess \alpha-\alpha\mathcal E_0\big)=0\quad\text{for all }\alpha .
\]
Consequently, for all compactly supported sections $\Psi_1,\Psi_2$ of $F$,
\[
(h,k)=\Big(\Rev\,D_KT(g_0,0)^*\Psi_1,\ \Rev\,D_gT(g_0,0)^*\Psi_2\Big)
\]
is a solution of the linearised constraint equations $D\Phi(g_0,0)(h,k)=0$.
\end{proposition}

\begin{proof}
The first identity follows from $T(\varphi_t^*g_0,0)=\varphi_t^*T(g_0,0)=0$, as in Proposition~\ref{prop:principe}. The second follows from Lemma~\ref{lem:variationK}, by differentiating $T(g_s,K_s)=0$ at $s=0$. By the first identity, $q=D_gT(g_0,0)^*\Psi_2$ is divergence free (Proposition~\ref{prop:principediv}), and $k=\Rev q$ satisfies $k-(\Tr k)g_0=q$. By \eqref{eq:KDR}, with $R(g_0)=(n-1)\bar\lambda$, the second identity reads $D_KT(g_0,0)\circ\Rev\circ DR(g_0)^*=0$; taking adjoints, $DR(g_0)\circ\Rev\circ D_KT(g_0,0)^*=0$. We conclude with \eqref{eq:DPhi}.
\end{proof}

The second identity expresses the fact that a deformation of $\Sigma$ in the normal direction does not change $T$, just as the first expresses invariance under diffeomorphisms of $\Sigma$; together, they are the infinitesimal version of the invariance of $T$ under diffeomorphisms of $\bar M$. The momentum part only involves a divergence-free machine for $g_0$ (Section~\ref{sec:divnulle}). The new ingredient is the Hamiltonian part: a machine for $\ker DR(g_0)$, given by the linearisation \emph{in $K$} of a tensor computed from the data. Natural tensors $T$ as in Proposition~\ref{prop:principecontr} are obtained from natural tensors $\bar T(\bar g)$ vanishing at $\bar g$, whose components along $\Sigma$ can be expressed in terms of $(g,K)$ by the Gauss and Codazzi equations and the Einstein equation.

\subsection{Space forms}

On a space form, the Codazzi tensor gives a machine of order $1$ for the Hamiltonian part.

\begin{theorem}\label{thm:contrformes}
Let $g_0$ have constant sectional curvature $\kappa$, $n\ge3$, so that $(g_0,0)$ satisfies $\Phi(g_0,0)=0$ with $\bar\lambda=n\kappa$. For $u\in S_2$ let $(d^\nabla u)_{ijk}=\nabla_ku_{ij}-\nabla_ju_{ik}$. Then:
\begin{enumerate}[label=(\roman*)]
\item $d^\nabla\circ\Rev\circ DR(g_0)^*=0$; hence $h=\Rev(d^\nabla)^*\Psi$ satisfies $DR(g_0)h=0$ for every $\Psi$, and, with $N_\kappa=\Rm-\frac\kappa2g\KN g$ as in Proposition~\ref{prop:divcc},
\[
(\Psi,Z)\longmapsto\big(\Rev(d^\nabla)^*\Psi,\ \Rev\,DN_\kappa(g_0)^*Z\big)
\]
is a machine for the linearised constraints at $(g_0,0)$, for both signatures; here $\Psi$ is a section of the bundle of $3$-tensors which are skew-symmetric in their last two indices, where $d^\nabla$ takes its values, and $Z$ a section of the bundle $\mathcal C$ of algebraic curvature tensors.
\item If $M$ is closed, $\ker DR(g_0)=\Ima\Rev(d^\nabla)^*\oplus\big(\ker DR(g_0)\cap\ker(d^\nabla\Rev)\big)$, the second summand being finite dimensional; together with Proposition~\ref{prop:divcc}, the machine produces all solutions of the linearised constraints at $(g_0,0)$ up to a finite-dimensional space.
\end{enumerate}
\end{theorem}

\begin{proof}
Here $\Ric_0=(n-1)\kappa g_0$, $R=n(n-1)\kappa$ and $\mathcal E=-\kappa g_0$. By \eqref{eq:KDR} and the Ricci identity of Lemma~\ref{lem:Lhess}(ii),
\[
d^\nabla\big(\Hess \alpha+\kappa \alpha g_0\big)_{ijk}=-R^m{}_{ikj}\nabla_m\alpha+\kappa\big(\nabla_k\alpha\,g_{ij}-\nabla_j\alpha\,g_{ik}\big)=0,
\]
since $R^m{}_{ikj}=\kappa(\delta^m_kg_{ij}-\delta^m_jg_{ik})$. This gives (i), with \eqref{eq:DPhi} and Proposition~\ref{prop:divcc} for the momentum part. For (ii), $\sigma_\xi(d^\nabla)u=i(\xi_ku_{ij}-\xi_ju_{ik})$ vanishes only for $u\in\mathbb R\,\xi\otimes\xi$ (contract with $\xi^k$), so $\sigma_\xi(d^\nabla\Rev)h=0$ forces $h=c(\xi\otimes\xi-|\xi|^2g_0)$, for which $|\xi|^2\Tr h-h(\xi^\sharp,\xi^\sharp)=-c(n-1)|\xi|^4$. Thus $\ker\sigma_\xi(DR)\cap\ker\sigma_\xi(d^\nabla\Rev)=\{0\}$, and Lemma~\ref{lem:fredholmD} applies with $f=0$, $E=S_2$, $D=DR(g_0):\Gamma(S_2)\to C^\infty(M)$ and $\mathcal P=\Rev(d^\nabla)^*$, since $D\mathcal P=0$ by (i) and $\mathcal P^*=d^\nabla\Rev$.
\end{proof}

\begin{remark}
The operator $d^\nabla K$ is the Codazzi expression of $\bar R(\nu,\cdot,\cdot,\cdot)$, and the theorem is Proposition~\ref{prop:principecontr} for $\bar T=\Rm(\bar g)-\frac{\kappa}{2}\bar g\KN\bar g$ on the space form of curvature $\kappa$ containing $(M,g_0)$ as a totally geodesic hypersurface; only the Codazzi component contributes to $D_KT(g_0,0)$, the others being quadratic in $K$. In dimension $3$, Beig and Chruściel \cite{BC20} have given potentials for the vector constraint and for the trace-free part of the linearised scalar constraint on Einstein manifolds, that is, on space forms, the latter of curl type, as $(d^\nabla)^*$ in dimension $3$; Theorem~\ref{thm:contrformes} also controls the trace, and holds in every dimension.
\end{remark}

\subsection{Locally symmetric static data}

For locally symmetric data with a static potential, a second-order operator built from the curvature plays the role of $d^\nabla$.

\begin{theorem}\label{thm:contrlocsym}
Let $(M,g_0)$ be connected and locally symmetric, $\nabla\Rm_0=0$, $n\ge3$, with a static potential $V\not\equiv0$, $DR(g_0)^*V=0$, and let $\bar\lambda=\frac{R}{n-1}$ and $\mathcal E=\Ric_0-\bar\lambda g_0$. Define the second-order operator $\mathcal Q:\Gamma(S_2)\to\Gamma(T^0_4)$ by
\[
\mathcal Q(u)_{lijk}:=\nabla_l\big(\nabla_ku_{ij}-\nabla_ju_{ik}\big)-u_l{}^mR_{mijk}-u_{lj}\mathcal E_{ik}+u_{lk}\mathcal E_{ij}.
\]
\begin{enumerate}[label=(\roman*)]
\item $\mathcal Q\circ\Rev\circ DR(g_0)^*=0$. Hence $h=\Rev \mathcal Q^*\Psi$ satisfies $DR(g_0)h=0$, and
\[
(\Psi,Z)\longmapsto\big(\Rev \mathcal Q^*\Psi,\ \Rev\,D(\nabla\Rm)(g_0)^*Z\big)
\]
is a machine for the linearised constraints at $(g_0,0)$, for both signatures; here $\Psi$ is a section of $T^0_4M$ and $Z$ a section of $T^0_5M$.
\item If $M$ is closed, $\ker DR(g_0)=\Ima\Rev \mathcal Q^*\oplus\big(\ker DR(g_0)\cap\ker(\mathcal Q\Rev)\big)$, the second summand being finite dimensional; together with Proposition~\ref{prop:divlocsym}, the machine produces all solutions of the linearised constraints at $(g_0,0)$ up to a finite-dimensional space.
\end{enumerate}
\end{theorem}

\begin{proof}
(i) Since $\Rm_0$ is parallel, so are $\Ric_0$ and $\mathcal E$, and $R$ is constant. By \eqref{eq:KDR} we must show $\mathcal Q(U)=0$ for $U=H-\alpha\mathcal E$, $H=\Hess \alpha$. By the Ricci identity, $(d^\nabla H)_{ijk}=-R^m{}_{ikj}\nabla_m\alpha$, and $d^\nabla(\alpha\mathcal E)_{ijk}=\nabla_k\alpha\,\mathcal E_{ij}-\nabla_j\alpha\,\mathcal E_{ik}$. Since $\nabla\Rm_0=0$ and $\nabla \mathcal E=0$,
\[
\nabla_l\big(d^\nabla U\big)_{ijk}=-R_{mikj}H_l{}^m-H_{lk}\mathcal E_{ij}+H_{lj}\mathcal E_{ik}.
\]
The second-order terms of $\mathcal Q(U)$ are therefore
\[
-R_{mikj}H_l{}^m-H_{lk}\mathcal E_{ij}+H_{lj}\mathcal E_{ik}-H_l{}^mR_{mijk}-H_{lj}\mathcal E_{ik}+H_{lk}\mathcal E_{ij}=-H_l{}^m\big(R_{mikj}+R_{mijk}\big)=0,
\]
there are no first-order terms, and the zero-order terms are $\alpha\,Z_{lijk}$ with
\[
Z_{lijk}:=\mathcal E_l{}^mR_{mijk}+\mathcal E_{lj}\mathcal E_{ik}-\mathcal E_{lk}\mathcal E_{ij}.
\]
We show $Z=0$ using $V$. The computation above shows that $\mathcal Q(\Hess\alpha-\alpha\mathcal E)=\alpha Z$ for every function $\alpha$. By \eqref{eq:KDR}, $\Hess V-V\mathcal E=0$, hence $VZ=0$. Since $\Hess V=V\mathcal E$ is a second-order linear ordinary differential equation along geodesics, $V$ cannot vanish on a non-empty open set of the connected manifold $M$ without vanishing identically; hence $Z=0$ on a dense set, and everywhere. The momentum part is Proposition~\ref{prop:divlocsym}.

(ii) The principal symbol of $\mathcal Q$ is $\sigma_\xi(\mathcal Q)u=-\xi_l(\xi_ku_{ij}-\xi_ju_{ik})$, which vanishes only for $u\in\mathbb R\,\xi\otimes\xi$. As in the proof of Theorem~\ref{thm:contrformes}, $\ker\sigma_\xi(DR)\cap\ker\sigma_\xi(\mathcal Q\Rev)=\{0\}$, and Lemma~\ref{lem:fredholmD} applies with $f=0$, $D=DR(g_0)$ and $\mathcal P=\Rev \mathcal Q^*$, whose adjoint is $\mathcal Q\Rev$.
\end{proof}

\begin{remark}[Origin of $\mathcal Q$]
In the setting of Proposition~\ref{prop:principecontr}, let $\bar g$ be locally symmetric and $\bar T=\bar\nabla\bar\Rm$. The tangential components $T_{lijk}=(\bar\nabla_{e_l}\bar\Rm)(\nu,e_i,e_j,e_k)$ are expressed in terms of $(g,K)$ by $\bar R(\nu,e_i,e_j,e_k)=(d^\nabla K)_{ijk}$ (Codazzi), $\bar\nabla_{e_l}\nu=K(e_l)$, $\bar\nabla_{e_l}e_i=\nabla_{e_l}e_i-\varepsilon K_{li}\nu$, the Gauss equation, and the formula for $\bar R_{\nu i\nu j}$ derived from the Gauss and Einstein equations, which reduces to \eqref{eq:Rnunu} at $K=0$. Linearising in $K$ at $K=0$ gives exactly $D_KT(g_0,0)=\mathcal Q$. Theorem~\ref{thm:contrlocsym} has been proved directly, without assuming that $g_0\oplus\varepsilon V^2d\theta^2$ is locally symmetric.
\end{remark}

\begin{example}\label{ex:contr}
Let $(M_1,g_1)$ be a locally symmetric Einstein manifold with $\Ric(g_1)=m\kappa\,g_1$, $m\ge1$, $\kappa\neq0$, and $(M_2,g_2)$ an $m$-dimensional space form of curvature $\kappa$. On $M=M_1\times M_2$, $g_0=g_1\oplus g_2$ is locally symmetric, $R=m\kappa(n-1)$, $\bar\lambda=m\kappa$ and $\mathcal E=0\oplus(-\kappa g_2)$. Every function $V$ on $M_2$ with $\Hess_{g_2}V=-\kappa Vg_2$ (a restriction of a linear function if $M_2$ is the round sphere $S^m(\kappa^{-1/2})$, $V=\cos(\sqrt\kappa\,x)$ if $m=1$) satisfies $\Hess V=V\mathcal E$, hence is a static potential, and the static extension is the Einstein product of $M_1$ with the space form of dimension $m+1$ and curvature $\kappa$. For instance:
\begin{itemize}
\item $M=S^2\times S^1$ with $g_0=g_{S^2}\oplus dx^2$ and $V=\cos x$: the Riemannian static extension is $S^2\times S^2$, the Lorentzian one the Nariai solution $dS_2\times S^2$;
\item $M=S^p\times S^{p-1}$ (unit spheres, $p\ge2$), with static extension $S^p\times S^p$, and $M=H^p\times H^{p-1}$, with static extension $H^p\times H^p$ (and $H^2\times AdS_2$ in the Lorentzian case for $p=2$); the factor $H^p$ may be replaced by any of its quotients;
\item $M=\mathbb{CP}^q\times S^m$ with appropriate normalisations.
\end{itemize}
In these examples the static extension is only defined where $V\neq0$; for $S^2\times S^1$ with $V=\cos x$, it consists of two bands $\{\cos x\neq0\}$, each of which completes to $S^2\times S^2$ or to a static patch of the Nariai solution. These data are not Einstein, and Theorem~\ref{thm:contrlocsym} applies; for $S^2\times S^1$, $S^p\times S^{p-1}$ and $\mathbb{CP}^q\times S^m$, $M$ is closed and (ii) applies.
\end{example}

\subsection{General vacuum data}\label{ss:contrgen}

We now drop the assumption $K_0=0$, in Lorentzian signature $\varepsilon=-1$. We use the Hamiltonian variables $x=(g,\pi)$, where $\pi=\tilde\pi\,dV_g$ is a symmetric contravariant tensor density and $\tilde\pi=\big(K-(\Tr_gK)g\big)^{\sharp\sharp}$, so that $K=\tilde\pi^\flat-\frac{1}{n-1}(\Tr_g\tilde\pi)\,g$; $(g,K)\mapsto(g,\pi)$ is a diffeomorphism. Tangent vectors are pairs $y=(h,p)$, with $h$ a covariant symmetric tensor and $p$ a contravariant symmetric density, and the symplectic form is
\[
\omega(y_1,y_2)=\int_M\big(\ip{h_1}{p_2}-\ip{h_2}{p_1}\big),
\]
which only involves the natural contraction of covariant tensors with contravariant densities. For a functional $\mathcal F(g,\pi)$ we define $\frac{\delta\mathcal F}{\delta g}$ (a contravariant density) and $\frac{\delta\mathcal F}{\delta\pi}$ (a covariant tensor) by $D\mathcal F(y)=\int_M\big(\ip{\frac{\delta\mathcal F}{\delta g}}{h}+\ip{\frac{\delta\mathcal F}{\delta\pi}}{p}\big)$, and we set $J(A,B)=(B,-A)$, so that $\omega\big(J(\frac{\delta\mathcal F}{\delta g},\frac{\delta\mathcal F}{\delta\pi}),y\big)=D\mathcal F(y)$. The constraint map is normalised as
\[
\Phi^{\mathrm H}_0(g,\pi)=\Big(|\tilde\pi|_g^2-\tfrac{1}{n-1}(\Tr_g\tilde\pi)^2-R_g+(n-1)\bar\lambda\Big)dV_g,\qquad \Phi^{\mathrm H}_1(g,\pi)=2\,(\Div_g\tilde\pi)\,dV_g ;
\]
since $|\tilde\pi|^2-\frac{1}{n-1}(\Tr\tilde\pi)^2=|K|^2-(\Tr K)^2$, we have $\Phi^{\mathrm H}=(-\Phi_0\,dV_g,2\Phi_1\,dV_g)$ for $\varepsilon=-1$, so that $\ker D\Phi^{\mathrm H}(x_0)$ and $\ker D\Phi(g_0,K_0)$ correspond to each other under the differential of $(g,K)\mapsto(g,\pi)$. We set $\mathcal H_{\alpha,X}(g,\pi)=\int_M\big(\alpha\,\Phi^{\mathrm H}_0(g,\pi)+\ip{X}{\Phi^{\mathrm H}_1(g,\pi)}\big)$.

Let $(\bar M,\bar g)$ be a vacuum spacetime, $\Ric(\bar g)=\bar\lambda\bar g$, and $\Sigma\subset\bar M$ a spacelike hypersurface with data $x_0$; for instance, $(\bar M,\bar g)$ is the development of a vacuum initial data set $x_0$ \cite{CB,CBG}. For compactly supported $(\alpha,X)$, let $Y$ be a vector field on $\bar M$ equal to $\alpha\nu+X$ along $\Sigma$, $\nu$ the future unit normal, $\varphi_s$ the flow of $Y$, and $x_s$ the data of $\varphi_s(\Sigma)$ pulled back to $\Sigma$ by $\varphi_s$. The first variation $\Gamma(\alpha,X):=\frac{d}{ds}\big|_{s=0}x_s$ only depends on $(\alpha,X)$ (Lemma~\ref{lem:variationgen}), and the vacuum evolution equations in Hamiltonian form state that $\Gamma(\alpha,X)=J\big(\frac{\delta\mathcal H_{\alpha,X}}{\delta g},\frac{\delta\mathcal H_{\alpha,X}}{\delta\pi}\big)$ at $x_0$, that is,
\begin{equation}\label{eq:hamilton}
\int_M\big(\alpha\,D\Phi^{\mathrm H}_0(x_0)y+\ip{X}{D\Phi^{\mathrm H}_1(x_0)y}\big)=\omega\big(\Gamma(\alpha,X),y\big)\qquad\text{for all compactly supported }y .
\end{equation}
For $n=3$ and $\bar\lambda=0$, this is the Hamiltonian formulation of Arnowitt, Deser and Misner \cite{ADM,FMa}; by duality, \eqref{eq:hamilton} also expresses the fact that $D\Phi^{\mathrm H}(x_0)^*$ is, up to $J$, the operator $(\alpha,X)\mapsto\Gamma(\alpha,X)$, whose kernel consists of the Killing initial data \cite{Mo,FMM,BCh}. The derivation is the same in every dimension and with a cosmological constant; for completeness and to fix conventions, we include a proof in Appendix~\ref{app:hamilton} (Proposition~\ref{prop:hamilton}).

With this identity at hand, the principle extends to arbitrary vacuum data.

\begin{proposition}[Principle for general vacuum data]\label{prop:principegen}
Let $T(g,\pi)$ be a differential operator with values in a tensor bundle $F$, natural in $(g,\pi)$, which vanishes on the data of all the hypersurfaces $\varphi(\Sigma)$, $\varphi$ a compactly supported diffeomorphism of $\bar M$ close to the identity. For $\Psi\in\Gamma_c(F)$ let $\mathcal P\Psi$ be defined by
\[
\omega(\mathcal P\Psi,y)=\int_M\ip{\Psi}{DT(x_0)y}\,dV_{g_0}\quad\text{for all }y,\qquad\text{i.e.}\quad \mathcal P\Psi=\big(D_\pi T(x_0)^*\Psi,\,-D_gT(x_0)^*\Psi\big).
\]
Then $D\Phi^{\mathrm H}(x_0)\circ\mathcal P=0$: $\mathcal P$ is a machine for the linearised constraints at $x_0$.
\end{proposition}

\begin{proof}
Since $T(x_s)=0$ for all $s$, $DT(x_0)\Gamma(\alpha,X)=0$. By \eqref{eq:hamilton} and the antisymmetry of $\omega$,
\[
\int_M\big(\alpha\,D\Phi^{\mathrm H}_0(\mathcal P\Psi)+\ip{X}{D\Phi^{\mathrm H}_1(\mathcal P\Psi)}\big)=\omega(\Gamma(\alpha,X),\mathcal P\Psi)=-\int_M\ip{\Psi}{DT(x_0)\Gamma(\alpha,X)}\,dV_{g_0}=0
\]
for all compactly supported $(\alpha,X)$.
\end{proof}

At $\pi_0=0$, $D_\pi T=D_KT\circ\Rev$ on tensors and $p=k-(\Tr k)g_0$: Proposition~\ref{prop:principegen} reduces to Proposition~\ref{prop:principecontr}, which has an elementary proof valid in both signatures. When $\pi_0\neq0$, $\mathcal P$ mixes the two components.

On spacelike hypersurfaces of constant curvature spacetimes, the Gauss and Codazzi tensors provide such a $T$.

\begin{theorem}\label{thm:contrgen}
Let $(\bar M^{n+1},\bar g)$, $n\ge3$, be a Lorentzian manifold of constant sectional curvature $\bar\kappa$ (locally Minkowski, de Sitter or anti-de Sitter), and $\Sigma$ an arbitrary spacelike hypersurface, with data $x_0=(g_0,\pi_0)$, corresponding to $(g_0,K_0)$. With $K=K(g,\pi)$, let
\[
T_1(g,\pi)=\Rm_g+\tfrac12K\KN K-\tfrac{\bar\kappa}2\,g\KN g,\qquad T_2(g,\pi)=d^\nabla K,\qquad T=(T_1,T_2).
\]
Thus $T$ takes its values in $F=\mathcal C\oplus\big(T^*M\otimes\Lambda^2T^*M\big)$, where $\mathcal C$ is the bundle of algebraic curvature tensors, and the machine $\mathcal P$ is defined on sections $(\Psi_1,\Psi_2)$ of $F$.
\begin{enumerate}[label=(\roman*)]
\item The operator $\mathcal P$ of Proposition~\ref{prop:principegen} associated with $T$ is a machine for the linearised constraints at $x_0$.
\item If $M$ is closed, $\Ima\mathcal P$ has finite codimension in $\ker D\Phi^{\mathrm H}(x_0)$ (in $C^\infty$): $\ker D\Phi^{\mathrm H}(x_0)=\Ima\mathcal P\oplus V$ with $\dim V<\infty$.
\end{enumerate}
\end{theorem}

\begin{proof}
(i) For a spacelike hypersurface of $(\bar M,\bar g)$ with data $(g,K)$, the Gauss equation $\bar R_{ijkl}=R_{ijkl}+K_{ik}K_{jl}-K_{il}K_{jk}$ and $\bar R=\frac{\bar\kappa}2\bar g\KN\bar g$ give $T_1=0$, and the Codazzi equation $\bar R(\nu,e_i,e_j,e_k)=(d^\nabla K)_{ijk}$ gives $T_2=0$. Apply Proposition~\ref{prop:principegen}.

(ii) We identify tangent vectors $(h,p)$ with pairs of symmetric $2$-tensors by means of $g_0$ and $dV_{g_0}$, take adjoints in $L^2(dV_{g_0})$, and apply Lemma~\ref{lem:fredholmD} with $f=0$, $E=S_2\oplus S_2$, $D=D\Phi^{\mathrm H}(x_0)$, whose components $D_0,D_1$ have orders $2$ and $1$, and the machine $\mathcal P$, whose restrictions to $\Psi_1$ and $\Psi_2$ have orders $2$ and $1$; $D\mathcal P=0$ by (i), and we write $\mathcal P^*=(\mathcal P_1^*,\mathcal P_2^*)$, where $\mathcal P_a$ is the restriction of $\mathcal P$ to $\Psi_a$. When $K_0\neq0$, the components have mixed orders: $D_0$ has order $2$ in $h$ and $0$ in $p$, while $D_1$ has order $1$ in $p$ and also in $h$, through $\pi_0$ and the Christoffel symbols, and similarly for $\mathcal P_1$ and $\mathcal P_2$. Lemma~\ref{lem:fredholmD} is applied with the maximal orders $d_0=2$, $d_1=1$, $m_1=2$, $m_2=1$, the principal symbols being the homogeneous parts of these orders. The principal symbols of all components at $(z,\xi)$, $z\in M$, only depend on $g_0(z)$, $K_0(z)$ and $\xi$, continuously.

\emph{The case $K_0(z)=0$.} Then $D_\pi T_1$ and $D_gT_2$ vanish at $z$ (the first is quadratic in $K$, the second only involves the Christoffel symbols times $K$), $D_gT_1=D\Rm$ up to a term of order $0$, and $D_\pi T_2=d^\nabla\circ\Rev$. Hence, at the level of principal symbols, $\mathcal P(\Psi_1,\Psi_2)=\big(\Rev(d^\nabla)^*\Psi_2,-D\Rm^*\Psi_1\big)$: the metric component comes from the Codazzi part $T_2$, through $D_\pi T$, and the momentum component from the Gauss part $T_1$, through $D_gT$; and
\[
\sigma_\xi(\mathcal P_1^*)(h,p)=-\sigma_\xi(D\Rm)p=-\tfrac12\,\xi\xi\KN p,\qquad \sigma_\xi(\mathcal P_2^*)(h,p)=\sigma_\xi(d^\nabla\Rev)h ,
\]
while $\sigma_\xi(D_0)(h,p)$ is a non-zero multiple of $\sigma_\xi(DR)h$ and $\sigma_\xi(D_1)(h,p)$ a non-zero multiple of $p(\xi^\sharp,\cdot)$. If all of them vanish, then $h=0$ by the proof of Theorem~\ref{thm:contrformes}; and $\xi\xi\KN p=0$ forces $p=\xi\otimes\eta+\eta\otimes\xi$ (take $\xi=e_1$ and look at the components $(1,a,1,b)$, $a,b\ge2$), so that $0=p(\xi^\sharp,\cdot)=|\xi|^2\eta+\ip{\eta}{\xi}\xi$ gives $\eta=0$.

\emph{Small $K_0$.} By continuity and compactness of the unit sphere, there is $\delta_n>0$ such that the assumption of Lemma~\ref{lem:fredholmD} holds at every $(z,\xi)$ with $|K_0(z)|_{g_0}\le\delta_n$.

\emph{Scaling.} For $c>0$, $\Sigma$ is also a spacelike hypersurface of $(\bar M,c^2\bar g)$, of constant curvature $\bar\kappa/c^2$, with data $(c^2g_0,cK_0)$, and $|cK_0|_{c^2g_0}=c^{-1}|K_0|_{g_0}$. The map $(g,\pi)\mapsto(c^2g,c^{n-3}\pi)$ sends the data of every hypersurface for $\bar g$ to its data for $c^2\bar g$; let $\rho_c(h,p)=(c^2h,c^{n-3}p)$ be its differential. By naturality, $\Phi^{\mathrm H}$ and $T$ for $c^2\bar g$, composed with this map, are constant multiples (component by component) of $\Phi^{\mathrm H}$ and $T$ for $\bar g$, and $\omega(\rho_cy_1,\rho_cy_2)=c^{n-1}\omega(y_1,y_2)$. Hence $\rho_c$ maps $\ker D\Phi^{\mathrm H}(x_0)$ onto the corresponding kernel for the rescaled data. Moreover, the pairing $\int\ip{\Psi}{DT\,y}\,dV$ which defines $\mathcal P$ through $\omega$ changes under $g\mapsto c^2g$ by constants which depend on the component $T_1$ or $T_2$; hence $\rho_c\,\mathcal P(\Psi_1,\Psi_2)=\mathcal P_c(a_1\Psi_1,a_2\Psi_2)$ for some constants $a_1,a_2>0$, where $\mathcal P_c$ is the machine of the rescaled data, and $\rho_c(\Ima\mathcal P)=\Ima\mathcal P_c$. Since $\rho_c$ rescales $h$ and $p$ by different factors, the $L^2$-orthogonality of the decomposition is not preserved, which is why (ii) is stated without it. For $c\ge\delta_n^{-1}\sup_M|K_0|_{g_0}$, Lemma~\ref{lem:fredholmD} applies to the rescaled data and gives a decomposition $\ker=\Ima\oplus V_c$ with $\dim V_c<\infty$; applying $\rho_c^{-1}$ gives (ii) with $V=\rho_c^{-1}V_c$.
\end{proof}
\begin{example}\label{ex:fermees}
Part (ii) applies to closed spacelike hypersurfaces of constant curvature spacetimes. They include the round slices of de Sitter spacetime, which are umbilic (Example~\ref{ex:ombilic}), with $K_0\neq0$ except at the waist, and all spacelike graphs over them; spacelike graphs over the slices of flat spacetimes $\mathbb R\times T^n$; and spacelike graphs over the slices of the locally anti-de Sitter spacetimes $-dt^2+\cos^2t\,g_H$, $|t|<\pi/2$, where $g_H$ is a compact hyperbolic metric. These data are not time-symmetric in general.
\end{example}

The simplest data with $K_0\neq0$ are umbilic slices, for which the machine of Theorem~\ref{thm:contrgen} can be written explicitly.

\begin{example}[Umbilic slices]\label{ex:ombilic}
Let $\Sigma$ be an umbilic spacelike hypersurface, $K_0=c\,g_0$ with $c$ constant, of a spacetime of constant curvature $\bar\kappa$; by the Gauss equation, $(M,g_0)$ is a space form of curvature $\kappa=\bar\kappa-c^2$. Examples are the hyperboloids of Minkowski spacetime ($\bar\kappa=0$, $c=1$, $\kappa=-1$) and the non-equatorial round slices of de Sitter spacetime (Example~\ref{ex:fermees}). For $\Psi=(\Psi_1,\Psi_2)$, set
\[
Q_{ij}=-\nabla^k\big((\Psi_2)_{ijk}+(\Psi_2)_{jik}\big),\qquad C=c(\Psi_1),\qquad D_{ad}=2\,\mathrm{Sym}_{ad}\,\nabla^b\nabla^c(\Psi_1)_{abcd},
\]
so that $Q=\mathrm{Sym}\,(d^\nabla)^*\Psi_2$ and $D=\Theta(\cdot)^*\Psi_1$ (Lemma~\ref{lem:FE1}). In terms of the variations $(h,k)$ of $(g,K)$, the machine $\mathcal P$ reads
\begin{equation}\label{eq:ombilic}
h=\Rev\big(Q+4c\,C\big),\qquad k=\Rev\big(-D+2(\kappa+2c^2)\,C+c\,Q\big).
\end{equation}
Indeed, with $\tilde p=p/dV$, the linearisation of $K(g,\pi)$ at $(g_0,cg_0)$ is $\delta K=\Rev(\tilde p)+(2-n)c\,h+\frac12c\,(\Tr h)\,g_0$. Since $\nabla K_0=0$, $\delta(\nabla K)=\nabla(\delta K-ch)$, and Lemma~\ref{lem:FE1}(a) gives
\[
\begin{aligned}
\delta T_1&=\Theta(h)+g_0\KN\Big(c\,\Rev(\tilde p)-\big(\tfrac\kappa2+(n-1)c^2\big)h+\tfrac12c^2(\Tr h)\,g_0\Big),\\
\delta T_2&=d^\nabla\Big(\Rev(\tilde p)+(1-n)c\,h+\tfrac12c\,(\Tr h)\,g_0\Big).
\end{aligned}
\]
Using $\ip{Z}{g_0\KN u}=4\ip{c(Z)}{u}$ and $\omega(\mathcal P\Psi,y)=\int_M\ip{\Psi}{DT(x_0)y}$, one obtains $\mathcal P\Psi=(h,p)$ with $h$ as in \eqref{eq:ombilic} and $\tilde p=-D+2(\kappa+2(n-1)c^2)C-2c^2(\Tr C)g_0+(n-1)cQ-\frac12c(\Tr Q)g_0$, and then $k=\delta K$ is given by \eqref{eq:ombilic}.

One can check \eqref{eq:ombilic} directly. At $(g_0,cg_0)$, the linearised constraints read $DR(g_0)h+2(n-1)c(\Tr k-c\Tr h)=0$ and $\Div(k-(\Tr k)g_0)=(n-1)c\big(\Div h+\frac12d\Tr h\big)$. The $\Psi_2$-part is $(h,k)=(\Rev Q,c\Rev Q)$, a variation through umbilic data with the same $c$, which solves the constraints since $DR(g_0)\Rev Q=0$ (the adjoint of $d^\nabla(\Hess\alpha+\kappa\alpha g_0)=0$, see Theorem~\ref{thm:contrformes}). For the $\Psi_1$-part, $DR(g_0)\Rev C=\nabla^i\nabla^jC_{ij}+\kappa\Tr C$ and $\Tr D=-2\nabla^b\nabla^cC_{bc}$, so the Hamiltonian constraint holds, while the momentum constraint reduces to $\Div(D-2\kappa C)=0$, which is Proposition~\ref{prop:divcc}.

For $c=0$ and $n=3$, \eqref{eq:ombilic} recovers the potentials of Beig and Chruściel \cite{BC20}. Writing $\Psi_1=S\KN g_0$, the commutation $[\nabla^b,\nabla_a]S_{bd}=3\kappa S_{ad}-\kappa(\Tr S)g_{0ad}$ gives
\[
k=2\Delta S-4\,\mathrm{Sym}\,\nabla(\nabla^bS_{b\cdot})+(\nabla^a\nabla^bS_{ab})g_0-(\Delta\Tr S)g_0-4\kappa S,\qquad \Delta=-\nabla^*\nabla,
\]
and with $\psi=-\frac12c(\Psi_1)$, that is $S=-2\psi+\frac12(\Tr\psi)g_0$, this is exactly the parametrisation $k=\check Q\psi+\frac13(\Tr k)g_0$ of the momentum constraint in \cite[(1.2)--(1.4)]{BC20}, with $R=6\kappa$. Writing $(\Psi_2)_{ijk}=\epsilon_{jkl}A_{il}$, we get $Q=-2\,\mathrm{Sym}\,\mathrm{curl}\,A$, with $(\mathrm{curl}\,A)_{ij}=\epsilon_{jkl}\nabla_kA_{il}$. For $A=\phi$ symmetric and trace free, $h=\Rev Q$ is, up to sign, the potential $h_{ij}=\epsilon_{i\ell k}\nabla^\ell\phi^k{}_j+\epsilon_{j\ell k}\nabla^\ell\phi^k{}_i$ of the Hamiltonian constraint in \cite[(1.8)]{BC20}; for $A$ antisymmetric, $A_{il}=\epsilon_{ilm}w_m$, one finds $h=-L_wg_0$, a gauge direction, and $A=fg_0$ gives $Q=0$. For $c\neq0$, \eqref{eq:ombilic} couples the two potentials, through the terms $4c\,\Rev C$ in $h$ and $c\,\Rev Q$ in $k$.
\end{example}

\begin{remark}
In the same way, Proposition~\ref{prop:principegen} applied to the components of $\bar\nabla\bar\Rm$ gives machines on all spacelike hypersurfaces of locally symmetric vacuum spacetimes, such as the Nariai solution; we have established the Fredholm property only for time-symmetric data (Theorem~\ref{thm:contrlocsym}). For hypersurfaces of the Schwarzschild or Kerr spacetimes, no suitable natural tensor vanishing on the spacetime is available. On Kerr, explicit operators $\mathfrak A$ with $\mathfrak A\circ L=0$ are known: they form the compatibility complex of the Killing operator constructed in \cite{AABKW}, and give divergence-free machines on the spacetime (Remark~\ref{rem:lorentz}). They are not natural tensors of the data, so Proposition~\ref{prop:principegen} does not apply to them directly; whether they induce machines for the linearised constraints on the slices of Kerr is an interesting question, which we leave open. In dimension $3+1$, the Bach tensor vanishes on Kerr, but it vanishes on every vacuum spacetime: once its normal derivatives are eliminated through the vacuum equations, it gives a tensor of the data which vanishes on all vacuum data, so that $DT(x_0)$ vanishes on $\ker D\Phi^{\mathrm H}(x_0)$; then $\omega(\mathcal P\Psi,y)=0$ for every solution $y$ of the linearised constraints, which by \eqref{eq:hamilton} means, formally, that $\mathcal P$ only produces gauge directions $\Gamma(\alpha,X)$. A useful $T$ must vanish on the given spacetime, but not on all vacuum spacetimes. The Mars--Simon tensor vanishes on Kerr and characterises it \cite{MS}, but it involves the Killing field $\xi$; it fits the relaxed framework of Proposition~\ref{prop:relache} with $\tau=\xi$ only if the term $D_\tau F(L_Y\tau)$ can be controlled, which we have not investigated.
\end{remark}

\appendix
\section{The Hamiltonian identity}\label{app:hamilton}

Identity \eqref{eq:hamilton} is the well-known Hamiltonian form of the vacuum Einstein equations \cite{ADM,FMa,Mo,FMM}. We include a proof in the setting of Section~\ref{ss:contrgen}, in every dimension $n\ge2$ and for every $\bar\lambda$, for completeness and to fix conventions; it follows the classical derivation, in which the dimension only enters through the coefficient $\frac1{n-1}$ and $\Tr g=n$, and $\bar\lambda$ only through the term $(n-1)\bar\lambda\,dV_g$. Let $(\bar M^{n+1},\bar g)$ be Lorentzian with $\Ric(\bar g)=\bar\lambda\bar g$, $\Sigma$ a spacelike hypersurface with future unit normal $\nu$, $\bar g(\nu,\nu)=-1$, and data $(g,K)$, $K(X,Y)=\bar g(\bar\nabla_X\nu,Y)$; we write $\tau=\Tr_gK$ and $(K\circ K)_{ij}=K_i{}^mK_{mj}$. With the conventions of the proof of Theorem~\ref{thm:contrgen}, the Gauss equation reads $\bar R_{kilj}=R_{kilj}+K_{kl}K_{ij}-K_{kj}K_{il}$ for tangent indices, and, as in the proof of Lemma~\ref{lem:variationK}, $\bar\lambda g_{ij}=\Ric(\bar g)_{ij}=-\bar R_{\nu i\nu j}+g^{kl}\bar R_{kilj}$. Hence
\begin{equation}\label{eq:Rnunugen}
\bar R_{\nu i\nu j}=\big(\Ric_g-\bar\lambda g+\tau K-K\circ K\big)_{ij},
\end{equation}
and, since $g^{ij}\bar R_{\nu i\nu j}=\Ric(\bar g)(\nu,\nu)=-\bar\lambda$, the trace of \eqref{eq:Rnunugen} gives the Hamiltonian constraint
\begin{equation}\label{eq:contrH}
R_g-(n-1)\bar\lambda=|K|^2-\tau^2 .
\end{equation}

\begin{lemma}\label{lem:variationgen}
Let $Y$ be a vector field on $\bar M$ with $Y=\alpha\nu+X$ along $\Sigma$, where $\alpha\in C^\infty_c(\Sigma)$ and $X$ is a compactly supported vector field tangent to $\Sigma$. Let $\varphi_s$ be the flow of $Y$, and $(g_s,K_s)$ the data of $\varphi_s(\Sigma)$ pulled back to $\Sigma$ by $\varphi_s$. Then the first variation only depends on $(\alpha,X)$, and
\[
\dot g=2\alpha K+L_Xg,\qquad \dot K=\Hess\alpha-\alpha\big(\Ric_g-\bar\lambda g+\tau K-2K\circ K\big)+L_XK .
\]
\end{lemma}

\begin{proof}
The data of the immersion $f_s=\varphi_s|_\Sigma$ are smooth functions of its $2$-jet in the space variables; their first variation is therefore a linear differential operator in $\partial_sf_s|_{s=0}=Y|_\Sigma=\alpha\nu+X$, and we may treat the two parts separately and choose any family of immersions with this variation field. For the tangential part, we extend $X$ to a vector field tangent to $\Sigma$, whose flow preserves $\Sigma$: the first variation is $(L_Xg,L_XK)$. For the normal part, we take $f_s(x)=\exp_x(s\alpha(x)\nu(x))$, with variation field $W=\alpha\nu$, and argue as in the proof of Lemma~\ref{lem:variationK}, without assuming $K=0$. First, $\bar\nabla_{e_i}W=(\partial_i\alpha)\nu+\alpha K(e_i)$, where $K(e_i)=K_i{}^me_m$, so that $\dot g_{ij}=2\alpha K_{ij}$. Next, differentiating $\bar g(\nu_s,f_{s*}e_i)=0$ and $\bar g(\nu_s,\nu_s)=-1$ gives $\bar g(\bar\nabla_s\nu,e_i)=-\bar g(\nu,\bar\nabla_{e_i}W)=\partial_i\alpha$ and $\bar g(\bar\nabla_s\nu,\nu)=0$, that is $\bar\nabla_s\nu=\nabla\alpha$ at $s=0$. Finally, $K_s(e_i,e_j)=\bar g(\bar\nabla_{e_i}\nu_s,f_{s*}e_j)$, and
\[
\dot K_{ij}=\bar g\big(\bar\nabla_{e_i}\bar\nabla_s\nu+\bar R(W,e_i)\nu,e_j\big)+\bar g\big(\bar\nabla_{e_i}\nu,\bar\nabla_{e_j}W\big)=\nabla_i\nabla_j\alpha-\alpha\bar R_{\nu i\nu j}+\alpha(K\circ K)_{ij},
\]
since the tangential part of $\bar\nabla_{e_i}\nabla\alpha$ is $\nabla_{e_i}\nabla\alpha$ and $\bar g(\bar R(\nu,e_i)\nu,e_j)=-\bar R_{\nu i\nu j}$, as in Lemma~\ref{lem:variationK}. We conclude with \eqref{eq:Rnunugen}.
\end{proof}

\begin{proposition}\label{prop:hamilton}
In the setting of Section~\ref{ss:contrgen}, for all compactly supported $(\alpha,X)$, $\Gamma(\alpha,X)=J\big(\frac{\delta\mathcal H_{\alpha,X}}{\delta g},\frac{\delta\mathcal H_{\alpha,X}}{\delta\pi}\big)$ at $x_0$. Consequently \eqref{eq:hamilton} holds.
\end{proposition}

\begin{proof}
Write $\mathcal H=\mathcal H_{\alpha,X}$. By definition of $J$, $\omega\big(J(\tfrac{\delta\mathcal H}{\delta g},\tfrac{\delta\mathcal H}{\delta\pi}),y\big)=D\mathcal H(x_0)y$, and, $\alpha$ and $X$ being fixed, $D\mathcal H(x_0)y$ is the left-hand side of \eqref{eq:hamilton}; so the second statement follows from the first. Both sides of the first statement are linear in $(\alpha,X)$, and we treat the two parts separately. We write $\mu$ for the density $dV_g$, so that $\pi=\tilde\pi\,\mu$.

\emph{The tangential part.} Since $L^*=2\Div$, $\mathcal H_{0,X}(g,\pi)=2\int_M\ip{X}{\Div_g\tilde\pi}\,\mu=\int_M\ip{L_Xg}{\pi}$, where the last integrand is the natural contraction of the covariant tensor $L_Xg$ with the density $\pi$. As $L_Xg$ is linear in $g$,
\[
\begin{aligned}
D\mathcal H_{0,X}(x_0)(h,p)&=\int_M\big(\ip{L_Xh}{\pi_0}+\ip{L_Xg_0}{p}\big)=\int_M\big(\ip{L_Xg_0}{p}-\ip{h}{L_X\pi_0}\big)\\
&=\omega\big((L_Xg_0,L_X\pi_0),(h,p)\big),
\end{aligned}
\]
since $\ip{h}{\pi_0}$ is a density, whose Lie derivative integrates to zero. On the other hand, $(g,K)\mapsto(g,\pi)$ is natural, so that, by Lemma~\ref{lem:variationgen}, $\Gamma(0,X)=(L_Xg_0,L_X\pi_0)$. Hence
\[
\Gamma(0,X)=J\Big(\frac{\delta\mathcal H_{0,X}}{\delta g},\frac{\delta\mathcal H_{0,X}}{\delta\pi}\Big).
\]

\emph{The normal part.} We take $X=0$ and drop the index $0$ of $x_0$. In terms of $\tilde\pi=(K-\tau g)^{\sharp\sharp}$ we have $\Tr\tilde\pi=-(n-1)\tau$,
\[
\tilde\pi\circ\tilde\pi-\tfrac{1}{n-1}(\Tr\tilde\pi)\,\tilde\pi=\big(K\circ K-\tau K\big)^{\sharp\sharp},\qquad |\tilde\pi|^2-\tfrac{1}{n-1}(\Tr\tilde\pi)^2=|K|^2-\tau^2 .
\]
\emph{Variation in $\pi$.} The variation of $\int\alpha\big(|\tilde\pi|^2-\frac{1}{n-1}(\Tr\tilde\pi)^2\big)\mu$ in the direction $p$ is $\int2\alpha\ip{\tilde\pi^\flat-\frac{1}{n-1}(\Tr\tilde\pi)g}{p}$, so $\frac{\delta\mathcal H_{\alpha,0}}{\delta\pi}=2\alpha K$, which is the first component $\dot g$ of $\Gamma(\alpha,0)$ by Lemma~\ref{lem:variationgen}.

\emph{Variation in $g$.} In coordinates, $\big(|\tilde\pi|^2-\frac{1}{n-1}(\Tr\tilde\pi)^2\big)\mu=\mu^{-1}\big(g_{ik}g_{jl}\pi^{ij}\pi^{kl}-\frac1{n-1}(g_{ij}\pi^{ij})^2\big)$, and $\delta\mu=\frac12(\Tr h)\mu$. Together with $\delta(R\mu)=\big(DR(g)h+\frac12R\,\Tr h\big)\mu$, $\int\alpha\,DR(g)h\,\mu=\int\ip{DR(g)^*\alpha}{h}\mu$ and $DR(g)^*\alpha=\Hess\alpha+(\nabla^*\nabla\alpha)g-\alpha\Ric$, this gives
\[
-\frac{\delta\mathcal H_{\alpha,0}}{\delta g}=\Big(DR(g)^*\alpha-2\alpha\big(K\circ K-\tau K\big)+\tfrac12\alpha\big(R-(n-1)\bar\lambda+|K|^2-\tau^2\big)g\Big)^{\sharp\sharp}\mu .
\]
By \eqref{eq:contrH}, the last term equals $\alpha\big(R-(n-1)\bar\lambda\big)g$.

\emph{The variation of $\pi$.} Let $U=K-\tau g$, so that $\pi=U^{\sharp\sharp}\mu$. By Lemma~\ref{lem:variationgen}, $\dot g=2\alpha K$, $\dot\mu=\alpha\tau\mu$, $\dot g^{-1}=-2\alpha K^{\sharp\sharp}$, and, using $\Tr\Hess\alpha=-\nabla^*\nabla\alpha$ and $\Tr\Ric=R$,
\[
\dot\tau=-2\alpha|K|^2+\Tr\dot K=-\nabla^*\nabla\alpha-\alpha\big(R-n\bar\lambda+\tau^2\big).
\]
Hence $\dot U=\dot K-\dot\tau g-2\alpha\tau K$, and
\[
\dot\pi=\Big(\dot U-2\alpha\big(K\circ U+U\circ K\big)+\alpha\tau U\Big)^{\sharp\sharp}\mu .
\]
Substituting $\dot K$ from Lemma~\ref{lem:variationgen} and $K\circ U+U\circ K=2K\circ K-2\tau K$, the terms in $\tau^2g$ cancel and we get
\[
\dot\pi=\Big(DR(g)^*\alpha+\alpha\big(R-(n-1)\bar\lambda\big)g+2\alpha\tau K-2\alpha K\circ K\Big)^{\sharp\sharp}\mu=-\frac{\delta\mathcal H_{\alpha,0}}{\delta g}.
\]
Thus $\Gamma(\alpha,0)=\big(\frac{\delta\mathcal H_{\alpha,0}}{\delta\pi},-\frac{\delta\mathcal H_{\alpha,0}}{\delta g}\big)=J\big(\frac{\delta\mathcal H_{\alpha,0}}{\delta g},\frac{\delta\mathcal H_{\alpha,0}}{\delta\pi}\big)$.
\end{proof}

\section*{Acknowledgements}
This article was partially funded by the French National Research Agency ANR-23-
CE40-0010: Einstein constraints: past, present, and future (EINSTEIN-PPF) and
ANR-24-CE40-0702: Orbifold, scalar curvature and rigidity (OrbiScaR).

\section*{AI usage statement}
The AI assistant Claude (Anthropic)  was used for the computations and for writing the English version of this article. The author takes full responsibility for its content.

\end{document}